\input ./style/arxiv-general.cfg
\documentclass[aap,MSNbibl,seceqn,dvips]{arximspdf}
\makeatletter
   \@ifpackageloaded{graphicx}{}{\usepackage{graphicx}}
\makeatother

\doi{10.1214/15-AAP1135}
\volume{26}
\issue{3}
\pubyear{2016}
\firstpage{1837}
\lastpage{1887}
\docsubty{FLA}

\makeatletter

\newcommand{\rrvert}{\vert}

\newcommand{\llvert}{\vert}
\def\xrightarrow{\longrightarrow}
\newcommand{\eqref}[1]{(\ref{#1})}
\newtheorem{theorem}{Theorem}
\newtheorem{prop}{Proposition}[section]
\newtheorem{lemma}[prop]{Lemma}
\newtheorem{coro}[prop]{Corollary}
\newproclaim{assump}{Assumption}
\newproclaim{rem}{Remark}[section]
\makeatother

\begin{document}
\begin{frontmatter}

\title{Gaussian fluctuations for linear spectral statistics of large
random covariance matrices}
\runtitle{Gaussian fluctuations for large covariance matrices}

\begin{aug}
\author[A]{\fnms{Jamal}~\snm{Najim}\corref{}\thanksref{t1}\ead[label=e1]{najim@univ-mlv.fr}}
\and
\author[B]{\fnms{Jianfeng}~\snm{Yao}\thanksref{t2}\ead[label=e2]{jianfeng.yao.sh@gmail.com}}
\runauthor{J. Najim and J. Yao}
\affiliation{Universit\'{e} Paris Est and CNRS and T\'{e}l\'{e}com Paristech}
\address[A]{Institut Gaspard Monge LabInfo\\
UMR 8049\\
Universit\'e Paris Est Marne-la-Vall\'ee\\
5 Boulevard Descartes\\
Champs sur Marne\\
77454 Marne-la-Vall\'ee Cedex 2\\
France\\
\printead{e1}}
\address[B]{T\'{e}l\'{e}com Paristech\\
46 rue Barrault\\
75634 Paris cedex\\
France\\
\printead{e2}}
\end{aug}
\thankstext{t1}{Supported by ANR Grant ANR-12-MONU-003 and Labex B\'{e}zout.}
\thankstext{t2}{Supported by program ``Futur et ruptures'' of Fondation
T\'{e}l\'{e}com.}

%
\received{\smonth{10} \syear{2014}}
%
\revised{\smonth{7} \syear{2015}}

%
\begin{abstract}
Consider a $N\times n$ matrix
$\Sigma_n= \frac{1}{\sqrt{n}} R_n^{1/2} X_n$, where $R_n$ is a
nonnegative definite Hermitian matrix and $X_n$ is a random
matrix with i.i.d. real or complex standardized entries. The
fluctuations of the linear statistics of the eigenvalues
\[
\operatorname{Trace} f\bigl(\Sigma_n \Sigma_n^*\bigr) =
\sum_{i=1}^N f(\lambda_i),
\qquad (\lambda_i)\mbox{ eigenvalues of } \Sigma_n
\Sigma_n^*,
\]
are shown to be Gaussian, in the regime where both dimensions of matrix
$\Sigma_n$ go to infinity at the same pace and in the case where $f$
is of class $C^3$, that is, has three continuous derivatives. The main
improvements with respect to Bai and Silverstein's CLT
[\textit{Ann. Probab.} \textbf{32} (2004) 553--605]
are twofold: First, we consider general
entries with finite fourth moment, but whose fourth cumulant is
nonnull, that is, whose fourth moment may differ from the moment of a
(real or complex) Gaussian random variable. As a consequence, extra
terms proportional to
\[
\vert \mathcal{V}\vert ^2= \bigl|\mathbb{E}\bigl(X_{11}^n
\bigr)^2 \bigr|^2\quad\mbox{and} \quad\kappa= \mathbb{E}\bigl
\llvert X_{11}^n\bigr\rrvert^4 -\llvert{
\mathcal V}\rrvert^2 - 2
\]
appear in the limiting variance and in the limiting bias, which not
only depend on the spectrum of matrix $R_n$ but also on its eigenvectors.
Second, we relax the analyticity assumption over $f$ by representing
the linear statistics with the help of
Helffer--Sj\"{o}strand's formula.

The CLT is expressed in terms of vanishing L\'{e}vy--Prohorov distance
between the linear statistics' distribution and a Gaussian probability
distribution, the mean and the variance of which depend upon $N$ and
$n$ and may not converge.
\end{abstract}

%
\end{frontmatter}

\maxcsec{10}
\tableofcontents

\section{Introduction}
Empirical random covariance matrices, whose probabilistic study may be
traced back to Wishart \cite{wishart-1928} in the late twenties, play
an important role in applied mathematics. After Mar\v{c}enko and
Pastur's seminal contribution \cite{mar-and-pas-67} in 1967, the large
dimensional setting (where the dimension of the observations is of the
same order as the size of the sample) has drawn a growing interest, and
important theoretical contributions \cite{bai-silverstein-aop-2004,silverstein-95,johnstone-2001} found many
applications in multivariate statistics, electrical engineering,
mathematical finance, etc.; cf. \cite{book-bai-chen-liang-2009,book-couillet-debbah,bouchaud-noise-dressing-99,guhr-credit-risk-2014}.
The aim of this paper is to describe the fluctuations for linear
spectral statistics of large empirical random covariance matrices. It
will complete the picture already provided by Bai and Silverstein \cite
{bai-silverstein-aop-2004} and will hopefully provide a generic result
of interest
for practitioners.

\subsection*{The model}
Consider a $N\times n$ random matrix $\Sigma
_n=(\xi_{ij}^{n})$ given by
%
\begin{equation}
\label{eqmodel-sigma}
\Sigma_n= \frac{1}{\sqrt{n}} R_n^{1/2}
{X_n},
\end{equation}
where $N=N(n)$ and $R_n$ is a $N\times N$ nonnegative definite
Hermitian matrix with spectral norm uniformly bounded in $N$. The
entries $(X_{ij}^{n}; i\le N,j\le n, n\ge1)$ of
matrices $(X_n)$ are real or complex, independent and
identically distributed (i.i.d.) with mean 0 and variance 1.
Matrix $\Sigma_n \Sigma_n^*$ models a sample covariance matrix,
formed from $n$ samples of the random vector $R_n^{1/2} X_{\cdot1}^n$,
with the population covariance matrix $R_n$. In the asymptotic regime where
%
\begin{equation}
\label{eqasymptotic}
N,n\rightarrow\infty\quad\mbox{and} \quad 0
<\liminf
\frac{N}n \le\limsup\frac{N}n < \infty
\end{equation}
(a condition that will be simply referred as $N,n\rightarrow\infty$
in the sequel), we study
the fluctuations of linear spectral statistics of the form:
%
\begin{equation}
\label{eqlinear-spectral-statistics}
\operatorname{tr}f\bigl(\Sigma_n
\Sigma_n^*
\bigr)= \sum_{i=1}^N f(\lambda
_i)\qquad\mbox{as } N,n\rightarrow\infty,
\end{equation}
where $\operatorname{tr}(A)$ refers to the trace of $A$ and the
$\lambda_i$'s are
the eigenvalues of $\Sigma_n \Sigma_n^*$. This subject has a rich
history with contributions by Arharov \cite{arharov}, Girko (see \cite
{girko-canonical-equations-I,girko-canonical-equation-II} and the
references therein), Jonsson \cite{jonsson-1982}, Khorunzhiy et al.
\cite{pastur-et-al-1996}, Johansson \cite{johansson-1998}, Sinai and
Soshnikov \cite{sinai-soshnikov-1998,soshnikov-2000},
Cabanal-Duvillard \cite{cabanal-duvillard-2001}, Guionnet \cite
{guionnet-2002-ihp}, Bai and Silverstein \cite
{bai-silverstein-aop-2004}, Anderson and Zeitouni \cite
{anderson-zeitouni-ptrf}, Pan et al. \cite
{pan-2014-IHP,pan-zhou-2008}, Chatterjee \cite{chatterjee-ptrf-2009},
Lytova and Pastur \cite{lytova-pastur-aop-2009}, Bai et al. \cite
{bai-wang-2010}, Shcherbina \cite{shcherbina-2011}, etc. There are
also contributions for heavy-tailed entries (see, e.g.,
Benaych-Georges et al. \cite{2014-benaych}).

In their 2004 article \cite{bai-silverstein-aop-2004}, Bai and
Silverstein established a CLT for the linear spectral statistics \eqref
{eqlinear-spectral-statistics}
as the dimensions $N$ and $n$ grow to infinity at the same pace
[$N/n\rightarrow c\in(0,\infty)$]  and under two important assumptions:
\begin{longlist}[2.]
\item[1.]  The entries $(X_{ij}^n)$ are centered with unit
variance and a finite fourth moment equal to the fourth moment of a
(real or complex) Gaussian standard variable.
\item[2.]  Function $f$ in \eqref{eqlinear-spectral-statistics} is analytic in a neighborhood of the
asymptotic spectrum of $\Sigma_n \Sigma_n^*$.
\end{longlist}
Such a result proved to be highly useful in probability theory,
statistics and various other fields.

The purpose of this article is to establish a CLT for linear spectral
statistics \eqref{eqlinear-spectral-statistics} for general entries
$X_{ij}^n$ with finite fourth moment and
for nonanalytic functions $f$, sufficiently regular, hence to relax
both assumptions (1) and (2) in \cite
{bai-silverstein-aop-2004}.

It is well known since the paper by Khorunzhiy et al. \cite
{pastur-et-al-1996} that if the fourth moment of the entries differs
from the fourth moment of a Gaussian random variable, then a term
appears in the variance of the trace of the resolvent, which is
proportional to the fourth cumulant of the entries. This term does not
appear if assumption (1) holds true because, in this case,
the fourth cumulant is zero.

In Pan and Zhou \cite{pan-zhou-2008}, assumption (1) has
been relaxed under an additional assumption on matrix $R_n$, which
somehow enforces structural conditions on $R_n$ (in particular, these
conditions are satisfied if matrix $R_n$ is diagonal).
In Hachem et al. \cite{kammoun-et-al-CLT-2008,hachem-et-al-clt-2012},
CLTs have been established for specific linear statistics of interest
in information theory, with general entries and (possibly noncentered)
covariance random matrices with a variance profile. In Bao et al. \cite
{pan-CLT-MIMO-preprint-2013}, the CLT is established for the white
model (where $R_n$ is equal to the identity matrix) with general
entries with finite fourth moment, featuring terms in the covariance
proportional to the square of the second nonabsolute moment and to the
fourth cumulant.

In Lytova and Pastur \cite{lytova-pastur-aop-2009} and Shcherbina
\cite{shcherbina-2011}, both assumptions have been relaxed for the
white model. In \cite{lytova-pastur-aop-2009}, it has been proved that
mild integrability conditions over the Fourier transform of $f$ was
enough to establish the CLT. In Bai et al. \cite{bai-wang-2010},
fluctuations for the white model are addressed as well, for
nonanalytic functions $f$. Following Shcherbina's ideas, Gu\'{e}don et
al. \cite{2014-guedon-et-al} establish a CLT for linear statistics of
large covariance matrices with vectors with log-concave distribution.
Following Lytova and Pastur, Yao \cite{yao-PhD-2013} relaxes the
analyticity assumption in \cite{bai-silverstein-aop-2004} by using
interpolation techniques and Fourier transforms. We follow here a
different approach, inspired from Bordenave \cite
{bordenave-personal-comm-2013}.

%
%

\subsection*{Non-Gaussian entries} The presence of matrix $R_n$ yields
interesting phenomena at the CLT level when considering entries with
non-Gaussian fourth moment: terms proportional to the fourth cumulant
and to $|\mathbb{E}(X_{11}^n)^2|^2$ appear in the asymptotic
variance (described in Section~\ref{secvariance-breakdown}); however,
their convergence is not granted under usual assumptions (roughly,
under the convergence of $R_n$'s spectrum), mainly because these extra
terms also depend on the eigenvectors of $R_n$.
As a consequence, such terms may not converge unless some very strong
structural assumption over $R_n$ (such as $R_n$ diagonal) is made.
This lack of convergence has consequences on the description of the
fluctuations.


Denote by $L_n(f)$ the (approximately) centered version of the linear
statistics~\eqref{eqlinear-spectral-statistics}, to be properly
defined below.
Instead of expressing the CLT in the usual way, that is
($\displaystyle\mathop{\rightarrow}^{\mathcal{D}}$ stands for the
convergence in distribution)
%
\begin{equation}
\label{equsual-way}
L_n(f) \mathop{\xrightarrow}_{N,n\rightarrow\infty
}^{\mathcal D}
{\mathcal N}\bigl({\mathcal B}^f_\infty,
\Theta^f_\infty\bigr),
\end{equation}
for some well-defined parameters ${\mathcal B}^f_\infty,\Theta
^f_\infty$,
we prove that the distribution of the linear statistics $L_n(f)$
becomes close to a family of Gaussian distributions, whose parameters
(mean and variance) may not converge. More precisely,\vspace*{1pt} we establish that
there exists a family of Gaussian random variables ${\mathcal
N}({\mathcal B^f_n}, \Theta^f_n)$, such that
%
\begin{equation}
\label{eqconvergence-LP}
d_{\mathcal{LP}} \bigl( L_n(f), {\mathcal
N}\bigl({
\mathcal B^f_n}, \Theta^f_n\bigr)
\bigr) \mathop{\xrightarrow}_{N,n\rightarrow\infty} 0,
\end{equation}
where $d_{\mathcal{LP}}$ denotes the L\'{e}vy--Prohorov distance (and in
particular metrizes the convergence of laws). Details are provided in
Section~\ref{secconvergence-preliminary} and the fluctuation results
are stated in Theorem~\ref{lemmamain} [for the resolvent $f(\lambda
)=(\lambda-z)^{-1}$] and Theorem~\ref{thnon-analytic} (for $f$ of
class $C^3$, the space of functions with third continuous derivative).

From a technical point of view, the analysis of the extra term
proportionnal to the fourth cumulant requires to cope with quadratic
forms of the resolvent (counterpart of isotropic Mar\v{c}enko--Pastur
law). We provide the needed results in Section~\ref{secproof}.



Expressing the CLT as in \eqref{eqconvergence-LP} makes it possible
to avoid any cumbersome assumption related to the joint convergence of
$R_n$'s eigenvectors and
eigenvalues; the technical price to pay however is the need to get
various uniform (in $N,n$) controls over the sequence ${\mathcal
N}({\mathcal B}_n,\Theta_n)$. This is
achieved by introducing a matrix meta-model in Section~\ref
{secmeta-model}. The case where matrix $R_n$ is diagonal is simpler
and the fluctuations express in the usual
way \eqref{equsual-way}; it is handled in Section~\ref{secdiagonal-R}.
Remarks on the white case ($R_n=I_N$) are also
provided in Sections~\ref{secwhite-case} and
\ref{secexplicit-formulas}.

This framework may also prove to be useful for other interesting models
such as large dimensional information-plus-noise type matrices \cite
{dozier-silverstein-07-a,hachem-et-al-2007} and more generally mixed
models combining large dimensional deterministic and random matrices.

\setcounter{footnote}{2}
\subsection*{Nonanalytic functions} In Section~\ref{sectionmain},
we establish the CLT for the trace of the resolvent
\[
\operatorname{tr} \bigl( \Sigma_n \Sigma_n^* -zI_N
\bigr)^{-1}.
\]
In order to transfer the CLT from the resolvent to the linear
statistics of the eigenvalues $ \operatorname{tr} f(\Sigma_n \Sigma
_n^*)$, we will use (Dynkin--)Helffer--Sj\"{o}strand's representation
formula\footnote{In \cite{helffer-book-2013}, Notes of Chapter~8, it
is written ``\textit{This formula is due to Dynkin but was popularized by
Helffer and Sj\"{o}strand in the context of spectral theory, leading
many authors to call it the Helffer--Sj\"{o}strand formula}.''} for a
function $f$ of class $C^{k+1}$ and with compact support \cite
{dynkin-1972,helffer-sjostrand-1989}. Denote by $\Phi_k(f):\mathbb
{C}^+\to\mathbb{C}$ the function
%
\begin{equation}
\label{eqalmost-analytic}
\Phi_k(f) (x+\mathbf{i}y) =\sum
_{\ell=0}^k \frac{(\mathbf
{i}y)^\ell}{\ell!} f^{(\ell)}(x)
\chi(y),
\end{equation}
where $\chi:\mathbb{R}\to\mathbb{R}^+$ is smooth, compactly
supported, with value 1
in a neighborhood of 0. Function $\Phi_k(f)$ coincides with $f$ on
the real line and is an appropriate extension of $f$ to the complex plane.
Let $\overline\partial= \partial_x +\mathbf{i}\partial_y$, then
Helffer--Sj\"{o}strand's formula writes
%
\begin{equation}
\label{eqhelffer-sjostrand} \operatorname{tr} f\bigl(\Sigma_n \Sigma_n^*
\bigr) = \frac{1}\pi\operatorname{Re} \int_{\mathbb{C}^+} \overline
\partial\Phi_k(f) (z) \operatorname{tr}\bigl(\Sigma_n
\Sigma_n^*-zI_N\bigr)^{-1}
\ell_2(dz),
\end{equation}
where $\ell_2$ stands for the Lebesgue measure over $\mathbb{C}^+$.
An elementary proof of formula \eqref{eqhelffer-sjostrand} can be found
in \cite{bordenave-lecture-notes}, Chapter~5.
Closest to our work are the papers by Pizzo, O'Rourke, Renfrew and
Soshnikov \cite{2013-orourke-et-al,soshnikov-et-al-2012} where the
fluctuations of the entries of regular functions of Wigner and large
covariance matrices are studied; see also the paper by Bao et al. \cite
{bao-pan-zhou-2013} where a CLT for partial linear eigenvalue
statistics is established for Wigner matrices. 

We believe that formula \eqref{eqhelffer-sjostrand} provides a very
streamlined way to handle nonanalytic functions and in fact enables us
to state the fluctuations for the linear statistics for functions of
class $C^3$, a lower regularity requirement than in \cite
{bai-wang-2010,lytova-pastur-aop-2009,yao-PhD-2013}; in Shcherbina's
article \cite{shcherbina-2011}, the requirements over the functions
are lower and expressed in terms of Sobolev norms $\|f\|
_{3/2+\varepsilon}<\infty$,\vspace*{1pt} a condition that is fulfilled if $f$ is
$C^2$ (with bounded derivatives in $L^2$).

\subsection*{Bias in the CLT and asymptotic expansion for the linear
spectral statistics}
Beside the fluctuations, a substantial part of this article is devoted
to the study of the bias that we describe hereafter.
In order to center the linear spectral statistics $\operatorname
{tr}f(\Sigma_n
\Sigma_n^*)$, we consider the (first-order) expansion of $ \frac{1}N
\mathbb{E} \operatorname{tr}f(\Sigma_n \Sigma_n^*)$
\[
\frac{1}N \mathbb{E} \operatorname{tr}f\bigl(\Sigma_n
\Sigma_n^*\bigr)= {\mathcal E}_{0,n}(f) + {\mathcal O}
\biggl( \frac{1}N \biggr),
\]
where ${\mathcal E}_{0,n}(f)$ is ${\mathcal O}(1)$ and does not depend
on the distribution of the entries of~$X_n$, and define $L_n(f)$ as
\[
L_n(f)= \operatorname{tr}f\bigl(\Sigma_n
\Sigma_n^*\bigr) -N {\mathcal E}_{0,n}(f).
\]
A precise description of $L_n(f)$ is provided in Section~\ref
{secrepresentation-LS}. In order to fully characterize the
fluctuations of $L_n(f)$, we must study the second-order expansion of $
\frac{1}N \mathbb{E} \operatorname{tr}f(\Sigma_n \Sigma_n^*)$,
\[
\frac{1}N \mathbb{E} \operatorname{tr}f\bigl(\Sigma_n
\Sigma_n^*\bigr)= {\mathcal E}_{0,n}(f) +\frac{{\mathcal
E}_{1,n}(f)}{N}+ o
\biggl( \frac{1}N \biggr),
\]
which will naturally yield the bias of $L_n(f)$, as $\mathbb{E}
L_n(f)={\mathcal E}_{1,n}(f) + o(1)$. Asymptotic expansions for various
matrix ensembles have already been studied; see, for instance, Pastur
et al. \cite{pastur-2001}, Bai and Silverstein \cite
{bai-silverstein-aop-2004}, Haagerup and Thorbj\o rnsen \cite
{haagerup-thorbjornsen-2005,haagerup-2012}, Schultz \cite
{schultz-ptrf-2005}, Capitaine and Donati-Martin \cite
{capitaine-donati-freeness-2007},
Vallet et al. \cite{vallet-2012}, Hachem et al. \cite
{hachem-et-al-aihp-2013}, etc.

The asymptotic bias is expressed in Theorem~\ref{lemmamain} for the resolvent.
In order to lift asymptotic expansions from the resolvent to smooth
functions, we combine ideas from Haagerup and Thorbj\o rnsen \cite
{haagerup-thorbjornsen-2005} and Loubaton et al. \cite
{hachem-et-al-aihp-2013,vallet-2012} together with some Gaussian
interpolation and the use of Helffer--Sj\"{o}strand's formula. For
smooth functions, the statement is given in Theorem~\ref
{thnon-analytic-bias}. Somehow surprisingly, the condition over
function $f$ is stronger for the asymptotic\vspace*{1pt} expansion to hold than for
the CLT as function $f$ needs to be of class $C^{18}$ (cf. Remark~\ref
{rem18}).

\section{General background---variance and bias formulas} \label
{secbackground}

\subsection{Assumptions} Recall the asymptotic regime where
$N,n\rightarrow\infty$, cf. \eqref{eqasymptotic}, and
denote by
\[
c_n=\frac{N}n,\qquad\bolds{\ell}^- = \liminf
\frac{N}n\quad\mbox{and}\quad\bolds{\ell}^+ = \limsup
\frac{N}n.
\]
\renewcommand{\theassump}{A-1}
\begin{assump}
\label{assX}
The random\vspace*{1pt} variables $(X_{ij}^n ; 1\le i\le N(n), 1\le j\le n,
n\ge1)$
are independent and identically distributed. They satisfy
\[
\mathbb{E}X_{ij}^n = 0,\qquad\mathbb{E}\bigl|X_{ij}^n\bigr|^2=1
\quad\mbox{and} \quad\mathbb{E}\bigl|X_{ij}^n\bigr|^4<
\infty.
\]
%
\end{assump}
\renewcommand{\theassump}{A-2}
\begin{assump}
\label{assR}
Consider a sequence $(R_n)$ of deterministic,
nonnegative definite Hermitian $N\times N$ matrices, with $N=N(n)$.
The sequence $(R_n,n\ge1)$
is bounded for the spectral norm as $N,n\rightarrow\infty$:
\[
\sup_{n\ge1} \|R_n\| <\infty.
\]
\end{assump}

In particular, we will have
\[
0 \le\bolds{\lambda}_R^- \stackrel{\triangle} {=} \liminf
_{N,n\rightarrow\infty} \| R_n\| \le\bolds{\lambda
}_R^+ \stackrel{\triangle} {=} \limsup_{N,n\rightarrow\infty} \|
R_n\| < \infty.
\]

\subsection{Resolvent, canonical equation and deterministic
equivalents}
Denote by $Q_n(z)$ (resp., $\widetilde{Q}_n$) the
resolvent of
matrix $\Sigma_n \Sigma_n^*$ (resp., of $\Sigma_n^*\Sigma_n$):
%
\begin{equation}
\label{eqresolvent}
Q_n(z)= \bigl( \Sigma_n
\Sigma_n^*- zI_N \bigr)^{-1},\qquad \widetilde{Q}_n(z)= \bigl( \Sigma^*_n \Sigma_n-
zI_n \bigr)^{-1},
\end{equation}
and by $f_n(z)$ and $\tilde f_n(z)$ their normalized traces which are
the Stieltjes transforms of the empirical distribution of $\Sigma_n
\Sigma_n^*$'s and $\Sigma_n^*\Sigma_n$'s eigenvalues:
%
\begin{equation}
\label{eqST}
f_n(z)= \frac{1}N\operatorname{tr}Q_n(z),\qquad \tilde f_n(z)= \frac{1}n\operatorname{tr}\widetilde
Q_n(z).
\end{equation}
The following canonical equation\footnote{We borrow the name
``canonical equation'' from V. L. Girko who established in \cite
{girko-canonical-equations-I,girko-canonical-equation-II} canonical
equations associated to various models of large random matrices.}
admits a unique solution $t_n$ in the class of Stieltjes transforms of
probability measures (see, e.g., \cite{bai-silverstein-aop-2004}):
%
\begin{equation}
\label{eqdeteq-ST}\hspace*{6pt}
t_n(z)= \frac{1}N \operatorname{tr} \bigl(
-zI_N + (1-c_n)R_n -zc_n
t_n(z) R_n \bigr)^{-1},\qquad  z\in\mathbb{C}\setminus\mathbb{R}^+.
\end{equation}
The function $t_n$ being introduced, we can define the following
$N\times N$ matrix:
%
\begin{equation}
\label{eqdeteq-RESOLVENT}
T_n(z)= \bigl( -zI_N +
(1-c_n)R_n -zc_n t_n(z)
R_n \bigr)^{-1}.
\end{equation}
Matrix $T_n(z)$ can be thought of as a \textit{deterministic equivalent}
of the resolvent $Q_n(z)$ in the sense that it
approximates the resolvent in various senses. For instance,
\[
\frac{1}N \operatorname{tr}T_n(z) - \frac{1}N
\operatorname{tr}Q_n(z)\mathop{\xrightarrow}_{N,n\rightarrow\infty} 0,\qquad
z\in  \mathbb{C}\setminus\mathbb{R}^+
\]
(in probability or almost surely). Otherwise stated, $t_n(z)=N^{-1}
\operatorname{tr}
T_n(z)$ is the deterministic equivalent of $f_n(z)$.
As we shall see later in this paper, the following property holds true:
%
\begin{equation}
\label{eqconv-quadra} u_n^* Q_n(z) v_n
-u_n^* T_n(z) v_n \mathop{
\xrightarrow}_{N,n\rightarrow\infty} 0,
\end{equation}
where $(u_n)$ and $(v_n)$ are deterministic $N\times1$ vectors with
uniformly bounded Euclidean norms in $N$. As a consequence of \eqref
{eqconv-quadra}, not only $T_n$ conveys information on the limiting
spectrum of the resolvent $Q_n$ but also on the eigenvectors of $Q_n$.

If $R_n=I_N$, then $t_n$ is simply the Stieltjes transform of Mar\v{c}enko--Pastur's distribution \cite{mar-and-pas-67} with parameter $c_n$.

\subsection{Entries with nonnull fourth cumulant and the limiting
covariance for the trace of the resolvent} \label{secvariance-breakdown}
As in \cite{bai-silverstein-aop-2004}, we first study the CLT for the
trace of the resolvent. Let ${\mathcal V}$ be the second moment of the random
variable $X_{ij}$ and $\kappa$ its fourth cumulant:
\[
{\mathcal V} = \mathbb{E}\bigl( X^n_{ij}
\bigr)^2 \quad\mbox{and}\quad \kappa= \mathbb{E}\bigl\llvert
X_{ij}^n\bigr\rrvert^4 -\llvert{\mathcal V}
\rrvert^2 - 2.
\]
If the entries are real or complex standard Gaussian, then ${\mathcal
V}=1$ or 0 and $\kappa=0$. Otherwise the fourth cumulant
is a priori no longer equal to zero. This induces extra terms in the
computation of the limiting variance, mainly due to the following
$(\mathcal{V},\kappa)$-dependent\vspace*{0pt} identity:
%
\begin{eqnarray}
&& \mathbb{E}\bigl(X^*_{\cdot1} A X_{\cdot1} -
\operatorname{tr}A\bigr) \bigl(X^*_{\cdot1} B X_{\cdot1} -
\operatorname{tr}B\bigr)
\nonumber
\\[-10pt]
\label{eqcovariance-identity}
\\[-10pt]
\nonumber
&&\qquad= \operatorname{tr}AB + \llvert\mathcal
{V}\rrvert
^2 \operatorname{tr}A B^T + \kappa\sum
_{i=1}^N A_{ii} B_{ii},
\end{eqnarray}
where $X_{\cdot1}$ stands for the first column (of dimension $N\times
1$) of matrix $X_n$ and where $A,B$ are deterministic $N\times N$ matrices.
As a consequence, there will be three terms in the limiting covariance
of the quantity \eqref{eqlinear-spectral-statistics}; one will raise
from the first term of the right-hand side (RHS) of \eqref
{eqcovariance-identity}, a second one will be proportional to
$|{\mathcal V}|^2$, and a third one to $\kappa$.
In order to describe these terms, let
%
\begin{eqnarray}
\label{defttilde}
\tilde{t}_n(z) &=& -\frac{1-c_n}z + c_n
t_n(z).
\end{eqnarray}
The\vspace*{1pt} quantity $\tilde{t}_n(z)$ is the deterministic equivalent associated
to $n^{-1} \operatorname{tr}(\Sigma_n^* \Sigma_n - zI_n)^{-1}$.
Denote by $R^T_n$
the transpose matrix of $R_n$ (notice that since $R_n$ is Hermitian,
$R_n^T=\overline{R}_n$ and we shall use this latter notation) and by
$T^T_n$, the transpose matrix\footnote{Beware that $T^T_n$ is not the
entry-wise conjugate of $T_n$, due to the presence of $z$.} of $T_n$:
%
\begin{equation}
\label{eqdeteq-RESOLVENT-conjugate}
T^T_n(z)= \bigl( -zI_N +
(1-c_n)\overline{R}_n -zc_n t_n(z)
\overline{R}_n \bigr)^{-1};
\end{equation}
notice that the definition of $t_n(z)$ in \eqref{eqdeteq-ST} does not
change if $R_n$ is replaced by $\overline{R}_n$ since the spectrum of both
matrices $R_n$ and $\overline{R}_n$ is the same. We can now describe the
limiting covariance of the trace of the resolvent
%
\begin{eqnarray}
&&  \operatorname{cov} \bigl( \operatorname{tr}Q_n(z_1),
\operatorname{tr}Q_n(z_2) \bigr)\nonumber\\
\label{eqdecomposition-covariance}
&&\qquad= \Theta
_{0,n}(z_1,z_2) + |{\mathcal V}|^2
\Theta_{1,n}(z_1,z_2) + \kappa
\Theta_{2,n}(z_1,z_2) + o(1)
\\
\nonumber
&&\qquad\stackrel{\triangle}{=} \Theta_n(z_1,z_2)
+ o(1),
\end{eqnarray}
where $o(1)$ is a term that converges to zero as $N,n\rightarrow\infty
$ and
%
\begin{eqnarray}
\label{deftheta0}
\Theta_{0,n}(z_1,z_2) &\stackrel{\triangle}
{=}& \biggl\{ \frac
{\tilde t_n'(z_1)
\tilde t'_n(z_2)}{(\tilde t_n(z_1) - \tilde t_n(z_2))^2} - \frac{1}{(z_1
- z_2)^2} \biggr\},
\\
\label{deftheta1}
\Theta_{1,n}(z_1,z_2) &\stackrel{\triangle}
{=}& \frac{\partial}{
\partial z_2} \biggl\{ \frac{\partial{\mathcal A}_n(z_1,z_2)} {\partial
z_1} \frac{1}{1-|\mathcal{V}
|^2{\mathcal A}_n(z_1,z_2)} \biggr\},
\\
\Theta_{2,n}(z_1,z_2) &\stackrel{\triangle}
{=}& \frac{ z_1^2 z_2^2
\tilde t'_n(z_1)
\tilde t'_n(z_2)}n
\nonumber
\\[-8pt]
\label{deftheta2}
\\[-8pt]
\nonumber
&&{}\times
 \sum_{i=1}^N \bigl(
R_n^{1/2} T^2_n(z_1)
R_n^{1/2} \bigr)_{ii} \bigl( R_n^{1/2}
T^2_n(z_2) R_n^{1/2}
\bigr)_{ii},
\end{eqnarray}
with
%
\begin{equation}
\label{defAn}
\hspace*{5pt}\quad{\mathcal A}_n(z_1,z_2) =
\frac{z_1 z_2} n \tilde t_n(z_1) \tilde
t_n(z_2) \operatorname{tr} \bigl\{ R_n^{1/2}T_n(z_1)R_n^{1/2}
\overline{R}_n^{1/2}T^T_n(z_2)
\overline{R}_n^{1/2} \bigr\}.
\end{equation}
For alternative formulas for $\Theta_{0,n}$ and $\Theta_{2,n}$, see
Remarks \ref{remaltern-theta0} and \ref{remaltern-theta2}.

At first sight, these formulas (established in Section~\ref{secproof})
may seem complicated; however, much information can be
inferred from them.

\subsubsection*{The term $\Theta_{0,n}$} This term is familiar as it
already appears in Bai and Silverstein's CLT \cite
{bai-silverstein-aop-2004}. Notice that the quantities $\tilde t_n$ and
$\tilde t_n'$ only depend on the spectrum of matrix $R_n$. Hence, under
the additional assumption that
%
\begin{equation}
\label{eqadd-assumption} c_n \mathop{\xrightarrow}_{N,n\rightarrow\infty
} c\in(0,
\infty) \quad\mbox{and}\quad F^{R_n} \mathop{\xrightarrow}_{N,n\rightarrow\infty
}^{\mathcal D}
F^{\mathbf R},
\end{equation}
where $F^{R_n}$ denotes the empirical distribution of $R_n$'s
eigenvalues and $F^{\mathbf R}$ is a probability measure, it can easily be
proved that
%
\begin{equation}
\label{eqdef-theta0}
\hspace*{12pt}\Theta_{0,n}(z_1,z_2) \mathop{
\xrightarrow}_{N,n\rightarrow\infty} \Theta_0(z_1,z_2)
= \biggl\{ \frac{\tilde t'(z_1) \tilde t'(z_2)}{(\tilde
t(z_1) - \tilde t(z_2))^2} - \frac{1}{(z_1 - z_2)^2} \biggr\},
\end{equation}
where $\tilde t,\tilde t'$ are the limits of $\tilde t_n,\tilde t_n'$
under \eqref{eqadd-assumption}.

\subsubsection*{The term $\Theta_{1,n}$} The interesting phenomenon
lies in the fact that this term involves products of matrices
$R_n^{1/2}$ and its conjugate $\overline{R}_n^{1/2}$. These matrices have
the same spectrum but conjugate eigenvectors. If $R_n$ is not real, the
convergence of $\Theta_{1,n}$ is not granted, even under \eqref
{eqadd-assumption}. If however $R_n$ and $X_n$'s entries are real,
that is, ${\mathcal V}=1$, then it can be easily proved that $\Theta
_{0,n}=\Theta_{1,n}$ hence the factor 2 in \cite
{bai-silverstein-aop-2004} between the complex and the real covariance.

\subsubsection*{The term $\Theta_{2,n}$} This term involves
quantities of the type $(R_n^{1/2} T_n R_n^{1/2})_{ii}$ which not only
depend on the spectrum of matrix $R_n$ but also on its eigenvectors. As
a consequence, the convergence of such terms does not follow from an
assumption such as \eqref{eqadd-assumption}, except in some
particular cases (e.g., if $R_n$ is diagonal) and any assumption
which enforces the convergence of such terms (as, e.g., in \cite{pan-zhou-2008}, Theorem~1.4) implicitly implies an asymptotic joint
behavior between $R_n$'s eigenvectors and eigenvalues. We shall adopt
a different point of view here and will not assume the convergence of
these quantities.


\subsection{Representation of the linear statistics and limiting
bias}\label{secrepresentation-LS} Recall that $t_n(z)$ is the
Stieltjes transform of a probability measure ${\mathcal F}_n$:
%
\begin{equation}
\label{eqstieltjes-basics}
t_n(z) = \int_{{\mathcal S}_n}
\frac{ {\mathcal F}_n(d\lambda
)}{\lambda-z}
\end{equation}
with support ${\mathcal S}_n$ included in a compact set. The purpose of
this article is to describe the fluctuations of the linear statistics
%
\begin{equation}
\label{eqlinear-stat-centered} L_n(f)=\sum_{i=1}^N
f(\lambda_i) - N\int f(\lambda) {\mathcal F}_n(d\lambda)
\end{equation}
as $N,n\rightarrow\infty$.

For a smooth enough function $f$ of class $C^{k+1}$ with bounded
support, one can rely on Helffer--Sj\"{o}strand's formula and write
\begin{eqnarray}
L_n(f)&=& \operatorname{tr}f\bigl(\Sigma_n
\Sigma_n^*\bigr) - N \int f(\lambda) {\mathcal F}_n(d
\lambda)
\nonumber
\\[-8pt]
\label{eqcauchy}
\\[-8pt]
\nonumber
&=& \frac{1}{\pi} \operatorname{Re}\int_\mathrm{\mathbb{C}^+}
\overline{\partial}\Phi_k(f) (z) \bigl\{ \operatorname{tr}Q_n(z)
- N t_n(z) \bigr\} \ell_2(d z),
\nonumber
\end{eqnarray}
where $\Phi_k(f)$ is defined in \eqref{eqalmost-analytic} and the
last equality follows from the fact that
\[
\int f(\lambda){\mathcal F}_n(d\lambda) = \frac{1}{\pi}
\operatorname{Re} \int_{\mathbb{C}^+} \overline{\partial}
\Phi_k(f) (z) t_n(z) \ell_2(dz).
\]
Based on \eqref{eqcauchy}, we shall first study the fluctuations of
\begin{eqnarray*}
&& \operatorname{tr}Q_n(z) - Nt_n(z) = \bigl\{
\operatorname{tr}Q_n(z) - \mathbb{E}\operatorname{tr}Q_n(z)
\bigr\} + \bigl\{ \mathbb{E}\operatorname{tr}Q_n(z) -
Nt_n(z) \bigr\}
\end{eqnarray*}
for $z\in\mathbb{C}^+$. The first difference in the RHS will yield the
fluctuations with a covariance $\Theta_n(z_1,z_2)$ described in \eqref
{eqdecomposition-covariance}
while the second difference, deterministic, will yield the bias
%
\begin{eqnarray}
\mathbb{E}\operatorname{tr}Q_n(z) - Nt_n(z) &=& |
\mathcal{V}|^2 {\mathcal B}_{1,n}(z) + \kappa{\mathcal
B}_{2,n}(z) + o(1)
\nonumber
\\[-8pt]
\label{defB}
\\[-8pt]
\nonumber
&\stackrel{\triangle} {=}& {\mathcal B}_n(z) + o(1),
\end{eqnarray}
where
\begin{eqnarray}
&& {\mathcal B}_{1,n}(z)
\nonumber
\\
&&\qquad\stackrel{\triangle} {=} - z^3
\tilde t_n^3 \biggl({\frac{1}n \operatorname{tr}
R_n^{1/2} T^2_n(z) R_n^{1/2} \overline{R}_n^{1/2} T^T_n(z) \overline
{R}_n^{1/2} }\biggr)
\nonumber
\\[-8pt]
\label{defB1}
\\[-8pt]
\nonumber
&&\qquad\quad{}\Big/\biggl(
\biggl( 1 - z^2 \tilde t_n^2 \frac{1}n \operatorname{Tr}R_n^2
T_n^2
\biggr)\\
\nonumber
&&\qquad\quad{}\times \biggl( 1 - |\mathcal{V}|^2 z^2 \tilde t_n^2 \frac{1}n \operatorname
{Tr}R_n^{1/2}
T_n(z) R_n^{1/2} \overline{R}_n^{1/2} T^T_n(z) \overline{R}_n^{1/2} \biggr)\biggr),
\nonumber\\
\label{defB2}
&& {\mathcal B}_{2,n}(z) \stackrel{\triangle} {=} - z^3
\tilde t_n^3 \frac{({1}/n) \sum_{i=1}^N ( R^{1/2}_n T_n R^{1/2}_n
)_{ii} ( R^{1/2}_n T^2_n R^{1/2}_n )_{ii}
}{
1 - z^2 \tilde t_n^2 ({1}/{n}) \operatorname{tr}R^2_n
T^2_n}.
\end{eqnarray}
The previous discussion on the terms $\Theta_{1,n}$ and $\Theta
_{2,n}$ also applies to the terms ${\mathcal B}_{1,n}$ and ${\mathcal
B}_{2,n}$ (whose expressions are established in Section~\ref{secproof})
which are likely not to converge for similar reasons.

\subsection{Gaussian processes and the central limit theorem}\label
{secconvergence-preliminary} A priori, the mean ${\mathcal B}_n$ and
covariance $\Theta_n$ of $ (\operatorname{tr}Q_n - Nt_n )$
do not
converge. Hence, we shall express the Gaussian fluctuations of the
linear statistics \eqref{eqlinear-stat-centered} in the following
way: we first prove the existence of a family $(G_n(z), z\in
{\mathcal C})_{n\in\mathbb{N}}$ of tight Gaussian processes with mean
and covariance
\begin{eqnarray*}
\mathbb{E}G_n(z) &=& {\mathcal B}_n(z),
\\
\operatorname{cov}\bigl(G_n(z_1),G_n(z_2)
\bigr)&=&\Theta_n(z_1,z_2).
\end{eqnarray*}
%
We then express the fluctuations of the centralized trace as
\[
d_{\mathcal{LP}} \bigl( \bigl(\operatorname{tr}Q_n(z) -
Nt_n(z) \bigr), G_n(z) \bigr) \mathop{
\xrightarrow}_{N,n\rightarrow\infty} 0,
\]
with $d_{\mathcal{LP}}$ the L\'{e}vy--Prohorov distance between $P$
and $Q$
probability measures over borel sets of $\mathbb{R},\mathbb
{R}^d,\mathbb{C}$ or $\mathbb{C}^d$:
%
\begin{equation}
\label{eqdef-LP}
\hspace*{6pt}d_{\mathcal{LP}}(P,Q)= \inf\bigl\{ \varepsilon>0,
P(A) \le Q
\bigl(A^\varepsilon\bigr) + \varepsilon\mbox{ for all Borel sets }A\bigr\},
\end{equation}
where $A^\varepsilon$ is an $\varepsilon$-blow up of $A$ (cf. \cite{book-dudley}, Section~11.3,  for more details).
If $X$ is a random variable and ${\mathcal L}(X)$ its distribution,
denote (with a slight abuse of notation) by $d_{\mathcal{LP}}(X,Y)
\stackrel{\triangle}{=}d_{\mathcal{LP}}
({\mathcal L}(X),{\mathcal L}(Y))$.

Similarly, we will express the fluctuations of $L_n(f)$ as
\[
d_{\mathcal{LP}} \bigl( L_n(f) 
, {\mathcal
N}_n(f) \bigr) \mathop{\xrightarrow}_{N,n\rightarrow\infty} 0,
\]
where ${\mathcal N}_n(f)$ is a well-identified Gaussian random variable.

\subsection{A meta-model argument}\label{secmeta-model} As we need
to cope with a sequence of Gaussian processes $(G_n)$ instead of
a single one, it will be necessary to establish various properties
uniform in $n,N$ such as:
\begin{longlist}[3.]
\item[1.] the tightness of the sequence $(G_n)$ (cf. Section~\ref
{secgaussian-process});
\item[2.] a uniform bound over the variances of $(\operatorname{Tr}G_n(z))$
(cf. Section~\ref{proofnon-analytic-I}), needed to extend the CLT to
nonanalytic functionals;
\item[3.] a uniform bound over the biases of $(\operatorname{Tr}G_n(z))$ (cf.
Section~\ref{secMM-bias}), needed to compute the bias for
nonanalytic functionals.
\end{longlist}

A direct approach based on the mere definition of process $G_n$'s
parameters seems difficult, mainly due to the definitions of $\Theta
_n$ and ${\mathcal B}_n$ which rely on quantities ($t_n$ and $\tilde
t_n$) defined as solutions of fixed-point equations. Since the previous
properties will be established for the processes $(\operatorname
{Tr}Q_n -Nt_n)$
anyway, the idea is to transfer them to $G_n$ by means of the
following matrix meta-model.

Let $N$, $n$ and $R_n$ be fixed and consider the $NM\times NM$ matrix
%
\begin{equation}
\label{eqdef-RM} R_n(M) = \pmatrix{
R_n & 0 & \cdots
\cr
0&\ddots&0
\cr
\cdots& 0 & R_n}.
\end{equation}
Matrix $R_n(M)$ is a block matrix with $N\times N$ diagonal blocks
equal to $R_n$, and zero blocks elsewhere; for all $M\ge1$ the
spectral norm of $R_n(M)$ is equal to the spectral norm of $R_n$ (which
is fixed).
In particular, the sequence $(R_n(M); M\ge1)$ with $N,n$ fixed
satisfies Assumption \ref{assR} with $(R_n(M); M\ge1)$ instead
of $(R_n)$. Consider now the random matrix model
%
\begin{equation}
\label{eqextended-model}
\Sigma_n(M) = \frac{1}{\sqrt{Mn}} R_n(M)^{1/2}
X_n(M),
\end{equation}
where $X_n(M)$ is a $MN\times Mn$ matrix with i.i.d. random entries
with the same distribution as the $X_{ij}$'s and satisfying  Assumption~\ref{assX}.
The interest of introducing matrix $\Sigma_n(M)$ lies in the
fact that matrices $\Sigma_n(M) \Sigma_n(M)^*$ and $\Sigma_n \Sigma
_n^*$ have loosely speaking the same deterministic equivalents. Denote
by $t_n$, $T_n$ and $\tilde t_n$ the deterministic equivalents of
$\Sigma_n \Sigma_n^*$ as defined in \eqref{eqdeteq-ST}, \eqref
{eqdeteq-RESOLVENT} and \eqref{eqdef-t-tilde}, and by $t_n(M)$,
$T_n(M)$ and $\tilde t_n(M)$ their counterparts for the model $\Sigma
_n(M) \Sigma_n(M)^*$. Taking advantage of the block structure of
$R_n(M)$, a straightforward computation yields ($N,n$ fixed)
\[
\forall M\ge1,\qquad t_n(M)=t_n, \qquad\tilde t_n(M) =
\tilde t_n \quad\mbox{and}\quad T_n(M) =  \pmatrix{
T_n & 0 & \cdots
\cr
&\ddots&
\cr
\cdots& 0 & T_n
}.
\]
Similarly, denote by ${\mathcal B}_{n,M}$ and $\Theta_{n,M}$ the
quantities given by formulas \eqref{defB} and~\eqref
{eqdecomposition-covariance} when replacing $N$, $t_n$, $T_n$ and
$\tilde t_n$ by $NM$, $t_n(M)$, $T_n(M)$ and $\tilde t_n(M)$.
Straightforward computation yields
\[
\forall M\ge1, \qquad {\mathcal B}_{n,M}={\mathcal B}_n \quad\mbox{and}\quad \Theta_{n,M}=\Theta_n.
\]
An interesting feature of this meta-model lies in the fact that all the
quantities associated to $\Sigma_n(M)\Sigma_n(M)^*$ converge as $M\to
\infty$ to the deterministic equivalents $t_n$, $\tilde t_n$, etc. As
a consequence,
one can easily transfer all the estimates obtained for
\[
\bigl( \operatorname{Tr}\bigl(\Sigma_n(M)\Sigma_n(M)^*-zI_{NM}
\bigr)^{-1} - NM t_n \bigr)
\]
to the process $(G_n)$.

\section{Statement of the CLT for the trace of the resolvent}\label
{sectionmain}

\subsection{Further notation}

If $A$ is a $N\times N$ matrix with real eigenvalues, denote by $F^A$
the empirical distribution of the eigenvalues ($\delta_i(A),i=1:N$) of
$A$, that is,
\[
F^A(dx)=\frac{1}n \sum_{i=1}^N
\delta_{\lambda_i(A)}(dx).
\]

Recall the definitions of $Q_n$, $t_n$, $T_n$ and $\tilde t_n$ [cf.
\eqref{eqresolvent}, \eqref{eqdeteq-ST}, \eqref
{eqdeteq-RESOLVENT} and \eqref{defttilde}].
The following relations hold true (see, e.g., \cite
{bai-silverstein-aop-2004}):
%
\begin{equation}\label{eqdef-t-tilde}
\qquad T_n(z) = -\frac{1}z \bigl( I_N + \tilde
t_n(z) R_n \bigr)^{-1} \quad\mbox{and}\quad \tilde
t_n(z) = -\frac{1}{z ( 1 +({1}/n) \operatorname{tr}R_nT_n(z) )}.\hspace*{-20pt}
\end{equation}

Recall the definition of ${\mathcal F}_n$ in \eqref
{eqstieltjes-basics} and let similarly $\widetilde{\mathcal F}_n$ be the
probability distribution associated to $\tilde t_n$. The central object
of study is the signed measure
\[
N \bigl( F^{\Sigma_n \Sigma_n^*} - {\mathcal F}_n \bigr)=n
\bigl( F^{\Sigma_n^* \Sigma_n} - \widetilde{\mathcal F}_n \bigr),
\]
%
and its Stieltjes transform
%
\begin{equation}
\label{eqprocess-def}
M_n(z) =N\bigl(f_n(z) - t_n(z)
\bigr)= n \bigl( \tilde f_n(z) - \tilde t_n(z) \bigr).
\end{equation}

Denote by $o_P(1)$ any random variable which converges to zero in
probability as $N,n\to\infty$.

\subsection{Truncation}\label{sectruncation} In this section, we
closely follow Bai and Silverstein \cite{bai-silverstein-aop-2004}. We
recall the framework developed there and introduce some additional notation.

Consider a sequence of positive numbers $(\delta_n)$ which satisfies
\[
\delta_n\rightarrow0,\qquad \delta_n n^{1/4}
\rightarrow\infty\quad\mbox{and}\quad \delta_n^{-4} \int
_{ \{ | X_{11}| \ge\delta
_n \sqrt{N} \}} |X_{11}|^4 \rightarrow0
\]
as $N,n\rightarrow\infty$. Let $\widehat\Sigma_n= n^{-1/2}
R_n^{1/2} \widehat{X}_n$ where $\widehat{X}_n$ is a $N\times n$
matrix having $(i,j)$th entry $X_{ij} 1_{\{ |X_{ij}|< \delta_n \sqrt
{N}\} }$. This truncation step yields
%
\begin{equation}
\label{eqtruncation-1}
\mathbb{P} \bigl(\Sigma_n \Sigma_n^* \neq
\widehat{\Sigma}_n \widehat{\Sigma}_n^* \bigr)\mathop{
\xrightarrow}_{N,n\rightarrow\infty} 0
\end{equation}
from which we deduce
%
\begin{equation}
\label{eqtruncated-variables}
\operatorname{tr} \bigl(\Sigma_n \Sigma_n^*
- z I_N\bigr)^{-1} - \operatorname{tr} \bigl(\widehat{
\Sigma}_n \widehat{\Sigma}_n^* - z I_N
\bigr)^{-1} \mathop{\xrightarrow}_{N,n\to\infty}^{\mathcal P} 0,
\end{equation}
where $ \mathop{\rightarrow}^{\mathcal P}$ stands for the convergence
in probability.
Define $\widetilde\Sigma_n=n^{-1/2} R_n^{1/2} \widetilde{X}_n$ where
$\widetilde{X}_n$ is a $N\times n$ matrix having $(i,j)$th entry
$(\widehat X_{ij} - \mathbb{E} \widehat{X}_{ij}) / \sigma_n$, where
$\sigma_n^2=\mathbb{E} | \widehat{X}_{ij} - \mathbb{E} \widehat
{X}_{ij}|^2$. Using the fact that $\lambda(\in\mathbb{R}) \mapsto
\frac{1}{\lambda-z}$ is Lipschitz with Lipschitz constant $|z|^{-2}$,
we obtain
\[
\mathbb{E}\bigl\llvert\operatorname{tr}\bigl( \widehat\Sigma_n
\widehat\Sigma_n^* - zI_N\bigr)^{-1} -
\operatorname{tr}\bigl( \widetilde\Sigma_n \widetilde\Sigma
_n^* - zI_N\bigr)^{-1}\bigr\rrvert\le
\frac{1}{|z|^2} \sum_{i=1}^N \mathbb{E}
\llvert\tilde\lambda_i - \hat\lambda_i \rrvert\mathop{
\longrightarrow}_{N,n\to\infty}^{(a)} 0,
\]
where $\tilde\lambda_i = \lambda_i(\widetilde\Sigma_n \widetilde
\Sigma_n^*)$, $\hat\lambda_i = \lambda_i(\widehat\Sigma_n
\widehat\Sigma_n^*)$ and $(a)$ follows from similar arguments as in
\cite{book-bai-silverstein}, Section~9.7.1. Hence,
%
\begin{equation}
\label{eqcentered-variables}
\operatorname{tr} \bigl(\widehat\Sigma_n
\widehat
\Sigma_n^* - z I_N\bigr)^{-1} -
\operatorname{tr} \bigl(\widetilde{\Sigma}_n \widetilde{
\Sigma}_n^* - z I_N\bigr)^{-1} \mathop{
\longrightarrow}_{N,n\to\infty}^{\mathcal P} 0.
\end{equation}

Combining \eqref{eqtruncated-variables} and \eqref
{eqcentered-variables}, we obtain
\[
\operatorname{tr} \bigl( \Sigma_n \Sigma_n^* - z
I_N\bigr)^{-1} - \operatorname{tr} \bigl(\widetilde{\Sigma
}_n \widetilde{\Sigma}_n^* - z I_N
\bigr)^{-1} \mathop{\longrightarrow}_{N,n\to\infty
}^{\mathcal P} 0.
\]
Moreover, the moments are asymptotically not affected by these
different steps:
%
\begin{equation}
\max\bigl( \bigl\llvert\mathbb{E} \widetilde X_{ij}^2 -
\mathbb{E} X_{ij}^2\bigr\rrvert; \bigl( \mathbb{E} |
\widetilde X_{ij}|^2 - 1 \bigr) ; \bigl( \mathbb{E}|
\widetilde X_{ij}|^4 - \mathbb{E} |X_{ij}|^4
\bigr) \bigr)\mathop{\longrightarrow}_{N,n\rightarrow\infty} 0.
\end{equation}
Note in particular that the fourth cumulant of $\widetilde X_{ij}$
converges to that of $X_{ij}$.
Hence, it is sufficient to consider variables truncated at $\delta
_n\sqrt{n}$, centralized and renormalized. This will be assumed in the
sequel (we shall simply write $X_{ij}$ and all related quantities with
$X_{ij}$'s truncated, centralized, renormalized with no superscript any more).

\subsection{The central limit theorem for the resolvent}

We extend below Bai and Silverstein's master lemma \cite{bai-silverstein-aop-2004}, Lemma~1.1. Let $A$ be such that
\[
A> \bolds{\lambda}^+_R \bigl( 1 + \sqrt{\bolds{\ell}^+}
\bigr)^2.
\]
%
%
Denote by $D$, $D^+$ and $D_\varepsilon$ the domains
%
\begin{eqnarray}
D & =& [0,A]+\mathbf{i}[0,1],
\nonumber
\\
\label{defsets-D}
D^+&=& [0,A]+\mathbf{i}(0,1],
\\
D_{\varepsilon}&=& [0,A] +\mathbf{i}[\varepsilon,1]\qquad (\varepsilon>0).
\nonumber
\end{eqnarray}

\begin{theorem} \label{lemmamain}
Assume that Assumption~\ref{assX} and
Assumption~\ref{assR} hold true, then:
\begin{longlist}[3.]
\item[1.] The process $\{M_n(\cdot)\}$ as defined in \eqref
{eqprocess-def} forms a tight sequence on $D_\varepsilon$, more
precisely,
\[
\sup_{z_1,z_2\in D_{\varepsilon},n\ge1} \frac{\mathbb{E}\llvert
M_n(z_1) - M_n(z_2)\rrvert^2}{|z_1- z_2|^2} <\infty.
\]
\item[2.]
There exists a sequence $(G_n(z), z\in D^+)$ of two-dimensional
Gaussian processes with mean
%
\begin{equation}
\mathbb{E}G_n(z) = |\mathcal{V}|^2 {\mathcal
B}_{1,n}(z) + \kappa{\mathcal B}_{2,n}(z),
\end{equation}
where ${\mathcal B}_{1,n}(z)$ and ${\mathcal B}_{2,n}(z)$ are defined
in \eqref{defB1} and \eqref{defB2}, and covariance
\begin{eqnarray*}
\operatorname{cov} \bigl( G_n(z_1),
G_n(z_2) \bigr) &=& \mathbb{E} \bigl(G_n(z_1)
- \mathbb{E}G_n(z_1) \bigr) \bigl(G_n(z_2)
- \mathbb{E}G_n(z_2) \bigr)
\nonumber
\\
&=& \Theta_{0,n}(z_1,z_2) +\llvert{\mathcal
V} \rrvert^2\Theta_{1,n}(z_1,z_2)
+ \kappa\Theta_{2,n}(z_1,z_2),
\end{eqnarray*}
and
\begin{eqnarray*}
\operatorname{cov} \bigl( G_n(z_1),
\overline{G_n(z_2)} \bigr) &=&\operatorname{cov} \bigl(
G_n(z_1), G_n(\overline{z_2})
\bigr),
\end{eqnarray*}
with $z_1,z_2\in D^+ \cup\overline{D^+}$ with $\overline{D^+}=\{
\bar{z}, z\in D^+\}$ and where $\Theta_{0,n}$, $\Theta_{1,n}$
and $\Theta_{2,n}$ are defined in \eqref{eqdecomposition-covariance},
\eqref{deftheta0}--\eqref{deftheta2}.
Moreover, $(G_n(z), z\in D_{\varepsilon})$ is tight.

\item[3.] For any continuous functional $F$ from $C (D_{\varepsilon
};\mathbb{C} )$ to $\mathbb{C}$,
\[
\mathbb{E}F( M_n) - \mathbb{E}F ( G_n ) \mathop{
\longrightarrow}_{N,n\rightarrow\infty} 0.
\]
\end{longlist}
\end{theorem}
\begin{rem}
1. The tightness of the process $\{ M_n\}$ immediately follows from
Bai and Silverstein's lemma as this result has been proved in \cite{bai-silverstein-aop-2004}, Lemma~1.1,
under Assumption~\ref{assX}
with no extra conditions on the moments of the entries.

2. Differences between Theorem~\ref{lemmamain} and \cite{bai-silverstein-aop-2004}, Lemma~1.1, appear in the bias and in the covariance
where there are respectively two terms instead of one and three terms
instead of one in \cite{bai-silverstein-aop-2004}, Lemma~1.1.

3. Since the extra terms may not converge, we need to consider a
sequence of Gaussian processes instead of a single Gaussian process as
in \cite{bai-silverstein-aop-2004}, Lemma~1.1.

4. In order to prove that the sequence of Gaussian processes is
tight, we introduce a meta-matrix model to transfer the tightness of
$\{ M_n\}$ to $\{ G_n\}$ (see, e.g., Section~\ref{secmeta-matrix}).

5. Following Bai and Silverstein \cite{bai-silverstein-aop-2004},
it is relatively straightforward with the help of Cauchy's formula to
describe the fluctuations of $L_n(f)$ for $f$ analytic with Theorem~\ref
{lemmamain} at hand. We skip this step since we will directly
extend the CLT to nonanalytic functions $f$ in Section~\ref{secnon-analytic}.
%
\end{rem}

\begin{rem} \label{remaltern-theta0} A closer look to Bai and
Silverstein's proof \cite{bai-silverstein-aop-2004}, Section~2,
page~578,
yields the following alternative expression for the term $\Theta_{0,n}$:
%
\begin{equation}\label{deftheta0-alternate}
\Theta_{0,n}(z_1,z_2) = \frac{\partial}{\partial z_2}
\biggl\{ \frac{\partial{\mathcal A}_{0,n}(z_1,z_2)} {\partial z_1} \frac
{1}{1-{\mathcal A}_{0,n}(z_1,z_2)} \biggr\},
\end{equation}
with
%
\begin{equation}
\label{defAn0}
{\mathcal A}_{0,n}(z_1,z_2) =
\frac{z_1 z_2} n \tilde t_n(z_1) \tilde
t_n(z_2) \operatorname{tr} \bigl\{ R_n
T_n(z_1) R_n T_n(z_2)
\bigr\}.
\end{equation}
Such an\vspace*{1pt} expression will be helpful in Section~\ref
{proofnon-analytic-I}. As an interesting consequence: In the case
where $R_n$ and $X_n$
have real entries [in particular ${\mathcal V} = \mathbb{E}(X_{ij})^2
= 1$],
then ${\mathcal A}_{0,n} = {\mathcal A}_n$ and $\Theta_{0,n}= \Theta_{1,n}$.
\end{rem}

\begin{rem}\label{remaltern-theta2}
A closer look to the proof below
[see, e.g., \eqref{alternative-formula-cov2-bis}] yields the
following formula for $\Theta_{2,n}$ which will be of help in the sequel:
%
\begin{equation}
\label{alternative-cov2}
\Theta_{2,n}(z_1,z_2) =
\frac{1}n \sum_{i=1}^N
\frac{ \partial
}{\partial z_1} \bigl[ z_1 T_n(z_1)
\bigr]_{ii} \frac{ \partial}{\partial z_2} \bigl[z_2
T_n(z_2) \bigr]_{ii}.
\end{equation}
\end{rem}

The proof of Theorem~\ref{lemmamain} is postponed to Section~\ref{secproof}.

The end of the section is devoted to various specializations of
Theorem~\ref{lemmamain} in the case where matrix $R_n$ is diagonal. In this
case, the results are simpler to express and comparisons can easily be
made with related works.

\subsection{Covariance and bias in the special case of diagonal
matrices $(R_n)$}\label{secdiagonal-R} This case partially falls into
the framework developed in Pan and Zhou \cite{pan-zhou-2008} (note
that the case $\mathcal{V}\neq0$ and $1$ is not handled there).
Matrix $R_n$ being nonnegative definite Hermitian, its entries are
real positive if $R_n$ is assumed to be diagonal. In this case, matrix
$T_n$ is diagonal as well [cf. \eqref{eqdeteq-RESOLVENT}],
$T_n=T^T_n$ and simplifications occur for the following terms:
\begin{eqnarray*}
{\mathcal A}_n(z_1,z_2) &=&
\frac{z_1 z_2}n \tilde t_n(z_1) \tilde
t_n(z_2) \operatorname{tr}R_n
T_n(z_1) R_n T_n(z_2),
\\
\Theta_{2,n}(z_1,z_2)&=& \frac{ z_1^2 z_2^2 \tilde t'_n(z_1) \tilde
t'_n(z_2)}n
\operatorname{tr} \bigl( R_n^2 T^2_n(z_1)
T^2_n(z_2) \bigr),
\\
{\mathcal B}_{1,n}(z)&=& - z^3 \tilde t_n^3
\frac{({1}/n)
\operatorname{tr}R_n^2 T_n^3}{
( 1 - z^2 \tilde t_n^2 ({1}/n) \operatorname{Tr}R_n^2
T_n^2
) ( 1 - |\mathcal{V}|^2 z^2 \tilde t_n^2 ({1}/n) \operatorname
{Tr}R_n^2 T^2_n
)},
\\
{\mathcal B}_{2,n}(z) &=& - z^3 \tilde
t_n^3 \frac{
({1}/n) \operatorname{tr}R^2_n T^3_n
}{
1 - z^2 \tilde t_n^2 ({1}/{n}) \operatorname{tr}R^2_n T^2_n}.
\end{eqnarray*}
As one may notice, all the terms in the variance and the bias now only
depend on the spectrum of $R_n$. Hence, the following convergence holds
true under the extra assumption \eqref{eqadd-assumption}:
\begin{eqnarray*}
{\mathcal A}_n(z_1,z_2) &\displaystyle \mathop{
\xrightarrow}_{N,n\rightarrow\infty} & {\mathcal A}(z_1,z_2) = c
\tilde t(z_1) \tilde t(z_2) \int\frac
{\lambda^2F^{\mathbf R}(d\lambda)}{(1+\lambda\tilde t(z_1))
(1+\lambda\tilde t(z_2))},
\\
\Theta_{1,n}(z_1,z_2) &\displaystyle \mathop{
\xrightarrow}_{N,n\rightarrow\infty}& \Theta_1(z_1,z_2)
=\frac{\partial}{\partial z_2} \biggl\{ \frac{\partial
{\mathcal A}(z_1,z_2)} {\partial z_1} \frac{1}{1-|\mathcal
{V}|^2{\mathcal
A}(z_1,z_2)} \biggr\},
\\
\Theta_{2,n}(z_1,z_2) &\displaystyle \mathop{
\xrightarrow}_{N,n\rightarrow\infty}& \Theta_2(z_1,z_2)
=c \tilde t'(z_1) \tilde t'(z_2)
\int\frac{\lambda
^2F^{\mathbf R}(d\lambda)}{(1+\lambda\tilde t(z_1))^2
(1+\lambda\tilde t(z_2))^2},
\\
{\mathcal B}_{1,n}(z) &\displaystyle \mathop{\xrightarrow}_{N,n\rightarrow\infty}&
{\mathcal
B}_1(z)= -\frac{cz^3 \tilde t^3(z)}{ (1-{\mathcal
A}(z,z))(1-|\mathcal{V}|^2{\mathcal A}(z,z))} {\int\frac{\lambda
^2F^{\mathbf
R}(d\lambda)}{(1+\lambda\tilde t(z))^3
}},
\\
{\mathcal B}_{2,n}(z) & \displaystyle\mathop{\xrightarrow}_{N,n\rightarrow\infty}& {
\mathcal B}_2(z)= -\frac{cz^3 \tilde t^3(z)}{ 1-{\mathcal A}(z,z)} {\int
\frac{\lambda^2F^{\mathbf R}(d\lambda)}{(1+\lambda\tilde t(z))^3
}},
\end{eqnarray*}
where $\tilde t,\tilde t'$ are the limits of $\tilde t_n,\tilde t_n'$
under \eqref{eqadd-assumption}. This can be packaged into the
following result.

\begin{coro}\label{corodiag} Assume that Assumptions \ref{assX} and \ref{assR} hold true.
Assume moreover that $R_n$ is diagonal and that the
convergence assumption \eqref{eqadd-assumption} holds true. Then
$M_n(\cdot)$ converges weakly on $D_\varepsilon$ [defined in \eqref
{defsets-D}] to a two-dimensional Gaussian process $N(\cdot)$ satisfying
\[
\mathbb{E}N(z) = {\mathcal B}(z) \qquad\mbox{where } {\mathcal B} = |
\mathcal{V}|^2 {\mathcal B}_1 +\kappa{\mathcal
B}_2,  z\in D_{\varepsilon}
\]
and ${\mathcal B}_1$ and ${\mathcal B}_2$ are defined above and covariance
%
\begin{eqnarray}
&& \operatorname{cov} \bigl( N(z_1), N_(z_2) \bigr) =
\Theta(z_1,z_2)
\nonumber\\
\eqntext{\displaystyle\mbox{where } \Theta=
\Theta_0 +|\mathcal{V}|^2\Theta_1 + \kappa
\Theta_2, z_1,z_2\in D_{\varepsilon} \cup
\overline{D_{\varepsilon}}}
\end{eqnarray}
and $\Theta_0$ defined in \eqref{eqdef-theta0} and $\Theta_1,\Theta
_2$ defined above.
\end{coro}

\subsection{Additional computations in the case where $R_n$ is the
identity}\label{secwhite-case} In this section, we assume that $R_n=I_N$.
\subsubsection*{The term proportional to $|\mathcal{V}|^2$} In this
case, the
quantity ${\mathcal A}(z_1,z_2)$ takes the simplified form
\[
{\mathcal A}(z_1,z_2) = \frac{c \tilde t_1 \tilde t_2}{(1+ \tilde
t_1)(1+\tilde t_2)},
\]
where we denote $\tilde t_i = \tilde t(z_i), i=1,2$. Straightforward
computations yield
\[
\frac{\partial}{\partial z_i} {\mathcal A} (z_1,z_2)=
\frac{\tilde
t'_i}{(1+\tilde t_i)\tilde t_i} {\mathcal A}(z_1,z_2),\qquad i=1,2
\]
and
\begin{eqnarray*}
\Theta_1(z_1,z_2) &=& \frac{c \tilde t'_1\tilde t'_2}{(1+\tilde
t_1)^2(1+\tilde t_2)^2 ( 1 - |\mathcal{V}|^2 {\mathcal
A}(z_1,z_2)
)^2}
\\
&=& \frac{c \tilde t'_1\tilde t'_2}{ ( (1+\tilde t_1)(1+\tilde
t_2) - |\mathcal{V}|^2 c \tilde t_1 \tilde t_2 )^2}.
\end{eqnarray*}
This formula is in accordance with \cite{pan-CLT-MIMO-preprint-2013},
formula (2.2) (use \cite{pan-CLT-MIMO-preprint-2013}, equation (3.4),
to equate both).
If needed, one can then use the explicit expression of the Stieltjes
transform of Mar\v cenko--Pastur distribution (cf. also Proposition~\ref
{propnew-variance-formula} below).

\section{Statement of the CLT for nonanalytic functionals}\label
{secnon-analytic}

In order to lift the CLT from the trace of the resolvent to a smooth
function $f$, the key ingredient is Helffer--Sj\"{o}strand's formula
\eqref{eqhelffer-sjostrand}. Let
\begin{eqnarray}
L_n(f) &\stackrel{(a)}=& \operatorname{Tr}f\bigl(
\Sigma_n\Sigma_n^*\bigr) - N\int f(\lambda) {\mathcal
F}_n(d\lambda)
\nonumber\\
&=& \bigl( \operatorname{Tr}f\bigl(\Sigma_n\Sigma_n^*
\bigr) - \mathbb{E}\operatorname{Tr}f\bigl(\Sigma_n
\Sigma_n^*\bigr) \bigr)
\nonumber
\\[-8pt]
\label{defLn}
\\[-8pt]
\nonumber
&& {}+ \biggl(\mathbb{E} \operatorname{Tr}f\bigl(\Sigma_n
\Sigma_n^*\bigr) - N\int f(\lambda) {\mathcal F}_n(d
\lambda) \biggr)
\nonumber
\\
&\stackrel{\triangle}=& L_n^1(f) + L_n^2(f),
\nonumber
\end{eqnarray}
where ${\mathcal F}_n$ in $(a)$ is defined in \eqref
{eqstieltjes-basics}. We
describe the fluctuations of $L_n^1(f)$ for nonanalytic functions $f$
in Section~\ref{secfluctuations-non-analytic} and study the bias $L_n^2(f)$
in Section~\ref{secbias-non-analytic}.

\subsection{Fluctuations for the linear spectral statistics}\label
{secfluctuations-non-analytic}
Denote by $C_c^\infty(\mathbb{R}^d)$ [resp., $C_c^{m}(\mathbb{R}^d)$]
the set of
infinitely differentiable (resp., $C^m$) functions from $\mathbb{R}^d$
to $\mathbb{R}
$ with compact support; by $C_c^{mp}(\mathbb{R}^2)$ the set of
functions from
$\mathbb{R}^2$ to $\mathbb{R}$ $m$ times differentiable with respect
to the first
coordinate and $p$ times with respect to the second one.
As usual, if the subscript $c$ is removed in the sets above, then the
corresponding functions may no longer have a compact support.


\begin{theorem}\label{thnon-analytic}
Assume that \ref{assX}  and \ref{assR}  hold true. Let
$f_1,\ldots, f_k$ be in $C^3_c(\mathbb{R})$.
Consider the centered Gaussian random vector
$Z^1_n({\mathbf f})\stackrel{\triangle}{=}(Z^1_n(f_1),\ldots
,\break Z^1_n(f_k))$ with covariance
%
\begin{eqnarray}
&& \operatorname{cov} \bigl( Z^1_n(f),
Z^1_n(g) \bigr)
\nonumber
\\
\label{eqcov-gaussian}
&&\qquad = \frac{1}{2\pi^2} \operatorname{Re} \int_{(\mathbb{C}^+)^2}
\overline{\partial} \Phi_2(f) (z_1) \overline{\overline{
\partial} \Phi_2(g) (z_2)} \Theta_n(z_1,
\bar{z}_2) \ell_2(dz_1)\ell_2(dz_2)
\\
&&\qquad\quad{}+ \frac{1}{2\pi^2} \operatorname{Re} \int_{(\mathbb{C}^+)^2} \overline{
\partial} \Phi_2(f) (z_1) \overline{\partial} \Phi
_2(g) (z_2) \Theta_n(z_1,z_2)
\ell_2(dz_1)\ell_2(dz_2),\nonumber
\end{eqnarray}
for $f,g\in\{ f_1,\ldots,f_k\}$, where $\Phi_2(f)$ and $\Phi_2(g)$
are defined as in \eqref{eqalmost-analytic}, and where $\Theta_n$ is
defined in \eqref{eqdecomposition-covariance}; let
\[
L^1_n({\mathbf f})=\bigl(L^1_n(f_1),
\ldots, L^1_n(f_k)\bigr) \qquad\mbox{with }
L_n^1(f) = \operatorname{tr}f\bigl(\Sigma_n
\Sigma_n^*\bigr) - \mathbb{E} \operatorname{tr} f\bigl(
\Sigma_n \Sigma_n^*\bigr).
\]
%
Then the sequence of $\mathbb{R}^k$-valued random vectors
$Z^1_n({\mathbf
f})$ is
tight and the following convergence holds true:
%
\begin{equation}
\label{eqfluctuations-1}
d_{\mathcal{LP}} \bigl( L^1_n({\mathbf
f}),Z^1_n({\mathbf f}) \bigr) \mathop{
\xrightarrow}_{N,n\rightarrow\infty} 0,
\end{equation}
%
or equivalently for every continuous bounded function $F:\mathbb{R}
^k\rightarrow\mathbb{C}$,
%
\begin{equation}
\label{eqfluctuations-2}
\mathbb{E} F\bigl(L^1_n({\mathbf f})
\bigr) - \mathbb{E} F\bigl(Z^1_n({\mathbf f})\bigr)
\mathop{\xrightarrow}_{N,n\rightarrow\infty} 0.
\end{equation}
%
%
\end{theorem}

The proof of Theorem~\ref{thnon-analytic} is postponed to Section~\ref
{secproof-non-analytic-I}.

We provide hereafter some information on the covariance operator.

Let $N_1,N_2\in\mathbb{N}$ and $f\in C^{N_1+1,N_2+1}_c(\mathbb
{R}^2)$; denote by
$z_1=x+\mathbf{i}u$, $z_2= y+\mathbf{i}v$ and let $\Phi
_{N_1,N_2}(f)$ be defined as
%
\begin{equation}
\label{defAAE-bidim}\quad
\Phi_{N_1,N_2}(f) (z_1,z_2) = \mathop{\sum
_{ n_1=0:N_1}}_{n_2=0:N_2} \frac{\partial^{n_1+n_2}}{\partial x^{n_1} \partial y^{n_2}} f(x,y)
\frac{(\mathbf{i}u)^{n_1}}{n_1!}
\frac{(\mathbf
{i}v)^{n_2}}{n_2!}\chi(u)\chi(v),
\end{equation}
where $\chi:\mathbb{R}\to\mathbb{R}^+$ is smooth, compactly
supported with value 1
in a neighborhood of the origin. Denote by $\overline{\partial}_1 =
\partial_x +\mathbf{i}\partial_u$ and $\overline{\partial}_2=
\partial_y
+\mathbf{i}\partial_v$.

\begin{prop}\label{propidentification-distribution} For every $f\in
C^{3,3}_c(\mathbb{R}^2)$, denote by
\begin{eqnarray*}
\Upsilon(f) &=& \frac{1}{2\pi^2} \operatorname{Re} \int_{(\mathbb{C}^+)^2}
\overline{\partial}_2\overline{\partial}_1
\Phi_{2,2} (f) (z_1,z_2)\Theta
_n(z_1,z_2) \ell_2(dz_1)
\ell_2(dz_2)
\\
&&{}+\frac{1}{2\pi^2} \operatorname{Re} \int_{(\mathbb{C}^+)^2}
\overline{
\partial}_2\overline{\partial}_1 \Phi_{2,2}
(f) (z_1,\overline{z_2})\Theta_n(z_1,
\overline{z_2}) \ell_2(dz_1)
\ell_2(dz_2).
\end{eqnarray*}
Then $\Upsilon(f)$ is a distribution (in the sense of L. Schwartz) on
$C_c^{3,3}(\mathbb{R}^2)$. Moreover, $\Upsilon$~admits the following boundary
value representation:
%
\begin{eqnarray}
\Upsilon(f) &=&-\frac{1}{4\pi^2} \lim
_{\varepsilon\searrow0}
\int_{\mathbb{R}^2} f(x,y) \bigl\{ \Theta_n(x+\mathbf{i}
\varepsilon, y+\mathbf{i} \varepsilon) +\Theta_n(x-\mathbf{i}
\varepsilon, y-\mathbf{i}\varepsilon)
\nonumber
\\[-8pt]
\label{boundary-value-2}
\\[-8pt]
\nonumber
&&{}- \Theta_n(x-\mathbf{i}\varepsilon, y+\mathbf{i}\varepsilon) -
\Theta_n(x+\mathbf{i}\varepsilon, y-\mathbf{i}\varepsilon)\bigr\}\,dx
\,dy.
\end{eqnarray}
%
Notice that for every $f,g\in C_c^3(\mathbb{R})$ then $f\otimes g\in
C^{3,3}_c(\mathbb{R}^2)$ [where $(f\otimes g) (x,y)= f(x)g(y)$] and
\[
\Upsilon(f\otimes g) = \operatorname{cov} \bigl( Z^1_n(f),Z^1_n(g)
\bigr).
\]
%
\end{prop}

The proof of Proposition~\ref{propidentification-distribution} is
postponed to Section~\ref{proofidentification-distribution}.

\begin{rem}
By relying on Tillmann's results \cite{Tillmann-1953}, one
may prove that the support of $\Upsilon$ (as a distribution) is
included in ${\mathcal S}_n\times{\mathcal S}_n$. We provide a more
direct approach in a slightly simpler case in Section~\ref
{secexplicit-formulas}.
\end{rem}

\subsection{More covariance formulas}\label{secexplicit-formulas} We
provide here more explicit formulas for the variance than those given
in Theorem~\ref{thnon-analytic} and
Proposition~\ref{propidentification-distribution}; we also verify
that these formulas are in agreement with other formulas available in
the literature.

Recall that by \cite{silverstein-choi-1995}, Theorem~1.1, the limit
$\lim_{\varepsilon\searrow0} \tilde t_n(x+\mathbf{i}\varepsilon)$ denoted
by $\tilde t_n(x)$ exists for all $x\in\mathbb{R}$, $x\neq0$; the
same holds
true for $t_n$.

\begin{prop}\label{propnew-variance-formula}
Assume that  Assumptions~\ref{assX} and~\ref{assR} hold true and let $f,g\in C_c^3(\mathbb{R})$;
assume moreover for simplicity that ${\mathcal V}= \mathbb{E}X_{ij}^2$ is
either equal to 0 or 1 and that $R_n$ has real entries. Then the
covariance of $(Z_n(f),Z_n(g))$ in Theorem~\ref{thnon-analytic} writes
%
\begin{eqnarray}
&& \operatorname{cov}\bigl( Z^1_n(f),Z^1_n(g)
\bigr)= \frac{1+|{\mathcal V}|^2}{2\pi^2} \int_{{\mathcal S}_n^2} f'(x)g'(y)
\ln\biggl\llvert\frac{{\tilde t}_n(x)
-\overline{{\tilde t}_n(y)}}{
{\tilde t}_n(x) - {\tilde t}_n(y)}\biggr\rrvert\,dx \,dy
\nonumber
\\[-8pt]
\label{eqnew-formula-cov}
\\[-8pt]
\nonumber
&&\qquad{}+ \frac{\kappa}{\pi^2 n} \sum_{i=1}^N \biggl(
\int_{{\mathcal S}_n} f'(x) \operatorname{Im} \bigl( x
T_n(x) \bigr) _{ii} \,dx \biggr) \biggl( \int
_{{\mathcal
S}_n} g'(y) \operatorname{Im} \bigl( y
T_n(y) \bigr) _{ii} \,dy \biggr).
\end{eqnarray}
\end{prop}

The proof for Proposition~\ref{propnew-variance-formula} is postponed to
Section~\ref{proofnew-variance-formula}.

\begin{rem}
Notice that the first term in the RHS matches with the expression
provided in \cite{bai-silverstein-aop-2004}, equation (1.17) (see also
\cite{book-bai-silverstein}, equation (9.8.8)).
\end{rem}
%
\begin{rem}
Concerning the cumulant term, we shall compare it with the explicit
formula provided in \cite{lytova-pastur-aop-2009} (see also \cite
{book-pastur-shcherbina}) in the case where $R_n= I_N$. Recall that in
the context of Mar\v cenko--Pastur's theorem where $R_n=I_N$, we have
${\mathcal S}_n=[\lambda^-,\lambda^+]$ where $\lambda^-=(1-\sqrt{c_n})^2$,
$\lambda^+=(1+\sqrt{c}_n)^2$ and $(T_n(x) )_{ii}= t_n(x)$. We will
prove hereafter that
%
\begin{eqnarray}
&& \frac{\kappa c_n}{\pi^2} \biggl( \int_{\lambda
^-}^{\lambda^+}
f'(x) \operatorname{Im}\bigl\{ x t_n(x)\bigr\} \,dx \biggr)
\biggl( \int_{\lambda^-}^{\lambda^+} g'(y)
\operatorname{Im}\bigl\{ y t_n(y)\bigr\} \,dy \biggr)
\nonumber
\\
\label{eqcucum-MP}
&&\qquad= \frac{\kappa}{4c_n\pi^2} \biggl( \int_{\lambda^-}^{\lambda^+} f(x)
\frac{x - (1+c_n)}{\sqrt{(\lambda^+ - x)(x-\lambda^-)}} \,dx \biggr)
\\
\nonumber
&&\qquad\quad{}\times\biggl( \int_{\lambda^-}^{\lambda^+}
g(y) \frac{y - (1+c_n)}{\sqrt
{(\lambda^+ - y)(y-\lambda^-)}} \,dy \biggr).
\end{eqnarray}
Notice that the LHS of the equation above is the cumulant term as
provided in~\eqref{eqnew-formula-cov} if $R_n=I_N$ while the RHS is
the cumulant term as provided\footnote{Denote by the superscript
$^{\tt LP}$ the quantities in \cite{lytova-pastur-aop-2009} and use
the correspondence $c^{\tt LP}\leftrightarrow1/c$, $a^{\tt
LP}\leftrightarrow c$ and $\kappa_4^{\tt LP}\leftrightarrow(a^{\tt
LP})^4 \kappa= c^2\kappa$ to check that the RHS of \eqref{eqcucum-MP} equates the formula provided in \cite
{lytova-pastur-aop-2009}.} in \cite{lytova-pastur-aop-2009}.
\end{rem}
In the case where $R_n=I_N$, the Stieltjes transform of Mar\v
cenko--Pastur's distribution has an explicit form given by
(see, e.g., \cite{book-pastur-shcherbina}, Chapter~7)
\[
t_n(z) = \frac{1}{2c_n z} \bigl\{ \sqrt{ \bigl(z-(1+c_n)
\bigr)^2 - 4c_n } -\bigl(z-(1-c_n)\bigr) \bigr\},
\]
where the branch of the square root is fixed by its asymptotics:
$z-(1+c) +o(1)$ as $z\rightarrow\infty$. In particular, if $x\in
[\lambda^-,\lambda^+]$ then
\[
\sqrt{ \bigl(z-(1+c)\bigr)^2 - 4c } |_{z=x + \mathbf{i}0} = \mathbf{i}
\sqrt{ \bigl(\lambda^+-x\bigr) \bigl(x-\lambda^-\bigr)}.
\]
Hence,
\[
\operatorname{Im}\bigl\{ x t_n(x)\bigr\} = \frac{\sqrt{(\lambda^+ -
x)(x-\lambda
^-)}}{2c_n }.
\]
It remains to perform an integration by parts to get
\begin{eqnarray*}
\int_{\lambda^-}^{\lambda^+} f'(x) \operatorname{Im}
\bigl\{ x t_n(x)\bigr\} \,dx &=& - \int_{\lambda^-}^{\lambda^+}
f'(x) \frac{\sqrt{(\lambda^+ -
x)(x-\lambda^-)}}{2c_n } \,dx
\\
&=& \frac{1}{2c_n}\int_{\lambda^-}^{\lambda^+} f(x)
\frac{(1+c_n)-
x}{\sqrt{(\lambda^+ - x)(x-\lambda^-)}} \,dx
\end{eqnarray*}
which yields \eqref{eqcucum-MP}.

As a corollary of Proposition~\ref{propnew-variance-formula}, we
obtain the following extension of Theorem~\ref{thnon-analytic}.

Recall that ${\mathcal S}_n$ is the support of the probability measure
${\mathcal F}_n$. Due to Assumption~\ref{assR}, it is clear that
%
\begin{equation}
\label{defS-infini}
{\mathcal S}_n \subset{\mathcal S}_\infty
\stackrel{\triangle} {=} \bigl[0,\bolds{\lambda}^+_R \bigl( 1
+ \sqrt{\bolds{\ell}^+} \bigr)^2 \bigr],
\end{equation}
uniformly in $n$. Denote by $\mathbf{h}\in C_c^\infty(\mathbb
{R})$ a
function whose value is 1 on a $\eta$-neighborhood ${\mathcal
S}_\infty^\eta$ of ${\mathcal S}_\infty$.

\begin{coro}\label{corounbounded-f}
Assume that Assumptions \ref{assX}~and~\ref{assR} hold true and let $f_\ell\in C^3(\mathbb{R})$ with
$1\le
\ell\le k$; assume moreover that ${\mathcal V}= \mathbb{E}X_{ij}^2$
is either
equal to 0 or 1 and that $R_n$ has real entries. Let $\mathbf{h}\in
C_c^\infty(\mathbb{R})$ be as above.
Then \eqref{eqfluctuations-1}--\eqref{eqfluctuations-2} remain true
with $L_n^1({\mathbf f})$ replaced by
\[
L^{1,{\mathbf{h}}}_n({\mathbf f})= \bigl( \operatorname{tr}f_\ell
\bigl(\Sigma_n \Sigma_n^*\bigr) - \mathbb{E}
\operatorname{tr}(f_\ell\mathbf{h}) \bigl(\Sigma_n
\Sigma_n^*\bigr) ; 1\le\ell\le k \bigr)
\]
and with the Gaussian random vector $Z_n^1({\mathbf f} \mathbf{h})$ as
in Theorem~\ref{thnon-analytic}.
\end{coro}

The proof of Corollary~\ref{corounbounded-f} is postponed to Section~\ref
{secproof-coro-unbounded-f}.

\subsection{First-order expansions for the bias in the case of
nonanalytic functionals}\label{secbias-non-analytic}


\begin{theorem}\label{thnon-analytic-bias}
Assume Assumptions \ref{assX}~and~\ref{assR}  hold true and let $f\in C_c^{18}(\mathbb{R})$.
Denote by
%
\begin{equation}
\label{defbias-distribution}
Z_n^2(f)= \frac{1}\pi\operatorname{Re}
\int_{\mathbb{C}^+} \overline{\partial} \Phi_{17}(f) (z) {
\mathcal B}_n(z) \ell_2( dz),
\end{equation}
%
where ${\mathcal B_n}$ is defined in \eqref{defB}. Then
\[
\mathbb{E} \operatorname{Tr}(f) \bigl(\Sigma_n\Sigma_n^*
\bigr) - N\int f(\lambda) {\mathcal F}_n(d\lambda) -
Z_n^2(f) \mathop{\xrightarrow}_{N,n\to\infty} 0.
\]
%
\end{theorem}

The proof of Theorem~\ref{thnon-analytic-bias} is postponed to Section~\ref
{proofnon-analytic-II}.

\begin{rem}[(Why eighteen?)]\label{rem18}
A quick sketch of the proof
of Theorem~\ref{thnon-analytic-bias} provides some hints. Let $f$
have a bounded support. By Gaussian interpolation (whose cost is $f\in
C^8$), we only need to prove
\[
\mathbb{E} \operatorname{Tr}f\bigl(\Sigma^\mathbb{C}_n\bigl(
\Sigma_n^\mathbb{C}\bigr)^*\bigr) - N\int f(\lambda) {
\mathcal F}_n(d\lambda)\to0,
\]
where $\Sigma_n^\mathbb{C}$ is the counterpart of $\Sigma_n$ with
${\mathcal
N}_\mathbb{C}(0,1)$ i.i.d. entries. The proof of the latter is based on
Helffer--Sj\"{o}strand's formula
\[
\mathbb{E} \operatorname{Tr}f\bigl(\Sigma^\mathbb{C}_n\bigl(
\Sigma_n^\mathbb{C}\bigr)^*\bigr) - N\int f(\lambda) {
\mathcal F}_n(d\lambda) =\frac{1}\pi\operatorname{Re}\int
_{\mathbb{C}^+} \overline{\partial} \Phi_k(f) \bigl\{
\operatorname{Tr}\mathbb{E} Q_n^\mathbb{C}- N t_n
\bigr\}\,d\ell_2,
\]
where $Q_n^\mathbb{C}=(\Sigma^\mathbb{C}_n(\Sigma_n^\mathbb
{C})^*-zI_N)^{-1}$, and on the
following estimate, stated in Proposition~\ref{propbias-hard}:
%
\begin{equation}
\label{eqestimate-singularity}\hspace*{4pt}
\bigl\llvert\mathbb{E}\operatorname
{Tr}\bigl(
\Sigma^\mathbb{C}_n \bigl(\Sigma_n^\mathbb{C}
\bigr)^* - zI_N\bigr)^{-1} - N t_n(z) \bigr
\rrvert\le\frac{1}{n} P_{12}\bigl(|z|\bigr)P_{17}\bigl(\bigl|
\operatorname{Im}(z)\bigr|^{-1}\bigr),
\end{equation}
where $P_k$ denotes a polynomial with degree $k$ and positive
coefficients. In view of Proposition~\ref{propcompensation},
$f$ needs to be of class $C^{18}$. If one can improve estimate \eqref
{eqestimate-singularity} and decrease the powers of $|\operatorname
{Im}(z)|^{-1}$,
then one will automatically lower the regularity assumption over $f$.
Notice that in the case of the Gaussian unitary ensemble, counterpart
of \eqref{eqestimate-singularity} features $|\operatorname{Im}(z)|^{-7}$
on its RHS
(cf. \cite{haagerup-thorbjornsen-2005}, Lemma~6.1), hence the needed
regularity is $f\in C^8$ in this case.
\end{rem}

\begin{prop}\label{propbias-as-distribution}
Let $Z_n^2(f)$ be defined as in \eqref{defbias-distribution}, then
$Z_n^2$ is a distribution (in the sense of L. Schwartz) on
$C^{18}_c(\mathbb{R}
)$ and
%
\begin{equation}
\label{eqboundary-bias}
Z_n^2(f)= \frac{-\mathbf{i}}{2\pi} \lim
_{\varepsilon\searrow0} \int_{\mathbb{R}
} f(x) \bigl\{ {\mathcal
B}_n(x+\mathbf{i}\varepsilon) - {\mathcal B}_n(x-
\mathbf{i}\varepsilon) \bigr\} \,dx.
\end{equation}
Moreover, the singular points of ${\mathcal B}_n(z)$ are included in
${\mathcal S}_n$ and so is the support of $Z_n^2$ (as a distribution).
In particular, one can extend $Z_n^2$ to $C^{18}(\mathbb{R})$ by
\[
\check{Z}_n^2(f) = Z_n^2(f
\mathbf{h}), \qquad f\in C^{18}(\mathbb{R}),
\]
where $\check{Z}_n^2$ is the extension to $C^{18}(\mathbb{R})$ and
$\mathbf{h}\in C^{\infty}_c(\mathbb{R})$ has value 1 on
${\mathcal S}_n$.
\end{prop}
The proof of Proposition~\ref{propbias-as-distribution} is postponed to
Section~\ref{proofprop-bias-distribution}.
%
\begin{coro}\label{corobias-unbounded}
Assume Assumptions \ref{assX}~and~\ref{assR} hold true. Let $f\in C^{18}(\mathbb{R})$
and $\mathbf{h}\in C_c^{\infty}(\mathbb{R})$ be a function whose value
is 1 on a neighborhood of ${\mathcal S}_\infty$,
then the following convergence holds true:
\[
\mathbb{E} \operatorname{Tr}(f\mathbf{h}) \bigl(\Sigma_n
\Sigma_n^*\bigr) - N\int f(\lambda) {\mathcal F}_n(d
\lambda) - \check{Z}_n^2(f) \mathop{\xrightarrow}_{N,n\to\infty}
0.
\]
\end{coro}
The proof is straighforward and is therefore omitted.

\section{Proof of Theorem~\texorpdfstring{\protect\ref{lemmamain}}{1} (CLT for the trace of the
resolvent)}\label{secproof}
Recall that $M_n(z) = \operatorname{tr}Q_n(z) - Nt_n(z)$.
It will be convenient to decompose $M_n(z)$ as
%
\begin{equation}
\label{eqdecomposition-Mn}
\hspace*{12pt} M_n(z)= M_n^1(z) +
M_n^2(z)\qquad \mbox{where } \cases{
\displaystyle M_n^1(z) = \operatorname{tr}Q_n(z) -
\operatorname{tr}\mathbb{E} Q_n(z),
\vspace*{3pt}\cr
\displaystyle M_n^2(z) = N \bigl( \mathbb{E} f_n(z) -
t_n(z) \bigr).}
\end{equation}
Denote by $\xi_j$ the $N\times1$ vector
\[
\xi_j= \Sigma_{\cdot j} = \frac{1}{\sqrt{n}} R^{1/2}
X_{\cdot j}
\]
and by $\mathbb{E}_j$ the conditional expectation
with respect to ${\mathcal G}_j$, the $\sigma$-field generated by $\xi
_1,\ldots,\xi_j$:
%
\begin{equation}
\label{defcond-expectation}
\mathbb{E}_j = \mathbb{E} ( \cdot\vert
{\mathcal
G}_j).
\end{equation}
By convention, $\mathbb{E}_0=\mathbb{E}$.
We split Theorem~\ref{lemmamain} into intermediate results. Recall
the definitions of $D_{\varepsilon}, D^+$ and $D$ in \eqref
{defsets-D}. Let
\[
\Gamma= D^+\cup\overline{D^+} \qquad\mbox{where } \overline{D^+}=\bigl\{
\bar{z}, z
\in D^+ \bigr\}.
\]

\begin{prop}\label{propconvergence-martingale}
Assume that  Assumptions \ref{assX}~and~\ref{assR} hold true; let $z_1,z_2 \in\Gamma$, then\vspace*{-6pt}
\[
M_n^{ 1}(z) = \sum_{j=1}^n
Z_j^n(z) + o_P(1),
\]
where the $Z_j^n$'s are martingale increments with respect to the
$\sigma$-field ${\mathcal G}_i$\vspace*{-3pt} and
%
\begin{eqnarray}\label{eqconv-martingale}
\sum_{j=1}^n \mathbb{E}_{j-1}
Z_j^n(z_1) Z_j^n(z_2)
- \Theta_n(z_1,z_2) &\displaystyle \mathop{
\xrightarrow}_{N,n\rightarrow\infty}^{\mathcal P}& 0,
\\[-2pt]
\label{eqconv-martingale-conjugate}
\sum_{j=1}^n \mathbb{E}_{j-1}
Z_j^n(z_1)\overline{Z_j^n(z_2)}
- \Theta_n(z_1,\overline{z_2}) & \displaystyle\mathop{
\xrightarrow}_{N,n\rightarrow\infty}^{\mathcal P}& 0,
\end{eqnarray}
where $\Theta_n$ is defined in \eqref{eqdecomposition-covariance}. Moreover,
\[
M_n^{ 2} (z) - {\mathcal B}_n(z) \mathop{
\xrightarrow}_{N,n\rightarrow\infty}0,
\]
where ${\mathcal B}_n$ is defined in \eqref{defB}.
\end{prop}

\begin{prop}\label{propgaussian-process}
$\!$There exists a sequence $(G_n(z), z\in\Gamma)$ of
two-dimen\-sional Gaussian processes with mean
$
\mathbb{E}G_n(z) = {\mathcal B}_n(z)
$
and covariance
\begin{eqnarray*}
\operatorname{cov} \bigl( G_n(z_1),
G_n(z_2) \bigr) &=& \mathbb{E} \bigl(G_n(z_1)
- \mathbb{E}G_n(z_1) \bigr) \bigl(G_n(z_2)
- \mathbb{E}G_n(z_2) \bigr)
\\[-2pt]
&=&\Theta_n(z_1,z_2).
\end{eqnarray*}
Moreover, $(G_n(z), z\in D_\varepsilon)$ is tight.
\end{prop}

\subsection{Proof for Proposition~\texorpdfstring{\protect\ref{propconvergence-martingale}}{5.1}}
The fact that $(M_n)$ is a tight sequence has already been established
in \cite{bai-silverstein-aop-2004} (regardless of the assumption
$\kappa=0$ and $|{\mathcal V}|=0/1$).
In order to proceed, we shall
rely heavily
on the proof of \cite
{bai-silverstein-aop-2004}, Lemma~1.1, which is the crux of Bai and
Silverstein's paper. In Section~\ref{secreview}, we recall the main
steps of Bai and Silverstein's computations of the variance/covariance.
In Sections~\ref{sectionvariance-vartheta} and \ref
{secdiagonal-terms-variance}, we compute the extra terms in the
limiting variance. In Section~\ref{secbias}, we compute the limiting
bias. 
In Section~\ref{secproof-lemma}, we finally conclude the proof of
Theorem~\ref{lemmamain} and address various subtleties which appear
due to the existence of a sequence of Gaussian limiting processes.

In the sequel, we shall drop subscript $n$ and write $Q$ and $R$
instead of $Q_n$ and $R_n$.
Denote by $Q_j(z)$ the resolvent of matrix $\Sigma\Sigma^* - \xi_j
\xi_j^*$,\vspace*{-3pt} that is,
\[
Q_j(z)= \bigl( -zI + \Sigma\Sigma^* - \xi_j
\xi_j^* \bigr)^{-1}.
\]
The following quantities will be\vspace*{-2pt} needed:
\begin{eqnarray*}
\beta_j(z)&=& \frac{1}{1+\xi_j^* Q_j(z) \xi_j},
\\[-2pt]
\bar{\beta}_j(z) &=& \frac{1}{1+({1}/n) \operatorname{tr}R_n
Q_j(z)},
\\
b_n(z)&=& \frac{1}{1 +({1}/n) \mathbb{E} \operatorname{tr}R_n
Q_1(z)},
\\
\varepsilon_j(z)&=& \xi_j^* Q_j(z)
\xi_j - ({1}/n) \operatorname{tr}R_nQ_j(z),
\\
\delta_j(z)&=&\xi_j^* Q_j^2(z)
\xi_j - \frac{1}N \operatorname{tr}R_n
Q_j^2(z) = \frac{d}{dz} \varepsilon_j(z).
\end{eqnarray*}

\subsubsection{Preliminary variance computations} \label{secreview}
We briefly review in this section the main steps related to the
computation of the variance/covariance as presented in \cite
{bai-silverstein-aop-2004}. These standard steps will finally lead to
equation \eqref{eqdecomp-prelim} which will be the starting point of
the computations associated to the $|{\mathcal V}|^2$- and $\kappa$-terms of the variance.

Let $z\in\Gamma$:
\begin{eqnarray*}
N \bigl( f_n(z) - \mathbb{E} f_n(z) \bigr) &=& -\sum
_{j=1}^n ( \mathbb{E}_j -
\mathbb{E}_{j-1} ) \beta_j(z) \xi_j^*
Q_j^2(z) \xi_j
\\
&=& - \sum_{j=1}^n \mathbb{E}_j
\biggl( \bar{\beta}_j(z) \delta_j(z) -\bar{
\beta}_j^2(z)\varepsilon_j(z)
\frac{1}n \operatorname{tr}R Q_j^2 \biggr) +
o_P(1),
\end{eqnarray*}
where $\mathbb{E}_j$ is introduced in \eqref{defcond-expectation}.
Denote by
\[
Z^n_j(z)= -\mathbb{E}_j \biggl( \bar{
\beta}_j(z) \delta_j(z) - \bar{\beta}_j^2(z)
\varepsilon_j(z) \frac{1}n \operatorname{tr}R
Q_j^2(z) \biggr) = - \mathbb{E}_j
\frac{d}{dz} \bigl( \bar{\beta}_j(z) \varepsilon_j(z)
\bigr).
\]
Hence,
\begin{eqnarray*}
M^1_n(z)&=& N \bigl( f_n(z) - \mathbb{E}
f_n(z) \bigr) = \sum_{j=1}^n
Z^n_j(z) + o_P(1).
\end{eqnarray*}
%
The RHS appears as a sum of martingale increments. Such a decomposition
is important since it will enable us to rely on powerful CLTs for
martingales (see \cite{book-billinsley-proba-and-measure-95},
Theorem~35.12, and the variations below
in Lemmas \ref{lemmavariations-martingales} and \ref
{lemmamartingales-multidim}). These CLTs rely on the study of the terms
\[
\sum_{j=1}^n \mathbb{E}_{j-1}
Z^n_j(z_1) Z^n_j(z_2)\quad
\mbox{and}\quad \sum_{j=1}^n
\mathbb{E}_{j-1} Z^n_j(z_1)
\overline{Z^n_j(z_2)}.
\]
Notice that since $\overline{Z^n_j(z)}=Z_j^n(\bar{z})$, we have $
\mathbb{E}_{j-1} Z^n_j(z_1) \overline{Z^n_j(z_2)}= \mathbb{E}_{j-1}
Z^n_j(z_1)\times  Z^n_j(\overline{z_2})$.
Since the set $\Gamma$ is stable by complex conjugation, it is
sufficient to study the limiting behavior of
\[
\sum_{j=1}^n \mathbb{E}_{j-1}
Z^n_j(z_1) Z^n_j(z_2),\qquad  z_1,z_2\in\Gamma
\]
in order to prove \eqref{eqconv-martingale} and \eqref
{eqconv-martingale-conjugate}. Now,
%
\begin{eqnarray}
&& \sum_{j=1}^n
\mathbb{E}_{j-1} Z^n_j(z_1)
Z^n_j(z_2)
\nonumber
\\[-8pt]
\label{eqdecomposition-importante}
\\[-8pt]
\nonumber
&&\qquad= \frac{\partial
^2}{\partial z_1 \partial z_2} \Biggl\{
\sum_{j=1}^n\mathbb{E}_{j-1}
\bigl[ \mathbb{E}_j \bigl( \bar{\beta}_j(z_1)
\varepsilon_j(z_1) \bigr) \mathbb{E}_j
\bigl( \bar{\beta}_j(z_2) \varepsilon_j(z_2)
\bigr) \bigr] \Biggr\}.
\end{eqnarray}
Following the same arguments as in \cite{bai-silverstein-aop-2004},
page~571, one can prove that it is sufficient to
study the convergence in probability of
\[
\sum_{j=1}^n\mathbb{E}_{j-1}
\bigl[ \mathbb{E}_j \bigl( \bar{\beta}_j(z_1)
\varepsilon_j(z_1) \bigr) \mathbb{E}_j
\bigl( \bar{\beta}_j(z_2) \varepsilon_j(z_2)
\bigr) \bigr].
\]
Moreover,
%
\begin{eqnarray}
&& \sum_{j=1}^n
\mathbb{E}_{j-1} \bigl[ \mathbb{E}_j \bigl( \bar{\beta
}_j(z_1) \varepsilon_j(z_1)
\bigr) \mathbb{E}_j \bigl( \bar{\beta}_j(z_2)
\varepsilon_j(z_2) \bigr) \bigr]
\nonumber
\\
\label{eqdecomposition-1}
&&\qquad= \sum_{j=1}^n b_n(z_1)
b_n(z_2) \mathbb{E}_{j-1} \bigl[ \mathbb
{E}_j \varepsilon_j(z_1)
\mathbb{E}_j \varepsilon_j(z_2) \bigr] +
o_P(1)
\\
&&\qquad= \sum_{j=1}^n z_1\tilde
t_n(z_1) z_2 \tilde t_n(z_2)
\mathbb{E}_{j-1} \bigl[ \mathbb{E}_j \varepsilon_j(z_1)
\mathbb{E}_j \varepsilon_j(z_2) \bigr] +
o_P(1).
\nonumber
\end{eqnarray}

Hence, it is finally sufficient to study the limiting behavior (in
terms of convergence in probability) of the quantity
%
\begin{equation}
\label{eqlimiting-covariance}
\sum_{j=1}^n
\mathbb{E}_{j-1} \bigl( \mathbb{E}_j \varepsilon
_j(z_1) \mathbb{E}_j \varepsilon_j
(z_2) \bigr),\qquad  z_1,z_2\in\Gamma.
\end{equation}

Denote by $A^T$ the transpose matrix of $A$. Applying \eqref
{eqcovariance-identity} yields
%
\begin{eqnarray}
&& \sum_{j=1}^n
\mathbb{E}_{j-1} \bigl( \mathbb{E}_j \varepsilon_j(z_1)
\mathbb{E}_j \varepsilon_j (z_2) \bigr)
\nonumber
\\
&&\qquad =
\frac{1}{n^2} \sum_{j=1}^n
\operatorname{tr} \bigl( R^{1/2} \mathbb{E}_j
Q_j (z_1) R \mathbb{E}_j
Q_j(z_2) R^{1/2} \bigr)
\nonumber
\\[-8pt]
\label{eqdecomp-prelim}
\\[-8pt]
\nonumber
&&\qquad\quad {}+\frac{|\mathcal{V}|^2}{n^2} \sum_{j=1}^n
\operatorname{tr} \bigl( R^{1/2} \mathbb{E}_j
Q_j (z_1) R^{1/2} \bigl( R^{1/2}
\mathbb{E}_j Q_j(z_2) R^{1/2}
\bigr)^T \bigr)
\\
&&\qquad\quad {}+ \frac{\kappa}{n^2} \sum_{j=1}^n
\sum_{i=1}^N \bigl( R^{1/2}
\mathbb{E}_j Q_j(z_1) R^{1/2}
\bigr)_{ii} \bigl( R^{1/2} \mathbb{E}_j
Q_j(z_2) R^{1/2} \bigr)_{ii}.
\nonumber
\end{eqnarray}
The limiting behavior of the first term of the RHS has been completely
described in \cite{bai-silverstein-aop-2004} where it has been shown that
%
\begin{eqnarray}
&& \frac{\partial^2}{\partial z_1 \partial z_2} \Biggl\{
z_1 z_2
\tilde t_n(z_1) \tilde t_n(z_2)
\frac{1}{n^2} \sum_{j=1}^n
\operatorname{tr} \bigl( R^{1/2} \mathbb{E}_j
Q_j (z_1) R \mathbb{E}_j
Q_j(z_2) R^{1/2} \bigr) \Biggr\}
\nonumber
\\[-8pt]
\label{eqtoto1}
\\[-8pt]
&&\qquad = \Theta_{0,n}(z_1,z_2) + o_P(1),
\nonumber
\end{eqnarray}
with $\Theta_{0,n}(z_1,z_2)$ defined in \eqref{deftheta0}.

We shall focus on the second and third terms.

\subsubsection{The term proportional to $|\mathcal{V}|^2$ in the
variance}\label{sectionvariance-vartheta} Notice first that the value
of $t_n$ and $\tilde t_n$ is the same whether $R$ is replaced by
$\overline
{R}$ in \eqref{eqdeteq-ST} and \eqref{eqdef-t-tilde} since $t_n$
and $\tilde t_n$ only depend on the spectrum of $R$ (which is the same
as the spectrum of $\overline{R}$). Notice also that
$(R^{1/2})^T=\overline
{R}^{1/2}$, hence
\[
\bigl( R^{1/2} \mathbb{E}_j Q_j(z_2)
R^{1/2} \bigr)^T = \overline{R}^{1/2}
\mathbb{E}_j Q^T_j(z_2)
\overline{R}^{1/2}.
\]
Recall the definition of $T^T_n(z)$ given by \eqref
{eqdeteq-RESOLVENT-conjugate}. Taking into account the fact that for a
deterministic matrix $A$,
%
\begin{equation}
\label{eqmixed-quadratic-forms}
\hspace*{10pt}\mathbb{E}\xi_j^T A \xi_j =
\frac{\mathcal{V}}n \operatorname{tr}\bigl(\overline{R}^{1/2} A
R^{1/2}\bigr) \quad \mbox{and}\quad \mathbb{E}\xi_j^*A \bar{
\xi}_j = \frac
{\overline{\mathcal{V}
}}n \operatorname{tr}\bigl(R^{1/2} A
\overline{R}^{1/2}\bigr),
\end{equation}
and following closely \cite{bai-silverstein-aop-2004}, Section~2, it
is a matter of bookkeeping
to establish that
%
\begin{eqnarray}
&& \frac{|\mathcal{V}|^2 z_1 z_2}{n^2} \tilde t_n(z_1) \tilde
t_n(z_2) \sum_{j=1}^n
\operatorname{tr} \bigl( R^{1/2} \mathbb{E}_j
Q_j (z_1) R^{1/2} \bigl( R^{1/2}
\mathbb{E}_j Q_j(z_2) R^{1/2}
\bigr)^T \bigr)
\nonumber
\\
\label{eqvar-moment2}
&&\qquad= |\mathcal{V}|^2 {\mathcal A}_n(z_1,z_2)
\times\frac{1}n \sum_{j=1}^n
\frac{1}{ 1 - (({j-1})/n) |\mathcal{V}|^2 {\mathcal
A}_n(z_1,z_2)} + o_P(1)
\\
&&\qquad= \int_0^{|\mathcal{V}|^2 {\mathcal A}_n(z_1,z_2)} \frac{dz}{1-z} +
o_P(1),
\nonumber
\end{eqnarray}
where
\[
{\mathcal A}_n(z_1,z_2) =
\frac{z_1 z_2}n \tilde t_n(z_1) \tilde
t_n(z_2) \operatorname{tr} \bigl\{ R^{1/2}
T_n(z_1) R^{1/2} \overline{R}^{1/2}
T^T_n(z_2) \overline{R}^{1/2} \bigr\}.
\]
Finally,
\begin{eqnarray}
\frac{\partial^2}{\partial z_1 \partial z_2} \eqref
{eqvar-moment2} &=&|{\mathcal
V}|^2 \Theta_{1,n}(z_1,z_2) +
o_P(1)
\nonumber
\\[-8pt]
\label{eqtoto2}
\\[-8pt]
\nonumber
& =& |\mathcal{V}|^2 \frac{\partial}{\partial z_2} \biggl\{ \frac
{\partial{\mathcal A}_n(z_1,z_2)/\partial z_1}{1-|\mathcal{V}
|^2{\mathcal A}_n(z_1,z_2)}
\biggr\} + o_P(1).
\nonumber
\end{eqnarray}
%

\subsubsection{The cumulant term in the variance} \label
{secdiagonal-terms-variance}

We now handle the term proportional to $\kappa$ in \eqref{eqdecomp-prelim}:
%
\begin{equation}
\label{eqcumul-term}
\frac{1}{n^2} \sum_{j=1}^n
\sum_{i=1}^N \bigl( R^{1/2}
\mathbb{E}_j Q_j(z_1) R^{1/2}
\bigr)_{ii} \bigl( R^{1/2} \mathbb{E}_j
Q_j(z_2) R^{1/2} \bigr)_{ii}.
\end{equation}
The objective is to prove that $\mathbb{E}_j Q_j(z)$ can be replaced by
$T_n(z)$ in the formula above, which boils down to prove a convergence
of quadratic forms of the type~\eqref{eqconv-quadra}. Such a
convergence has already been established in \cite
{hachem-et-al-aihp-2013} for large covariance matrices based on a
noncentered matrix model with separable variance profile.

Let $\bolds{\delta}_z$ be the distance between the point $z\in
\mathbb{C}$
and the real nonnegative axis $\mathbb{R}^+$:
%
\begin{equation}
\label{defdzz}
\bolds{\delta}_z = \operatorname{dist}\bigl(z,
\mathbb{R}^+\bigr).
\end{equation}

\begin{prop} \label{propvariance-forme-quadra}
Assume that Assumptions \ref{assX} and \ref{assR}  hold true and let $u_n$ be a
deterministic $N\times1$ vector, then
\[
\mathbb{E}\bigl\llvert u_n^* Q(z) u_n -
u_n^* \mathbb{E}Q(z) u_n \bigr\rrvert^2 \le
\frac{1}n \Phi\bigl(|z|\bigr)\Psi\biggl( \frac{1}{\bolds{\delta}_z} \biggr)\|
u_n\|^2,
\]
where $\Phi$ and $\Psi$ are fixed polynomials with coefficients
independent from $N,n,z$ and $(u_n)$.
\end{prop}

Proof of Proposition~\ref{propvariance-forme-quadra} is an easy adaptation
of \cite{hachem-et-al-aihp-2013}, Proposition~2.7; see also the proof
of Proposition~\ref{propvariance-estimate-sharp} below. It is
therefore omitted.

\begin{prop}\label{propcumulant} Assume that Assumptions \ref{assX} and
\ref{assR}  hold true, then the following convergence holds true:
\begin{eqnarray*}
&& \frac{1}{n^2} \sum_{j=1}^n \sum
_{i=1}^N \bigl( R^{1/2}
\mathbb{E}_j Q_j(z_1) R^{1/2}
\bigr)_{ii} \bigl( R^{1/2} \mathbb{E}_j
Q_j(z_2) R^{1/2} \bigr)_{ii}
\\
&&\qquad{}- \frac{1}n \sum_{i=1}^N
\bigl( R^{1/2} T(z_1) R^{1/2}
\bigr)_{ii} \bigl( R^{1/2} T(z_2) R^{1/2}
\bigr)_{ii} \mathop{\xrightarrow}_{n,N\rightarrow\infty}^{\mathcal P} 0.
\end{eqnarray*}
\end{prop}

The proof below has been suggested by a referee whom we thank; it
substantially simplifies the initial one.

\begin{pf*}{Proof of Proposition~\protect\ref{propcumulant}}
We first transform the sum to be calculated:
%
\begin{equation}
\label{eqsum-cumulant-term}
\frac{1}{n^2} \sum_{j=1}^n
\sum_{i=1}^N \bigl( R^{1/2}
\mathbb{E}_j Q_j(z_1) R^{1/2}
\bigr)_{ii} \bigl( R^{1/2} \mathbb{E}_j
Q_j(z_2) R^{1/2} \bigr)_{ii}.
\end{equation}
Using Proposition~\ref{propvariance-forme-quadra} enables us to
replace the conditional expectation $\mathbb{E}_i$ by the true
expectation in every term $ (R^{1/2} \mathbb{E}_j Q_j(z)
R^{1/2}
)_{ii}$. Now using the fact that
%
\begin{equation}
\label{eqrank-one-perturbation-bis}
Q = Q_j - \frac{Q_j \xi_j \xi_j^*
Q_j}{1+\xi_j^* Q_j \xi_j},
\end{equation}
%
one can replace $\mathbb{E}Q_j$ by $\mathbb{E}Q$. We now prove the
following estimate:
%
\begin{equation}
\label{eqreviewer}
\bigl\llvert\mathbb{E}u^* Q(z) v - u^* T(z) v\bigr
\rrvert\le
\frac{C}{\sqrt
{n}
\operatorname{Im}^k(z)} \| u\| \|v\|,
\end{equation}
where neither $K$ nor $k$ depend on $N,n$. Clearly, Proposition~\ref
{propcumulant} follows\break from~\eqref{eqreviewer}.

Using \eqref{eqrank-one-perturbation-bis} and the associated fact
that $(Q(z)\xi_j)_i=\beta_j(z) (Q_j(z)\xi_j)_i$, we get
%
\begin{eqnarray}
\mathbb{E}u^* Q(z) \Sigma\Sigma^* T(z) v &=& \sum_j
\mathbb{E}\beta_j \bigl( u^* Q_j(z) \xi_j
\xi_j^* T(z)v \bigr)
\nonumber
\\
&\stackrel{(a)}=& -z\tilde t_n(z) \sum_j
\mathbb{E} \bigl( u^* Q_j(z) \xi_j
\xi_j^* T(z)v \bigr)+ O \biggl( \frac{1}{\sqrt{n} \operatorname{Im}
^k(z)} \biggr)
\nonumber
\\[-8pt]
\label{eqreviewer-1}
\\[-8pt]
\nonumber
&=& -z\frac{\tilde t_n(z)}n \sum_j \mathbb{E}
\bigl( u^* Q_j(z) R T(z)v \bigr)+ O \biggl( \frac{1}{\sqrt{n}
\operatorname{Im}^k(z)}
\biggr)
\\
&\stackrel{(b)}=& -z\tilde t_n(z) \mathbb{E}u^* Q(z) RT(z) v + O
\biggl( \frac
{\|u\| \|v\|}{\sqrt{n} \operatorname{Im}^k(z)} \biggr),
\nonumber
\end{eqnarray}
where we used that $\mathbb{E}| \beta_j(z) + z\tilde t_n(z)|^2 \le K n^{-1}
|\operatorname{Im}^{-k}(z)|$ (see, e.g., \cite{bai-silverstein-aop-2004})
to prove $(a)$ and we used \eqref{eqrank-one-perturbation-bis} to
replace $Q_j$ by $Q$ in $(b)$.

On the other hand,
\begin{eqnarray}
 \mathbb{E}u^* Q(z) \Sigma\Sigma^* T(z) v &=&
\mathbb{E}u^* Q(z)
\bigl( \Sigma\Sigma^* -zI_N + zI_N \bigr) T(z) v
\nonumber
\\[-8pt]
\label{eqreviewer-2}
\\[-8pt]
\nonumber
&=& u^* T(z) v + z\mathbb{E}u^* Q(z) T(z) v.
\nonumber
\end{eqnarray}
Taking into account that by Definitions \eqref{eqdeteq-RESOLVENT} and
\eqref{defttilde}
\[
T(z) = \bigl(-z\tilde t_n(z) c_n R -zI_N
\bigr)^{-1},
\]
we get
\begin{eqnarray*}
&& u^* T(z) v - u^* \mathbb{E}Q(z) v \\
&&\qquad= u^* T(z) v - \mathbb{E}u^* Q(z)
\bigl(-z
\tilde t_n(z) c_n R -zI_N\bigr)T(z) v
\\
&&\qquad= u^* T(z) v + z\tilde t_n(z) c_n \mathbb{E}u^* Q(z) R
T(z) v + z \mathbb{E}u^* Q(z) T(z) v
\\
&&\qquad\stackrel{(a)}=  \mathbb{E}u^* Q(z) \Sigma\Sigma^* T(z) v + z\tilde
t_n(z) c_n \mathbb{E}u^* Q(z) R T(z) v
\\
&&\qquad\stackrel{(b)}= O \biggl( \frac{\|u\| \|v\|}{\sqrt{n} \operatorname{Im}
^k(z)} \biggr),
\end{eqnarray*}
where $(a)$ follows from \eqref{eqreviewer-2} and $(b)$ from \eqref
{eqreviewer-1}.
This completes the proof of~\eqref{eqreviewer}, hence the proof of
Proposition~\ref{propcumulant}.
\end{pf*}

Combining the result in Proposition~\ref{propcumulant} together with
\eqref{eqdecomposition-1} and \eqref{eqdecomp-prelim}, we have
proved so far that
%
\begin{eqnarray}
&& \frac{\partial^2}{\partial z_1
\partial z_2} \Biggl\{ \frac{z_1 z_2 \tilde t_n(z_1) \tilde t_n(z_2)}{
n^2}\nonumber \\
&&\quad{}\times\sum
_{j=1}^n \sum_{i=1}^N
\bigl( R^{1/2} \mathbb{E}_j Q_j(z_1)
R^{1/2} \bigr)_{ii} \bigl( R^{1/2}
\mathbb{E}_j Q_j(z_2) R^{1/2}
\bigr)_{ii} \Biggr\}
\nonumber
\\[-8pt]
\label{alternative-formula-cov2}
\\[-8pt]
\nonumber
&&\qquad= \frac{1}n \sum_{i=1}^N
\frac{ \partial^2}{\partial z_1\partial z_2} \bigl\{z_1 z_2 \tilde
t_n(z_1) \tilde t_n(z_2) \bigl(
R^{1/2}_n T_n(z_1)
R^{1/2}_n \bigr)_{ii} \bigl( R^{1/2}_n
T_n(z_2) R^{1/2}_n
\bigr)_{ii} \bigr\}\\
\nonumber
&&\qquad\quad {}+ o_P(1).
\end{eqnarray}

Taking into account \eqref{eqdef-t-tilde} and the matrix identity $
U(I+VU)^{-1} V = 1- (I+UV)^{-1}$, we obtain
%
\begin{eqnarray}
\eqref{alternative-formula-cov2} 
&=&
\frac{1}n \sum_{i=1}^N
\frac{ \partial^2}{\partial z_1\partial
z_2} \bigl( I_N - \bigl(I_N + \tilde
t_n(z_1)R_n\bigr)^{-1}
\bigr)_{ii} \bigl( I_N - \bigl(I_N + \tilde
t_n(z_2)R_n\bigr)^{-1}
\bigr)_{ii} \nonumber\\
&&{}+ o_P(1)
\nonumber
\\
\label{alternative-formula-cov2-bis}
&=& \frac{1}n \sum_{i=1}^N
\frac{ \partial}{\partial z_1} \bigl[ z_1 T_n(z_1)
\bigr]_{ii} \frac{ \partial}{\partial z_2} \bigl[z_2
T_n(z_2) \bigr]_{ii} + o_P(1)
\\
\nonumber
&=& \frac{ z_1^2 z_2^2 \tilde t'_n(z_1) \tilde t'_n(z_2)}n \sum_{i=1}^N
\bigl( R_n^{1/2} T^2_n(z_1)
R_n^{1/2} \bigr)_{ii} \bigl( R_n^{1/2}
T^2_n(z_2) R_n^{1/2}
\bigr)_{ii} + o_P(1)
\\
&=& \Theta_{2,n}(z_1,z_2) +
o_P(1),
\nonumber
\end{eqnarray}
where $\Theta_{2,n}$ is given by formula \eqref{deftheta2}.

Now gathering \eqref{eqtoto1}, \eqref{eqtoto2} and \eqref
{alternative-formula-cov2-bis}, we have established so far:
\[
\sum_{j=1}^n \mathbb{E}_{j-1}
Z_j^n(z_1) Z_j^n(z_2)
= \Theta_n(z_1,z_2) +o_P(1)
\]
which is the first part of Proposition~\ref{propconvergence-martingale}.

\subsubsection{Computations for the bias}\label{secbias}
In this
section, we are interested in the computation of $N(\mathbb{E}f_n(z) -
t_n(z))$. As
\[
\tilde f_n(z) = -\frac{(1-c_n)} z + c_n
f_n(z)\quad\mbox{and} \quad\tilde t_n(z) = -
\frac{(1-c_n)} z + c_n t_n(z),
\]
we immediately obtain $N(\mathbb{E}f_n(z) - t_n(z))= n(\mathbb
{E}\tilde f_n(z) -
\tilde t_n(z))$. Combining \eqref{defttilde} and \eqref
{eqdef-t-tilde} yields
%
\begin{equation}
\label{eqequation-t-tilde}
-z -\frac{1}{\tilde t_n(z)} + \frac{1}n
\operatorname{tr}R_n
\bigl( I_N + \tilde t_n(z) R_n
\bigr)^{-1} = 0.
\end{equation}
Following Bai and Silverstein \cite{bai-silverstein-aop-2004},
Section~4, we introduce the quantity $A_n(z)$
defined as
\begin{eqnarray*}
A_n(z) &=& z\mathbb{E}\tilde f_n(z) +1+
\frac{1}n \operatorname{tr} \bigl( I_N + \mathbb{E}\tilde
f_n(z)R_n \bigr)^{-1} -c_n
\\
&=& z\mathbb{E}\tilde f_n(z) +1 + \frac{1}n
\operatorname{tr} \bigl( I_N + \mathbb{E}\tilde f_n(z)R_n
\bigr)^{-1} -\frac{1}n \operatorname{tr}I_N^{-1}
\\
&=&- \mathbb{E}\tilde f_n(z) \biggl( - z -\frac{1}{\mathbb{E}\tilde
f_n(z)}+
\frac{1}n \operatorname{tr}R_n\bigl( I_N +
\mathbb{E}\tilde f_n(z) R_n\bigr)^{-1} \biggr),
\end{eqnarray*}
hence
%
\begin{equation}
\label{eqequation-f-tilde}
-\frac{A_n(z)}{\mathbb{E}\tilde f_n(z)} = -z
-\frac{1}{\mathbb
{E}\tilde f_n(z)} +\frac{1}n
\operatorname{tr}R_n\bigl( I_N + \mathbb{E}\tilde
f_n(z) R_n\bigr)^{-1}.
\end{equation}
Subtracting \eqref{eqequation-t-tilde} to \eqref
{eqequation-f-tilde} finally yields
\begin{eqnarray*}
&& \mathbb{E}\tilde f_n(z) - \tilde t_n(z)
\\
&&\qquad= -A_n(z) \tilde t_n(z) \\
&&\qquad\quad{}\times\biggl[ 1 -
\frac{\tilde t_n(z) \mathbb{E}\tilde f_n(z)}n \operatorname{tr}R_n^2
\bigl(
I_N + \mathbb{E} \tilde f_n(z) R_n
\bigr)^{-1} \bigl( I_N + \tilde t_n(z)
R_n \bigr)^{-1} \biggr]^{-1},
\end{eqnarray*}
which is the counterpart of \cite{bai-silverstein-aop-2004}, equation
(4.12). The same arguments as in \cite
{bai-silverstein-aop-2004} now yield
%
\begin{eqnarray}
&& n \bigl(\mathbb{E}\tilde f_n(z) - \tilde
t_n(z) \bigr)
\nonumber
\\[-8pt]
\label{eqbias-computation}
\\[-8pt]
\nonumber
&&\qquad= -n A_n(z) \tilde t_n(z)
\biggl[ 1 - \frac{\tilde t_n^2(z)}n \operatorname{tr}R_n^2
\bigl( I_N + \tilde t_n(z) R_n
\bigr)^{-2} \biggr]^{-1} + o(1).
\end{eqnarray}
It remains to study the behavior of $nA_n(z)$. Following \cite{bai-silverstein-aop-2004}, equation (4.10), we obtain
\begin{eqnarray*}
&& nA_n(z)\\
&&\qquad =
\frac{b_n^2}n \mathbb{E}\operatorname{tr}Q_1 ( \mathbb{E}
\tilde f_n R_n + I_N )^{-1}R_n
Q_1 R_n - b_n^2 n \mathbb{E}
\biggl[ \biggl( \xi_1^* Q_1 \xi_1 -
\frac{1}n \operatorname{tr} Q_1 R_n \biggr)
\\
&&\qquad\quad{}\times\biggl( \xi_1^* Q_1 ( \mathbb{E}\tilde
f_n R_n + I_N )^{-1} \xi
_1 - \frac{1}n \operatorname{tr}Q_1 (
\mathbb{E}\tilde f_n R_n + I_N
)^{-1} R_n \biggr) \biggr] + o(1).
\end{eqnarray*}


Applying\vspace*{1pt} \eqref{eqcovariance-identity} to the right term to the RHS
of the previous equation (recall that $R^T=\overline{R}$), we obtain
%
\begin{eqnarray}
nA_n(z) &=& - |\mathcal{V}|^2
\frac{ b_n^2}{n} \mathbb{E}\operatorname{tr}R_n^{1/2}
Q_1 ( \mathbb{E}\tilde f_n R_n +
I_N )^{-1}R_n^{1/2}
\overline{R}_n^{1/2} Q^T_1 \overline
{R}^{1/2}_n
\nonumber
\\[-8pt]
\label{defdecomposition-bias}
\\[-8pt]
&&{}-\kappa\frac{b_n^2}n \sum_{i=1}^N
\bigl( R_n^{1/2} Q_1 R_n^{1/2}
\bigr)_{ii} \bigl( R_n^{1/2} Q_1 (
\mathbb{E}\tilde f_n R_n + I_N
)^{-1} R_n^{1/2} \bigr)_{ii} +o(1).\hspace*{-15pt}
\nonumber
\end{eqnarray}
The first term of the RHS has been fully analyzed in \cite
{bai-silverstein-aop-2004} in the case where $R_n$ and $X_n$ are real
matrices. We can
adapt these computations to the general case and get the following
identity\footnote{\label{footpreprintv3} Details can be found in the
previous version of this article, \arxivurl{arxiv:1309.3728v3}.}:
%
\begin{eqnarray}
&& - |\mathcal{V}|^2\frac{ b_n^2}{n}
\mathbb{E}
\operatorname{tr}R_n^{1/2} Q_1 ( \mathbb{E}
\tilde f_n R_n + I_N )^{-1}R_n^{1/2}
\overline{R}_n^{1/2} Q^T_1 \overline{R}^{1/2}_n
\nonumber
\\[-8pt]
\label{eqestimate-gaussian-term}
\\[-8pt]
\nonumber
&&\qquad= |\mathcal{V}|^2\frac{\frac{z^3 \tilde t_n^2}{n} \operatorname
{tr}R_n^{1/2} T^2_n(z)
R_n^{1/2} \overline{R}_n^{1/2} T^T_n(z) \overline{R}_n^{1/2}}{ 1 - \frac
{|\mathcal{V}
|^2z^2 \tilde t_n^2} n \operatorname{tr}R_n^{1/2} T_n(z) R_n^{1/2}
\overline{R}_n^{1/2}
T^T_n(z) \overline{R}_n^{1/2}}+ o(1),
\nonumber
\end{eqnarray}
where $T^T_n(z)$ is defined in \eqref{eqdeteq-RESOLVENT-conjugate}.
The term proportional to the cumulant in \eqref
{defdecomposition-bias} can be analyzed as in Section~\ref
{secdiagonal-terms-variance}, and one can prove that
%
\begin{eqnarray}
&&
-\kappa\frac{b_n^2}n \sum_{i=1}^N
\bigl( R_n^{1/2} Q_1 R_n^{1/2}
\bigr)_{ii} \bigl( R_n^{1/2} Q_1 (
\mathbb{E}\tilde f_n R_n + I_N
)^{-1} R_n^{1/2} \bigr)_{ii}
\nonumber
\\[-8pt]
\label{eqestimate-diag-term}
\\[-8pt]
\nonumber
&&\qquad= -\kappa\frac{z^2 \tilde t_n^2}n \sum_{i=1}^N
\bigl( R_n^{1/2} T_n R_n^{1/2}
\bigr)_{ii} \bigl( R_n^{1/2} T_n (
\tilde t_n R_n + I_N )^{-1}
R_n^{1/2} \bigr)_{ii} + o(1).
\end{eqnarray}
We now plug \eqref{eqestimate-gaussian-term} and \eqref
{eqestimate-diag-term} into \eqref{eqbias-computation} to conclude
\begin{eqnarray*}
&&\! n \bigl(\mathbb{E}\tilde f_n(z) - \tilde t_n(z)
\bigr)
\\
&&\!\qquad= - |\mathcal{V}|^2\frac{z^3 \tilde t_n^3}{n} \frac{ \operatorname
{tr}R_n^{1/2} T^2_n(z)
R_n^{1/2} \overline{R}_n^{1/2} T^T_n(z) \overline{R}_n^{1/2}}{ (1 - \frac
{|\mathcal{V}|^2z^2 \tilde t_n^2} n \operatorname{tr}R_n^{1/2} T_n(z)
R_n^{1/2} \overline{R}_n^{1/2} T^T_n(z) \overline{R}_n^{1/2} ) (1 - \frac
{\tilde
t_n^2} n \operatorname{tr}R_n^2 T_n^2 ) }
\\
&&\! \qquad\quad{}-\kappa\frac{z^3 \tilde t_n^3}n \sum
_{i=1}^N \frac{ (
R_n^{1/2} T_n R_n^{1/2} )_{ii} ( R_n^{1/2} T_n^2
R_n^{1/2} )_{ii}}{1 - \frac{z^2 \tilde t_n^2} n \operatorname{tr}R_n^2
T_n^2}+ o(1).
\end{eqnarray*}
The proof of Proposition~\ref{propconvergence-martingale} is completed.
\subsection{Proof of Proposition~\texorpdfstring{\protect\ref{propgaussian-process}}{5.2}}\label
{secgaussian-process}
Recall the meta-model introduced in Section~\ref{secmeta-model}.


\subsubsection{The Gaussian process $G_n$} \label{secmeta-matrix}

Let
%
\[
M_{n,M}(z)= \operatorname{tr}\bigl( {\Sigma}_n(M) {
\Sigma}_n(M)^*- zI_{NM}\bigr)^{-1} -
MNt_n(z).
\]
Applying Proposition~\ref{propconvergence-martingale} to the matrix
model $\Sigma_n(M) \Sigma_n(M)^*$ yields
\[
\forall z\in\Gamma, \qquad M_{n,M}^{ 1}(z) = \sum
_{j=1}^{nM} Z_j^M(z) +
o_P(1),
\]
where the $Z_j^M$'s are martingale increments and
\begin{eqnarray*}
 \sum_{j=1}^{nM} \mathbb{E}_{j-1}
Z_j^M(z_1) Z_j^M(z_2)
&\displaystyle \mathop{\hbox to 60pt{\rightarrowfill}}_{N,n \ \mathrm{fixed}, M\rightarrow\infty
}^{\mathcal P} & \Theta
_n(z_1,z_2),
\\
M_{n,M}^{ 2} (z) &\displaystyle \mathop{\hbox to 60pt{\rightarrowfill}}_{N,n \ \mathrm{fixed}, M
\rightarrow\infty}& {\mathcal B}_n(z).
\end{eqnarray*}

Notice that there is a genuine limit in the previous convergence.
Applying the central limit theorem for martingales \cite
{book-billinsley-proba-and-measure-95}, Theorem~35.12, plus the tightness
argument for $(M_{n,M}(z),z\in\Gamma)$ provided by Proposition~\ref
{propconvergence-martingale} immediately yields the fact that
$M_{n,M}$ converges in distribution to a Gaussian process $(G
_n(z),z\in\Gamma)$ with mean ${\mathcal B}_n(z)$ and covariance
function $\Theta_n(z_1,z_2)$.

\subsubsection{Tightness of the sequence of Gaussian processes
$(G_n)$}\label{sectightness-gaussian}

In order to prove that the sequence of Gaussian processes $(G_n)$
is tight, we shall prove, according to Prohorov's theorem, that it is
relatively compact in distribution. Consider the set of matrices
\[
\bigl\{ \bigl( R_n(M), M\ge1 \bigr); R_n \mbox{ is a } N
\times n \mbox{ matrix}, N=N(n); n\ge1 \bigr\},
\]
where $R_n(M)$ is defined in \eqref{eqdef-RM}. Since $\| R_n(M)\| =
\| R_n\|$ for all $M\ge1$, we have
\[
\sup_{M\ge1,N,n\rightarrow\infty} \bigl\| R_n(M)\bigr\| = \sup
_{N,n\rightarrow\infty}\| R_n\| <\infty
\]
by Assumption~\ref{assR}.
Hence, by Proposition~\ref
{propconvergence-martingale}, the family $\{ M_{n,M};M\ge1 \}
_{N,n\rightarrow\infty}$ is tight, hence relatively compact in
distribution. As the distribution ${\mathcal L}(G_n)$ of the
Gaussian process $G_n$ is the limit (in $M$) of the distribution
${\mathcal L}(M_{n,M})$ of $M_{n,M}$, ${\mathcal L}(G_n)$ belongs
to the closure of $\{ {\mathcal L}(M_{n,M})\}$, which is compact.
Finally, $\{ {\mathcal L}(G_n) \}$ is included in a compact set,
hence is relatively compact. In particular, the family of Gaussian
processes $(G_n)$ is tight.

\subsection{Proof of Theorem~\texorpdfstring{\protect\ref{lemmamain}}{1}}\label{secproof-lemma}


The two propositions below are minor variations of known results. They
will be helpful to conclude the proof of Theorem~\ref{lemmamain}.


\begin{lemma}[(CLT for martingales I)]\label{lemmavariations-martingales}
Suppose that for each $n$ $Y_{n1},
Y_{n2}, \ldots,\break Y_{n r_n}$ is a real martingale difference sequence
with respect to the increasing $\sigma$-field $\{ {\mathcal G}_{n,j}\}
$ having second moments. Assume moreover that $(\Theta_n^2)$ is a
sequence of nonnegative real numbers, uniformly bounded. If
\[
\sum_{j=1}^{r_n} \mathbb{E} \bigl(
Y_{nj}^2\vert{\mathcal G}_{n,j-1} \bigr) -
\Theta_n^2 \mathop{\xrightarrow}_{n\rightarrow\infty}^{\mathcal P}0,
\]
and for each $\varepsilon>0$,
\[
\sum_{j=1}^{r_n} \mathbb{E} \bigl(
Y_{nj}^2 1_{|Y_{nj}|>\varepsilon
} \bigr) \mathop{
\xrightarrow}_{n\rightarrow\infty} 0,
\]
then, for every bounded continuous function $f:\mathbb{R}\rightarrow
\mathbb{R}$,
%
\begin{equation}
\label{eqfonc-carac} \mathbb{E}f \Biggl( \sum_{j=1}^{r_n}
Y_{nj} \Biggr) - \mathbb{E}f(Z_n) \mathop{
\xrightarrow}_{n\rightarrow\infty} 0,
\end{equation}
where $Z_n$ is a centered Gaussian random variable with variance
$\Theta_n^2$.
\end{lemma}

Hereafter is the multidimensional and complex extension of Lemma~\ref
{lemmavariations-martingales} we shall rely on in the sequel.

\begin{lemma}[(CLT for martingales II)]\label{lemmamartingales-multidim}
Suppose that for each $n$ $(Y_{nj}; 1\le j\le r_n)$ is a $\mathbb{C}^d$-valued
martingale difference sequence with respect to the increasing $\sigma
$-field $\{ {\mathcal G}_{n,j}; 1\le j \le r_n\}$ having second
moments. Write
\[
Y_{nj}^T=\bigl(Y_{nj}^1,
\ldots,Y_{nj}^d\bigr).
\]

Assume moreover that $(\Theta_n(k,\ell))_n$ and $(
\widetilde\Theta_n(k,\ell))_n$ are uniformly bounded sequences of
complex numbers, for $1\le k,\ell\le d$. If
%
\begin{eqnarray}\label{eqconv-martingale-variance}
\sum_{j=1}^{r_n} \mathbb{E} \bigl(
Y^k_{nj} \overline{Y}_{nj}^d\vert {
\mathcal G}_{n,j-1} \bigr) - \Theta_n(k,\ell) & \displaystyle\mathop{
\xrightarrow}_{n\rightarrow\infty}^{\mathcal P}& 0,
\\
\label{eqconv-martingale-covariance}
\sum_{j=1}^{r_n} \mathbb{E} \bigl(
Y^k_{nj} Y_{nj}^\ell\vert {\mathcal
G}_{n,j-1} \bigr) - \widetilde\Theta_n(k,\ell) &\displaystyle \mathop{
\xrightarrow}_{n\rightarrow\infty}^{\mathcal P}& 0,
\end{eqnarray}
and for each $\varepsilon>0$,
%
\begin{equation}
\label{eqlyapounov}
\sum_{j=1}^{r_n} \mathbb{E}
\bigl( | Y_{nj}| ^2 1_{|Y_{nj}|>\varepsilon
} \bigr) \mathop{\xrightarrow}_{n\rightarrow\infty} 0,
\end{equation}
then, for every bounded continuous function $f:\mathbb{C}^d
\rightarrow\mathbb{R}$,
%
\begin{equation}
\label{eqfonc-caracc} \mathbb{E}f \Biggl( \sum_{j=1}^{r_n}
Y_{nj} \Biggr) - \mathbb{E}f(Z_n) \mathop{
\xrightarrow}_{n\rightarrow\infty} 0,
\end{equation}
where $Z_n$ is a $\mathbb{C}^d$-valued centered Gaussian random vector
with parameters
\[
\mathbb{E}Z_n Z_n^*= \bigl(\Theta_n(k,\ell)
\bigr)_{k,\ell}\quad \mbox{and}\quad \mathbb{E}Z_n
Z_n^T = \bigl(\widetilde\Theta_n(k,\ell)
\bigr)_{k,\ell}.
\]
\end{lemma}

Lemmas \ref{lemmavariations-martingales} and \ref
{lemmamartingales-multidim} are variations around the central limit
theorem for martingales (see Billingsley \cite
{book-billinsley-proba-and-measure-95}, Theorem~35.12) which enables us
to prove
(in the real case)
\[
\forall t\in\mathbb{R},\qquad \mathbb{E}e^{\mathbf{i}t\sum
_{j=1}^{r_n} Y_{nj}} - e^{-({t^2 \Theta_n^2})/2}
\rightarrow0
\]
and L\'{e}vy theorem for the weak convergence criterion via
characteristic functions (see Kallenberg \cite
{book-kallenberg-second-edition}, Theorems~5.3 and~5.5) which yields
\eqref
{eqfonc-caracc} from the above convergence. Details of the proof are omitted.

\begin{lemma}[(Tightness and weak convergence)]\label{lemmatightness}
Let $K$ be a compact set in $\mathbb{C}$; let $X_1,X_2,\ldots$ and
$Y_1,Y_2,\ldots$ be random elements in $C(K,\mathbb{C})$. Assume that
for all
$d\ge1$, for all $z_1,\ldots,z_d\in K$, for all $f\in C(\mathbb
{C}^d,\mathbb{C})$
we have
\[
\mathbb{E}f\bigl(X_n(z_1),\ldots,X_n(z_d)
\bigr) - \mathbb{E}f\bigl(Y_n(z_1),\ldots
,Y_n(z_d)\bigr) \mathop{\xrightarrow}_{n\rightarrow\infty} 0.
\]
Assume moreover that $(X_n)$ and $(Y_n)$ are tight, then for every
continuous and bounded functional $F:C(K,\mathbb{C})\rightarrow
\mathbb{C}$, we have
\[
\mathbb{E}F(X_n) - \mathbb{E}F(Y_n) \mathop{
\xrightarrow}_{n\rightarrow\infty} 0.
\]
\end{lemma}

Lemma~\ref{lemmatightness} can be proved as \cite
{book-kallenberg-second-edition}, Lemma~16.2; the proof is therefore omitted.

We are now in position to conclude.

In order to apply Lemma~\ref{lemmamartingales-multidim}, it remains
to check that $\Theta_n$ as defined in \eqref
{eqdecomposition-covariance} is uniformly bounded for $z_1,z_2\in
\Gamma$ fixed but this is an easy byproduct of Proposition~\ref
{propgaussian-process}.

Proposition~\ref{propconvergence-martingale} together with Lemma~\ref
{lemmamartingales-multidim} (notice that condition \eqref
{eqlyapounov} can be proved as in \cite{bai-silverstein-aop-2004})
yield the fact that for every $z_1,\ldots, z_d\in\Gamma$ and for
every bounded continuous function $f:\Gamma^d\rightarrow\mathbb{C}$
\[
\mathbb{E}f\bigl(M_n(z_1),\ldots,M_n(z_d)
\bigr) - \mathbb{E}f\bigl(G_n(z_1),\ldots,G
_n(z_d)\bigr) \mathop{\xrightarrow}_{N,n\rightarrow\infty} 0,
\]
where $G_n$ is well defined by Proposition~\ref
{propgaussian-process}. Now the tightness of $M_n$ and $G_n$
together with Lemma~\ref{lemmatightness}
yield the last statement of Theorem~\ref{lemmamain}.

%

\section{Proof of Theorem~\texorpdfstring{\protect\ref{thnon-analytic}}{2} (fluctuations for
nonanalytic functionals)}\label{secproof-non-analytic-I}

In this section, we will assume that the random variables $(X_{ij}^n)$
are truncated, centered and normalized, following Section~\ref{sectruncation}.

\subsection{Useful properties}\label{secHS-properties}

Recall that ${\mathcal S}_n \subset{\mathcal S}_\infty\stackrel
{\triangle}{=}
[0,\bolds{\lambda}^+_R ( 1 + \sqrt{\bolds{\ell
}^+} )^2 ]$
uniformly in $n$. Denote by $\mathbf{h}\in C_c^{\infty}(\mathbb
{R})$ a
function whose value is 1 on a $\eta$-neighborhood ${\mathcal
S}_\infty^\eta$ of ${\mathcal S}_\infty$.

\begin{prop}\label{propreduction}
1. Assume that Assumptions \ref{assX}  and \ref{assR}  hold true; let
the random variables\vspace*{1pt} $(X_{ij}^n)$ be truncated as in Section~\ref
{sectruncation}, function $\mathbf{h}$ be defined as above and
$f:\mathbb{R}\to\mathbb{R}$ be a continuous function. Then
\begin{eqnarray*}
&& \operatorname{tr}f\bigl(\Sigma_n \Sigma_n^*\bigr) -
\operatorname{tr}(f\mathbf{h}) \bigl(\Sigma_n \Sigma
_n^*\bigr) \mathop{\xrightarrow}_{N,n\to\infty}^{\mathrm{a.s.}} 0.
\end{eqnarray*}
%
%
2. Let $\mathbf{h}_n$ be a smooth function on $\mathbb{R}$
with compact support and whose value is 1 on a $\eta$-neighborhood
${\mathcal S}_n^\eta$ of ${\mathcal S}_n$; then
\begin{eqnarray*}
\int_\mathbb{R}f(\lambda) {\mathcal F}_n(d\lambda)
&=& \int_\mathbb{R} (f\mathbf{h}_n) (\lambda) {
\mathcal F}_n(d\lambda).
\end{eqnarray*}
%
%
\end{prop}

The proof of Proposition~\ref{propreduction} is straightforward and is
based on the fact that almost surely
\[
\limsup_{N,n\to\infty} \bigl\| \Sigma_n \Sigma_n^*
\bigr\| < \bolds{\lambda}^+_R \bigl( 1 + \sqrt{\bolds{\ell}^+}
\bigr)^2 +\eta,
\]
a fact that can be found in \cite{book-bai-silverstein} for instance.
Details are left to the reader.

%


The following proposition underlines how a sufficient regularity of
function $f$ compensates the singularity in $\operatorname{Im}(z)^{-1}$
near the
real axis.

\begin{prop}\label{propcompensation}
Let $\mu, \nu$ be two probability measures on $\mathbb{R}$ and
$g_\mu$ and
$g_\nu$ their associated Stieltjes transforms. Assume that
\[
\bigl| g_\mu(z)- g_\nu(z)\bigr| \le\frac{|h(z)|}{\operatorname{Im}(z)^k}, \qquad z\in
\mathbb{C} ^+,
\]
where $h$ is a continuous function over $\mathrm{cl}(\mathbb{C}^+)$, the
closure of $\mathbb{C}^+$.

Let $f:\mathbb{R}\to\mathbb{R}$ be a function of order $C^{k+1}$
with bounded
support; recall the definition of $\Phi_k(f)$ in \eqref{eqalmost-analytic}
and
denote by
\[
\| f\|_{k+1}= \sup_{0\le\ell\le k+1} \bigl\| f^{(\ell)}
\bigr\|_{\infty} \qquad\mbox{where } \|g\|_\infty=\sup_{x\in\mathbb{R}}
\bigl|g(x)\bigr|.
\]
Then
%
\begin{eqnarray}
\biggl\llvert\int f \,d\mu- \int f \,d\nu\biggr\rrvert &\le& \frac
{1}\pi
\biggl\llvert\int_{\mathbb{C}^+} \overline{\partial}
\Phi_k(f) (z) \bigl\{ g_\mu(z) - g_\nu(z)
\bigr\} \ell_2( dz)\biggr\rrvert
\nonumber
\\
\label{estimate-HS}
&\le& K \| f\|_{k+1} \int_{\mathrm{supp}(f)\times\operatorname{supp}(\chi)} \bigl|h(z)\bigr|
\ell_2(dz)
\\
&\le& K' \| f\|_{k+1}.
\nonumber
\end{eqnarray}
%
\end{prop}

\begin{pf}
Write
\begin{eqnarray*}
\overline{\partial} \Phi_k(f) (x+ \mathbf{i}y) &=&
\partial_x \Phi_k(f) (x+ \mathbf{i}y) +\mathbf{i}
\partial_y \Phi_k(f) (x+ \mathbf{i}y)
\\
&=& \frac{(\mathbf{i}y)^k f^{(k+1)}(x)}{k!} \chi(y) + \mathbf{i}\sum
_{\ell=0}^k
\frac{(\mathbf{i}y)^\ell f^{(\ell)}(x)}{\ell!} \chi'(y).
\end{eqnarray*}
From this and the fact that $\chi$ is equal to 1 for $y$ small enough,
we deduce that
\[
\overline{\partial} \Phi_k(f) (x+ \mathbf{i}y) = \frac{(\mathbf
{i}y)^k f^{(k+1)}(x)}{k!}
\]
near the real axis. Hence, $\llvert \overline{\partial} \Phi_k(f)(x+
\mathbf{i}y) \rrvert \le1_{\mathrm{supp}(f)\times\mathrm{supp}(\chi)
}(x,y) K \|f\|_{k+1} y^k$ near the real axis, which yields \eqref{estimate-HS}.
\end{pf}

\subsection{Proof of Theorem~\texorpdfstring{\protect\ref{thnon-analytic}}{2}}\label
{proofnon-analytic-I}
Recall the definition of the sets $D$, $D^+$ and $D_\varepsilon$ given
in \eqref{defsets-D} and the fact that constant $A> \bolds{\lambda}^+_R
( 1 + \sqrt{\bolds{\ell}^+} )^2$.

\begin{lemma}\label{lemmaabstract}
Let $(\varphi_n(z), z\in D^+\cup\overline{D^+})_{n\in\mathbb{N}}$ and
$(\psi_n(z),z \in D^+\cup\overline{D^+})_{n\in\mathbb{N}}$ be centered
complex-valued continuous random processes and such that $\varphi(\bar
{z})=\overline{\varphi(z)}$ and $\psi(\bar{z})=\overline{\psi(z)}$.
Assume that:
\begin{longlist}[(iii)]
\item[(i)] The following convergence in distribution holds true: for
all $d\ge1$ and $(z_1,\ldots,z_d) \in D^+$,
\[
d_{\mathcal{LP}} \bigl( \bigl( \varphi_n(z_1),\ldots,
\varphi_n(z_d)\bigr), \bigl(\psi_n(z_1),
\ldots,\psi_n(z_d)\bigr) \bigr)\mathop{
\xrightarrow}_{n\rightarrow\infty} 0.
\]

\item[(ii)] For all $\varepsilon>0$, $\varphi_n(z)$ and $\psi_n(z)$
are tight on $D_\varepsilon$.
\item[(iii)] The process $(\psi_n(z))$ is Gaussian with covariance
matrix $\kappa_n(z_1,z_2)$, $(z_1,z_2\in D^+\cup\overline{D^+})$.
\item[(iv)] The following estimates hold true:
\[
\forall n\in\mathbb{N},\forall z\in D^+, \qquad\operatorname{var} \varphi
_n(z) \le\frac{1}{\operatorname{Im}(z)^{2k}} \quad\mbox{and}\quad \operatorname{var}
\psi_n(z) \le\frac{1}{\operatorname{Im}(z)^{2k}}.
\]
\item[(v)] Let functions $g_\ell:\mathbb{R}\to\mathbb{R}$ $(1\le
\ell\le L$) be
$C^{k+1}$ and have compact support.
\end{longlist}
Then
\begin{eqnarray*}
&& d_{\mathcal{LP}} \biggl( \frac{1}\pi\operatorname{Re} \int
_{\mathbb{C}^+} \overline\partial\Phi_k({\mathbf g}) (z)
\varphi_n(z) \ell_2(dz), \frac{1}\pi\operatorname{Re}
\int_{\mathbb{C}^+} \overline\partial\Phi_k ({\mathbf g})
(z) \psi_n(z) \ell_2(dz) \biggr) \\
&&\qquad\mathop{
\xrightarrow}_{n\rightarrow\infty} 0,
\end{eqnarray*}
where
\begin{eqnarray*}
\overline\partial\Phi_k(g_j) (z)  &=& (
\partial_x+\mathbf{i}\partial_y) \sum
_{\ell=0}^k \frac{(\mathbf{i}y)^\ell}{\ell!} g_j^{(\ell)}(x)
\chi(y) \quad\mbox{and}\\
\overline\partial\Phi_k({\mathbf g})&=& \bigl(
\overline\partial\Phi_k(g_j) ; 1\le j\le L \bigr)
\end{eqnarray*}
with $\chi$ being smooth, compactly supported with value 1 in a
neighborhood of 0. Moreover,
\[
\frac{1}\pi\operatorname{Re} \int_{\mathbb{C}^+} \overline\partial
\Phi_k({\mathbf g}) (z) \psi_n(z) \ell_2(dz)
\]
is centered Gaussian with covariance matrix
%
\begin{eqnarray}
&& \operatorname{cov} \biggl( \frac{1}\pi
\operatorname{Re} \int
_{\mathbb{C}^+} \overline\partial\Phi_k(g_k)
(z) \psi_n(z) \ell_2(dz), \frac{1}\pi\operatorname{Re}
\int_{\mathbb{C}^+} \overline\partial\Phi_k(g_\ell)
(z) \psi_n(z) \ell_2(dz) \biggr)
\nonumber\\
\label{eqcov-non-analytic}
&&\qquad= \frac{1}{2\pi^2} \operatorname{Re} \int_{(\mathbb{C}^+)^2} \overline{
\partial} \Phi_k(g_k) (z_1) \overline{
\overline{\partial} \Phi_k(g_\ell) (z_2)}
\kappa_n(z_1,\bar{z}_2) \ell_2(dz_1)
\ell_2 (dz_2)
\\
&&\qquad\quad{}+ \frac{1}{2\pi^2} \operatorname{Re} \int_{(\mathbb{C}^+)^2} \overline{
\partial} \Phi_k(g_k) (z_2) \overline{
\partial} \Phi_k(g_\ell) (z_2)
\kappa_n(z_1,z_2) \ell_2(dz_1)
\ell_2(dz_2),
\nonumber
\end{eqnarray}
for $1\le k,\ell\le L$.
\end{lemma}

The proof of Lemma~\ref{lemmaabstract} is provided in Appendix~\ref{prooflemmaabstract}.

The strategy to prove Theorem~\ref{thnon-analytic} closely follows
this lemma. Denote by
\[
\varphi_n(z) =\operatorname{tr}Q_n(z) - \mathbb{E}
\operatorname{tr}Q_n(z)\quad \mbox{and}\quad \psi_n(z)=
G_n(z) - \mathbb{E}G_n(z),
\]
the process $G_n$ being defined in Theorem~\ref{lemmamain},
then conditions (i), (ii) and (iii) are immediate consequences of
Theorem~\ref{lemmamain}.
In order to check condition (iv), we establish the following proposition.
%
\begin{prop}\label{propvariance-estimate-sharp}
Assume that
Assumptions \ref{assX}  and \ref{assR}  hold true, then:
\begin{longlist}[(ii)]
\item[(i)] (Bordenave \cite{bordenave-personal-comm-2013}, Hachem et
al. \cite{hachem-et-al-2008}, Lemma~6.3,
Shcherbina \cite{shcherbina-2011}). For all $z\in\mathbb{C}^+$,
\[
\operatorname{var} \operatorname{tr}Q_n(z)\le\frac{C}{\operatorname
{Im}(z)^4}.
\]
\item[(ii)] For all $z\in\mathbb{C}^+$,
\[
\operatorname{var} G_n(z)\le\frac{C}{\operatorname{Im}(z)^4},
\]
\end{longlist}
where $C$ is a constant that may depend polynomially on $|z|$.
\end{prop}

The first part of the proposition is classical and its proof is omitted
(for the details, see footnote \ref{footpreprintv3}).
Proof of Proposition~\ref{propvariance-estimate-sharp}(ii) is
postponed to Appendix~\ref{appprop-variance-sharp}.

Taking into account the estimates established in Proposition~\ref
{propvariance-estimate-sharp} immediately yields the first part of
Theorem~\ref{thnon-analytic} in the case where functions $(g_\ell)$
have a bounded support and satisfy (v) with $k=2$, that is, are $C^3$.
It remains to prove the equivalence between \eqref{eqfluctuations-1}
and \eqref{eqfluctuations-2}, but this immediately follows from the following.
%
\begin{prop}\label{propequivlevy} Let $(X_n)$ and $(Y_n)$ be
$\mathbb{C}
^d$-valued random variables and assume that both sequences are tight,
then the following are equivalent:
\begin{longlist}[(ii)]
\item[(i)] the following convergence holds true:
$
d_{\mathcal{LP}}(X_n,Y_n) \displaystyle\mathop{\xrightarrow}_{n\to\infty}
0$;
\item[(ii)] for every continuous bounded function $f:\mathbb{C}^d\to
\mathbb{C}$,
$ \mathbb{E} f(X_n) -\break \mathbb{E} f(Y_n) \displaystyle\mathop{\xrightarrow}_{n\to
\infty} 0$.
\end{longlist}
\end{prop}

Proposition~\ref{propequivlevy} can be proved easily by
contradiction using the fact that $d_{\mathcal{LP}}$ meterizes the
convergence of
laws; its proof is hence omitted.

\subsection{Proof of Proposition~\texorpdfstring{\protect\ref{propidentification-distribution}}{4.1}}\label{proofidentification-distribution}
Let $f\in C_c^\infty(\mathbb{R}^2)$. A simple but lengthy computation yields
the fact that
%
\begin{eqnarray}
&& \overline{\partial}_2\overline{\partial}_1
\Phi_{N_1,N_2}(f) (x+\mathbf{i} u,y+\mathbf{i}v)
\nonumber
\\[-8pt]
\label{eqderivative-2D}
\\[-8pt]
\nonumber
&&\qquad= \frac{\partial
^{N_1+N_2+2}}{\partial
x^{N_1+1}\partial
y^{N_2+1}} f
(x,y)\times\frac{(\mathbf{i}u)^{N_1}}{N_1!} \frac{(\mathbf{i}v)^{N_2}}{N_2!}
\end{eqnarray}
for $u,v$ small enough. Let now $N_1=N_2=2$. Since $|\Theta
_n(z_1,z_2)|\le K|z_1 z_2|^{-2}$ for any $z_1,z_2\in\mathbb{C}^+$ and
$z_1,z_2$ in a compact set (use Cauchy--Schwarz and apply
Proposition~\ref{propvariance-estimate-sharp}), $\Upsilon(f)$ is well defined.
Let $K$ be a compact set in $\mathbb{R}^2$ and let $f\in C_c^{\infty
}(\mathbb{R}^2)$
with support included in $K$, then one can easily prove that
\[
\bigl|\Upsilon(f)\bigr| \le C_K \| f\|_{3,3} \qquad\mbox{with } \|f\|
_{3,3} = \mathop{\sup_{\ell,p\le3}}_{(x,y)\in K}
\bigl\| \partial^\ell_x
\partial^p_y f(x,y)\bigr\|_{\infty}.
\]
This in\vspace*{1pt} particular implies that $\Upsilon$ is a distribution on
$C_c^\infty(\mathbb{R}^2)$, of finite order $(3,3)$, and hence uniquely
extends as a distribution on $C^{3,3}_c(\mathbb{R}^2)$.

Moreover, $\Upsilon(f)$ can be written as
\begin{eqnarray*}
\Upsilon(f) &=& \lim_{\varepsilon\downarrow0} \frac{1}{2\pi^2}
\operatorname{Re} \int
_{(\mathbb{C}_\varepsilon^+)^2} \overline{\partial}_2\overline{
\partial}_1 \Phi_{2,2} (f) (z_1,z_2)
\Theta_n(z_1,z_2) \ell_2
(dz_1)\ell_2(dz_2)
\\
&&{}+\lim_{\varepsilon\downarrow0} \frac{1}{2\pi^2} \operatorname{Re} \int
_{(\mathbb{C}_\varepsilon^+)^2} \overline{\partial}_2\overline{\partial
}_1 \Phi_{2,2} (f) (z_1,\overline{z_2})
\Theta_n(z_1,\overline{z_2})
\ell_2(dz_1)\ell_2(dz_2),
\end{eqnarray*}
where $\mathbb{C}^+_\varepsilon=\{ z\in\mathbb{C},\operatorname
{Im}(z)\ge\varepsilon\}$.
Taking into account the facts that
\[
\overline{\overline{\partial}_2\overline{\partial}_1
\Phi_{n_1,n_2}(f) (z_1,z_2)} = \overline{
\partial}_2\overline{\partial}_1 \Phi_{n_1,n_2}(f)
(\overline{z_1},\overline{z_2})\quad \mbox{and} \quad\overline{
\Theta_n(z_1,z_2)}= \Theta_n(
\overline{z_1},\overline{z_2}),
\]
we can expand $\Upsilon(f)$ as
\begin{eqnarray*}
\Upsilon(f) &=&\lim_{\varepsilon\downarrow0}\frac{1}{4\pi^2} \int
_{(\mathbb{C}_\varepsilon^+)^2} \overline{\partial}_2\overline{\partial
}_1 \Phi_{2,2} (f) (z_1,z_2)
\Theta_n(z_1,z_2) \ell_2(dz_1)
\ell_2(dz_2)
\\
&& {}+\lim_{\varepsilon\downarrow0}\frac{1}{4\pi^2} \int_{(\mathbb{C}
_\varepsilon^+)^2}
\overline{\partial}_2\overline{\partial}_1 \Phi
_{2,2} (f) (\overline{z_1},\overline{z_2})
\Theta_n(\overline{z_1},\overline{z_2})
\ell_2(dz_1)\ell_2(dz_2)
\\
&&{}+\lim_{\varepsilon\downarrow0}\frac{1}{4\pi^2} \int_{(\mathbb{C}
_\varepsilon^+)^2}
\overline{\partial}_2\overline{\partial}_1 \Phi
_{2,2} (f) (z_1,\overline{z_2})
\Theta_n(z_1,\overline{z_2})
\ell_2 (dz_1)\ell_2(dz_2)
\\
&& {}+\lim_{\varepsilon\downarrow0}\frac{1}{4\pi^2} \int_{(\mathbb{C}
_\varepsilon^+)^2}
\overline{\partial}_2\overline{\partial}_1 \Phi
_{2,2} (f) (\overline{z_1},z_2)
\Theta_n(\overline{z_1},z_2)
\ell_2 (dz_1)\ell_2(dz_2).
\end{eqnarray*}
We now apply twice Green's formula to each integral and obtain
\begin{eqnarray*}
\Upsilon(f) &=&-\lim_{\varepsilon\downarrow0}\frac{1}{4\pi^2} \int
_{\mathbb{R}^2} \Phi_{2,2} (f) (x_1+\mathbf{i}
\varepsilon,x_2+\mathbf{i}\varepsilon)\Theta_n(x_1+
\mathbf{i}\varepsilon,x_2+\mathbf{i}\varepsilon) \,dx_1
\,dx_2
\\
&& {}- \lim_{\varepsilon\downarrow0}\frac{1}{4\pi^2} \int_{\mathbb{R}
^2}
\Phi_{2,2} (f) (x_1-\mathbf{i}\varepsilon,x_2-
\mathbf{i}\varepsilon)\Theta_n (x_1-\mathbf{i}
\varepsilon,x_2-\mathbf{i}\varepsilon) \,dx_1
\,dx_2
\\
&& {}+ \lim_{\varepsilon\downarrow0}\frac{1}{4\pi^2} \int_{\mathbb{R}
^2}
\Phi_{2,2} (f) (x_1+\mathbf{i}\varepsilon,x_2-
\mathbf{i}\varepsilon)\Theta_n (x_1+\mathbf{i}
\varepsilon,x_2-\mathbf{i}\varepsilon) \,dx_1
\,dx_2
\\
&& {}+\lim_{\varepsilon\downarrow0}\frac{1}{4\pi^2} \int_{\mathbb{R}^2}
\Phi_{2,2} (f) (x_1-\mathbf{i}\varepsilon,x_2+
\mathbf{i}\varepsilon)\Theta_n (x_1-\mathbf{i}
\varepsilon,x_2+\mathbf{i}\varepsilon) \,dx_1
\,dx_2.
\end{eqnarray*}
Notice that the sign changes in the two last integrals follow from the
contour orientations in Green's formula. We now prove
%
\begin{eqnarray}
&&\lim_{\varepsilon\downarrow0} \int_{\mathbb{R}^2}
\Phi_{2,2} (f) (x_1+\mathbf{i} \varepsilon,x_2+
\mathbf{i}\varepsilon)\Theta_n(x_1+\mathbf{i}\varepsilon
,x_2+\mathbf{i}\varepsilon) \,dx_1 \,dx_2
\nonumber
\\[-8pt]
\label{boundary-value-eop}
\\[-8pt]
\nonumber
&&\qquad =\lim_{\varepsilon\downarrow0} \int_{\mathbb{R}^2} f
(x_1,x_2)\Theta_n(x_1+
\mathbf{i}\varepsilon,x_2+\mathbf{i}\varepsilon) \,dx_1
\,dx_2.
\end{eqnarray}
The three other integrals can be handled similarly, and this will
achieve the boundary value representation \eqref{boundary-value-2} for
$\Upsilon(f)$.

Using the mere definition of $\Phi_{N_1,N_2}(f)$ [cf. \eqref
{defAAE-bidim}] and Green's formula, we get
\begin{eqnarray*}
&&\int_{(\mathbb{C}_\varepsilon^+)^2} \overline{\partial}_2
\overline{\partial}_1 \Phi_{1,0} (f) (z_1,z_2)
\Theta_n(z_1,z_2) \ell_2(dz_1)
\ell_2(dz_2)
\\
&&\qquad= -\int_{\mathbb{R}^2} \Phi_{1,0} (f) (x_1+
\mathbf{i}\varepsilon,x_2+\mathbf{i} \varepsilon)
\Theta_n(x_1+\mathbf{i}\varepsilon,x_2+
\mathbf{i}\varepsilon) \,dx_1 \,dx_2
\\
&&\qquad= -\int_{\mathbb{R}^2} f (x_1,x_2)
\Theta_n(x_1+\mathbf{i}\varepsilon,x_2+
\mathbf{i} \varepsilon) \,dx_1 \,dx_2
\\
&&\qquad\quad {}- \mathbf{i}\varepsilon\int_{\mathbb{R}^2} \partial_x
f (x_1,x_2)\Theta_n(x_1+\mathbf
{i}\varepsilon,x_2+\mathbf{i}\varepsilon) \,dx_1
\,dx_2.
\end{eqnarray*}
We extract\vspace*{1pt} the first term of the RHS from the equation above.
Taking into account~\eqref{eqderivative-2D} and the fact that
$|\Theta_n(z_1,z_2)|\le|z_1 z_2|^{-2}$ for $z_1,z_2$ in a compact set
of $\mathbb{C}\setminus\mathbb{R}$, we obtain
\[
\limsup_{\varepsilon\downarrow0} \biggl\llvert\varepsilon^3 \int
_{\mathbb{R}
^2} f (x_1,x_2)
\Theta_n(x_1+\mathbf{i}\varepsilon,x_2+
\mathbf{i}\varepsilon) \,dx_1 \,dx_2\biggr\rrvert<\infty.
\]
By applying the same argument to the quantity
\[
\int_{(\mathbb{C}_\varepsilon^+)^2} \overline{\partial}_2\overline{
\partial}_1 \Phi_{4-\ell,0} (f) (z_1,z_2)
\Theta_n(z_1,z_2) \ell_2
(dz_1)\ell_2(dz_2)
\]
for $\ell=2$ then $\ell=1$ and $\ell= 0$, we can similarly prove that
%
\begin{eqnarray}
&& \limsup_{\varepsilon\downarrow0} \biggl\llvert\varepsilon^\ell\int
_{\mathbb{R}^2} f (x_1,x_2)
\Theta_n(x_1+\mathbf{i}\varepsilon,x_2+
\mathbf{i}\varepsilon) \,dx_1 \,dx_2\biggr\rrvert<\infty
\nonumber\\
\eqntext{\displaystyle\mbox{for } \ell= 2, 1, 0.}
\end{eqnarray}
We finally obtain
\[
\limsup_{\varepsilon\downarrow0} \biggl\llvert\int_{\mathbb{R}^2} f
(x_1,x_2)\Theta_n(x_1+
\mathbf{i}\varepsilon,x_2+\mathbf{i}\varepsilon) \,dx_1
\,dx_2\biggr\rrvert<\infty.
\]
Expanding $\Phi_{2,2}(f)$ into \eqref{boundary-value-eop} and using
the above estimate immediately\break yields~\eqref{boundary-value-eop}.

The proof of Proposition~\ref{propidentification-distribution} is complete.

\subsection{Proof of Proposition~\texorpdfstring{\protect\ref{propnew-variance-formula}}{4.2}}\label{proofnew-variance-formula} The
covariance writes (in short)
\begin{eqnarray*}
&& \operatorname{cov}\bigl(Z_n^1(f),Z_n^1(g)
\bigr) \\
&&\qquad= -\frac{1}{4\pi^2} \lim_{\varepsilon
\downarrow0} \sum
_{\pm_1,\pm_2} (\pm_1\pm_2)\int f(x) g(y)
\Theta_{n}(x\pm_1\mathbf{i}\varepsilon,y\pm_2
\mathbf{i}\varepsilon) \,dx \,dy,
\end{eqnarray*}
where\vspace*{1pt} $\pm_1,\pm_2\in\{ +,-\}$ and $\pm_1\pm_2$ is the sign
resulting from the product $\pm_1 1$ by $\pm_2 1$. Unfolding $\Theta
_n=\Theta_{0,n}+|{\mathcal V}|^2 \Theta_{1,n} + \kappa\Theta
_{2,n}$, we have three terms to compute. According to the assumptions of
Proposition~\ref{propnew-variance-formula}, either $|{\mathcal V}|^2$
equals 1 or 0. In the latter case, the term corresponding to $\Theta
_{1,n}$ vanishes; if $|\mathcal V|^2=1$, then the quantities ${\mathcal
A}_n$ and ${\mathcal A}_{0,n}$ [resp., defined in \eqref{defAn} and
\eqref{defAn0}] are equal, and so are $\Theta_{0,n}$ and $\Theta
_{1,n}$. We first establish
\begin{eqnarray}
&&-\frac{1}{4\pi^2} \lim_{\varepsilon\downarrow0} \sum
_{\pm_1,\pm
_2} (\pm_1\pm_2) \int f(x) g(y)
\Theta_{0,n}(x\pm_1\mathbf{i}\varepsilon,y
\pm_2 \mathbf{i}\varepsilon) \,dx \,dy
\nonumber
\\[-8pt]
\label{eqcov0}
\\[-8pt]
\nonumber
&&\qquad= \frac{1}{2\pi^2} \int_{{\mathcal S}_n^2} f'(x)g'(y)
\ln\biggl\llvert\frac{{\tilde t}_n(x) -\overline{{\tilde t}_n(y)}}{
{\tilde t}_n(x) - {\tilde t}_n(y)}\biggr\rrvert\,dx \,dy.
\nonumber
\end{eqnarray}
The proof relies on formula \eqref{deftheta0-alternate} and the
following expression of ${\mathcal A}_{0,n}$:
%
\begin{equation}
\label{defA0-alternate} 1-{\mathcal A}_{0,n}(z_1,z_2) =
\frac{(z_1-z_2)\tilde t_n(z_1) \tilde
t_n(z_2)}{\tilde t_n(z_1) - \tilde t_n(z_2)}
\end{equation}
which can be obtained using \eqref{eqdef-t-tilde}. Using \eqref
{deftheta0-alternate} and performing a double integration by parts yields
\begin{eqnarray*}
&& \int f(x) g(y) \Theta_{0,n}(x+\mathbf{i}\varepsilon, y +
\mathbf{i} \varepsilon) \,dx \,dy
\\
&&\qquad= \int f'(x)g'(y) \ln\bigl\llvert1- {\mathcal
A}_{0,n} (x+\mathbf{i} \varepsilon, y +\mathbf{i}\varepsilon)\bigr
\rrvert\,dx \,dy
\\
&&\qquad\quad {}+ \mathbf{i}\int f'(x)g'(y) \operatorname{Arg} \bigl(
1- {\mathcal A}_{0,n} (x+\mathbf{i}\varepsilon, y +\mathbf{i}
\varepsilon) \bigr) \,dx \,dy.
\end{eqnarray*}
Following \cite{bai-silverstein-aop-2004}, Section~5, we need only to
consider the logarithm term and show its convergence since the $\textrm
{Arg}$ term will eventually disappear (functions $f$ and $g$ being
real, the covariance will be real as well). Using \eqref
{defA0-alternate}, we obtain
\begin{eqnarray*}
&& \int f'(x)g'(y) \ln\bigl\llvert1- {\mathcal
A}_{0,n} (x+\mathbf{i}\varepsilon, y +\mathbf{i}\varepsilon)\bigr
\rrvert\,dx \,dy
\\
&&\qquad= \int f'(x)g'(y) \ln\biggl\llvert
\frac{(x-y) \tilde t_n(x+\mathbf
{i}\varepsilon)
\tilde t_n(y+\mathbf{i}\varepsilon) }{\tilde t_n(x+\mathbf
{i}\varepsilon) -
\tilde t_n(y+\mathbf{i}\varepsilon)}\biggr\rrvert\,dx \,dy
\end{eqnarray*}
and the sum writes
\begin{eqnarray*}
&& \sum_{\pm_1,\pm_2} (\pm_1
\pm_2)\int f(x) g(y) \Theta_{n}(x\pm_1
\mathbf{i}\varepsilon,y\pm_2 \mathbf{i}\varepsilon) \,dx \,dy
\\
&&\qquad = 2\int f'(x)g'(y) \ln\biggl\{ \biggl\llvert
\frac{(x-y) \tilde
t_n(x+\mathbf{i}
\varepsilon) \tilde t_n(y+\mathbf{i}\varepsilon) }{\tilde
t_n(x+\mathbf{i}
\varepsilon) - \tilde t_n(y+\mathbf{i}\varepsilon)}\biggr\rrvert
\\
&& 
\qquad\quad{}\times\biggl\llvert\frac{\tilde t_n(x+\mathbf{i}\varepsilon) - \tilde
t_n(y-\mathbf{i}
\varepsilon)}{(x-y+ 2\mathbf{i}\varepsilon) \tilde t_n(x+\mathbf
{i}\varepsilon)
\tilde t_n(y-\mathbf{i}\varepsilon) }
\biggr\rrvert\biggr\} \,dx \,dy
\\
&&\qquad\stackrel{(a)}= 2\int f'(x)g'(y) \biggl\{ \ln
\biggl\llvert\frac{x-y
}{x-y+ 2\mathbf{i}\varepsilon} \biggr\rrvert+ \ln\biggl\llvert
\frac{\tilde
t_n(x+\mathbf{i}\varepsilon) -\tilde t_n(y-\mathbf{i}\varepsilon)
}{\tilde
t_n(x+\mathbf{i}\varepsilon) - \tilde t_n(y+\mathbf{i}\varepsilon
)}\biggr\rrvert\biggr\} \,dx \,dy,
\end{eqnarray*}
where $(a)$ follows from the fact that $\tilde t_n(\bar z) =\overline
{\tilde t_n(z)}$ and $|z| =|\bar{z}|$. It is straightforward to prove
that the first integral of the RHS vanishes as $\varepsilon\to0$.
Using similar arguments as in \cite{bai-silverstein-aop-2004},
Section~5, one can prove that
%
\begin{eqnarray*}
&& \sum_{\pm_1,\pm_2} (\pm_1\pm_2)
\int f(x) g(y) \Theta_{n}(x\pm_1\mathbf{i}\varepsilon,y
\pm_2 \mathbf{i}\varepsilon) \,dx \,dy
\\
&&\qquad = 2\int f'(x)g'(y) \ln\biggl\llvert
\frac{\tilde t_n(x) -\overline{\tilde
t_n(y)}}{\tilde t_n(x) - \tilde t_n(y)}\biggr\rrvert\,dx \,dy,
\end{eqnarray*}
which is the desired result.
We now establish
%
\begin{eqnarray}
&& -\frac{\kappa}{4\pi^2} \lim_{\varepsilon\downarrow0} \sum
_{\pm
_1,\pm_2} (\pm_1\pm_2) \int f(x) g(y)
\Theta_{2,n}(x\pm_1\mathbf{i} \varepsilon,y
\pm_2 \mathbf{i}\varepsilon) \,dx \,dy
\nonumber
\\[-8pt]
\label{eqcov2}
\\[-8pt]
\nonumber
&&\qquad =\frac{\kappa}{\pi^2 n} \sum_{i=1}^N \biggl(
\int_{{\mathcal S}_n} f'(x) \operatorname{Im} \bigl( x
T_n(x) \bigr) _{ii} \,dx \biggr) \biggl( \int
_{{\mathcal S}_n} g'(y) \operatorname{Im} \bigl(y
T_n(y) \bigr) _{ii} \,dy \biggr).
\end{eqnarray}
Due to formula \eqref{alternative-cov2}, we only need to prove
\begin{eqnarray}
&&
\frac{\mathbf{i}}{2\pi}\lim_{\varepsilon
\downarrow0}\sum
_{\tiny
\pm\in
\{ +,-\}} \pm\int f(x) \frac{\partial}{\partial x} \bigl[ (x\pm\mathbf{i}
\varepsilon)T_n(x\pm\mathbf{i}\varepsilon) \bigr]_{ii} \,dx
\nonumber
\\[-8pt]
\label{eqcov2-partial}
\\[-8pt]
\nonumber
&&\qquad= \frac{1}\pi\int_{{\mathcal S}_n} f'(x)
\operatorname{Im} \bigl(x T_n(x) \bigr) _{ii} \,dx.
\nonumber
\end{eqnarray}
Performing an integration by parts and taking into account the fact that
$T_n(\bar{z})=\overline{T_n(z)}$ yields
\begin{eqnarray*}
&&
\frac{\mathbf{i}}{2\pi}\lim_{\varepsilon\downarrow
0}\sum
_{\tiny\pm\in\{ +,-\}} \pm\int f(x) \frac{\partial}{\partial x} \bigl[
(x\pm\mathbf{i}
\varepsilon) T_n(x\pm\mathbf{i}\varepsilon) \bigr]_{ii} \,dx
\\
&&\qquad= - \frac{\mathbf{i}}{2\pi}\lim_{\varepsilon\downarrow0} \int f'(x)
2\mathbf{i} \operatorname{Im} \bigl[ (x+ \mathbf{i}\varepsilon) T_n(x+
\mathbf{i}\varepsilon) \bigr]_{ii} \,dx
\\
&&\qquad\stackrel{(a)}= \frac{1}\pi\int_{{\mathcal S}_n}
f'(x) \operatorname{Im} \bigl( x T_n(x)
\bigr)_{ii} \,dx,
\end{eqnarray*}
where step $(a)$ follows from the fact that
%
\begin{equation}
\label{lower-unif-bounded}
\mathop{\inf_{1\le i\le N,}}_{z\in(0,A]\times(0,B]}
\bigl\llvert\bigl(1+\tilde
t_n(z) \lambda_i\bigr)\bigr\rrvert>0,
\end{equation}
where the $\lambda_i$'s stand for $R_n$'s eigenvalues.
In fact, assume that \eqref{lower-unif-bounded} holds true, then using
the spectral decomposition of $R_n$,
the pointwise convergence of $\tilde t_n(z)$ to $\tilde t_n(x)$ as
$\mathbb{C}
^+ \ni z\to x\in\mathbb{R}$ (see, e.g., \cite
{silverstein-choi-1995}) and formula \eqref{eqdef-t-tilde}, then one
obtains the pointwise convergence
\[
\operatorname{Im} \bigl[ (x+ \mathbf{i}\varepsilon) T_n(x+ \mathbf{i}
\varepsilon) \bigr]_{ii}\mathop{\xrightarrow}_{\varepsilon\to0}
\operatorname{Im} \bigl[ x T_n(x) \bigr]_{ii}
\]
for $x>0$. Since $\operatorname{Im}(\tilde t(x))=0$ outside ${\mathcal
S}_n$, so is
$ \operatorname{Im} [ x T_n(x) ]_{ii}$. Finally, \eqref
{lower-unif-bounded}
provides a uniform bound for $\operatorname{Im} [ (x+ \mathbf
{i}\varepsilon) T_n(x+
\mathbf{i}\varepsilon) ]_{ii}$ and $(a)$ follows from the dominated
convergence theorem. It remains to prove \eqref{lower-unif-bounded}.
Assume that the infimum is zero, then there exists $\lambda^*\in\{
\lambda_1,\ldots,\lambda_N\}$
with $\lambda^*\neq0$ and a sequence $(z_\ell)$ such that $\tilde
t_n(z_{\ell}) \to- \frac{1}{\lambda^*}$ and $z_{\ell}\to x^*\in
\mathbb{R}$.
Formula \eqref{eqdef-t-tilde} yields
\begin{eqnarray*}
&&\forall z\in\mathbb{C}^+,\qquad \tilde t_n(z) =\frac{1}{ - z +\frac{1}n
\sum_{i=1}^N \frac{\lambda_i}{1+\tilde t_n(z) \lambda_i}}
\\
&&\qquad\Leftrightarrow\quad\frac{1}n \sum_{i=1}^N
\frac{\lambda_i}{1+\tilde t_n(z) \lambda_i} =\frac{1}{\tilde t_n(z)} +
z.
\end{eqnarray*}
Taking $z=z_\ell$ yields a contradiction since the LHS goes to
infinity while the RHS remains bounded. Necessarily,
\eqref{lower-unif-bounded} holds true and \eqref{eqcov2} is proved.

The proof of Proposition~\ref{propnew-variance-formula} is complete by
gathering \eqref{eqcov0}, \eqref{eqcov2} and using the fact that
$\Theta_{0,n}=\Theta_{1,n}$.

\subsection{Proof of Corollary~\texorpdfstring{\protect\ref{corounbounded-f}}{4.3}}\label
{secproof-coro-unbounded-f}
In order to establish the fluctuations in the case where functions
$(f_\ell)$ are $C^3$ in a neighborhood of ${\mathcal S}_{\infty}$
but may not have a bounded support, we proceed as following: Write
\begin{eqnarray*}
&&\operatorname{tr}f_\ell\bigl(\Sigma_n
\Sigma_n^*\bigr) - \mathbb{E} \operatorname{tr}(f_\ell
\mathbf{h}) \bigl(\Sigma_n \Sigma_n^*\bigr)
\\
&&\qquad= \underbrace{\operatorname{tr}f_\ell\bigl(\Sigma_n
\Sigma_n^*\bigr) - \operatorname{tr}(f_\ell\mathbf{h})
\bigl(\Sigma_n \Sigma_n^*\bigr)}_{\Gamma^1_\ell} +
\underbrace{\operatorname{tr}(f_\ell\mathbf{h}) \bigl(
\Sigma_n \Sigma_n^*\bigr) - \mathbb{E}
\operatorname{tr}(f_\ell\mathbf{h}) \bigl(\Sigma_n
\Sigma_n^*\bigr)}_{\Gamma_\ell^2}.
\end{eqnarray*}
By Proposition~\ref{propreduction}, the vector $(\Gamma_\ell^1)$
almost surely converges to zero while the fluctuations for vector
$(\Gamma_\ell^2)$
are described by Theorem~\ref{thnon-analytic} with covariance given
by Proposition~\ref{propnew-variance-formula}, where functions $f_k$
and $f_\ell$
must be replaced by $(f_k\mathbf{h})$ and $(f_\ell\mathbf{h})$. The
variance formula provided in this proposition shows that
$\operatorname{cov}(Z_ 1^n(f_{\ell} \mathbf{h}),Z_1^n(f_{\ell
'}\mathbf{h}))$ does not depend on function $\mathbf{h}$ as
long as $\mathbf{h}$ has value 1 on ${\mathcal S}_n$.

\section{Proof of Theorem~\texorpdfstring{\protect\ref{thnon-analytic-bias}}{3} (bias for
nonanalytic functionals)}\label{proofnon-analytic-II}
\subsection{Proof of Theorem~\texorpdfstring{\protect\ref{thnon-analytic-bias}}{3}}
Denote by $X^\mathbb{C}_n$ a $N\times n$ matrix whose entries are independent
standard complex circular Gaussian r.v. [i.e., $X^\mathbb
{C}_{ij}=U+iV$ where
$U,V$ are independent ${\mathcal N}(0,2^{-1})$ random variables];
denote accordingly $\Sigma^\mathbb{C}_n= n^{-1/2} R^{1/2} X^\mathbb
{C}_n$, $\xi^\mathbb{C}
_j = ( \Sigma^\mathbb{C}_n ){\cdot j}$ and
\[
Q^\mathbb{C}_n(z)= \bigl(-zI_N +
\Sigma^\mathbb{C}_n \bigl(\Sigma^\mathbb{C}_n
\bigr)^* \bigr)^{-1}.
\]
We split the bias into two terms:
%
\begin{eqnarray*}
&& \mathbb{E} \operatorname{Tr}f\bigl(\Sigma_n
\Sigma_n^*\bigr) - N\int f(\lambda) {\mathcal F}_n(d
\lambda)
\\
&&\qquad= \mathbb{E} \operatorname{Tr}f\bigl(\Sigma_n\Sigma_n^*
\bigr) - \mathbb{E} \operatorname{Tr} f\bigl(\Sigma^\mathbb{C}_n
\bigl(\Sigma_n^\mathbb{C}\bigr)^*\bigr)
\\
&& \qquad\quad{}+ \mathbb{E} \operatorname{Tr}f\bigl(\Sigma^\mathbb{C}_n
\bigl(\Sigma_n^\mathbb{C}\bigr)^*\bigr) - N\int f(\lambda) {
\mathcal F}_n(d\lambda)
\\
&&\qquad\stackrel{\triangle} = T_1 +T_2.
\end{eqnarray*}
%
We will prove the following. Provided that function $f$ is of class
$C^8$ with bounded support, then
%
\begin{eqnarray}
&& \mathbb{E} \operatorname{Tr}f\bigl(\Sigma_n
\Sigma_n^*\bigr) - \mathbb{E} \operatorname{Tr} f\bigl(
\Sigma^\mathbb{C}_n\bigl(\Sigma_n^\mathbb{C}
\bigr)^*\bigr) - \frac{1}\pi\operatorname{Re} \int_{\mathbb{C}^+}
\overline{\partial}\Phi_7(f) (z) {\mathcal B}_n(z)
\ell_2(dz)
\nonumber
\\[-8pt]
\label{eqbias-interpolation}
\\[-8pt]
\nonumber
&&\qquad \mathop{\xrightarrow}_{N,n\to\infty} 0.
\end{eqnarray}
Provided that function $f$ is of class $C^{18}$ with bounded support, then
%
\begin{equation}
\label{eqbias-hard} \mathbb{E} \operatorname{Tr}f\bigl(\Sigma^\mathbb{C}_n
\bigl(\Sigma_n^\mathbb{C}\bigr)^*\bigr) - N\int f(\lambda) {
\mathcal F}_n(d\lambda) \mathop{\xrightarrow}_{N,n\to\infty} 0.
\end{equation}
As one can check, it is much more demanding in terms of assumptions to
prove \eqref{eqbias-hard} than \eqref{eqbias-interpolation}.
Convergence in \eqref{eqbias-hard} should be compared to the results
in Haagerup and Th\o rbjornsen \cite{haagerup-thorbjornsen-2005}
(counterpart in the GUE case), Schultz \cite{schultz-ptrf-2005} (GOE),
Capitaine and Donati-Martin \cite{capitaine-donati-freeness-2007},
Loubaton et al. \cite{vallet-2012} (``signal plus noise'' model), etc.

\subsubsection{Proof of \texorpdfstring{\protect\eqref{eqbias-interpolation}}{(7.1)}}\label
{secMM-bias} The heart of the proof lies in Helffer--Sj\"{o}strand's
formula, in Theorem~\ref{lemmamain} (bias part) and in a dominated
convergence argument.
By Theorem~\ref{lemmamain},
\[
\mathbb{E} \operatorname{Tr} \bigl( \Sigma_n\Sigma_n^*
-zI_N \bigr)^{-1} - N t_n(z) - {\mathcal
B}_n(z) \mathop{\xrightarrow}_{N,n\to\infty} 0.
\]
The same argument yields
\[
\mathbb{E} \operatorname{Tr} \bigl( \Sigma^\mathbb{C}_n\bigl(
\Sigma^\mathbb{C}_n\bigr)^* - zI_N
\bigr)^{-1} - N t_n(z) \mathop{\xrightarrow}_{N,n\to\infty}
0,
\]
because in the later case ${\mathcal V}=\kappa=0$, hence the bias is
zero for the matrix model $\Sigma_n^\mathbb{C}(\Sigma^\mathbb{C}_n)^*$.
Subtracting yields
\[
\mathbb{E} \operatorname{Tr}Q_n(z) - \mathbb{E} \operatorname
{Tr}Q_n^\mathbb{C}(z) - {\mathcal B}_n(z) \mathop{
\xrightarrow}_{N,n\to\infty} 0.
\]
The following proposition will be of help.
%
\begin{prop}\label{propapprox-gaussienne} Assume that
Assumptions \ref{assX} and \ref{assR} hold true, then
%
\begin{equation}
\label{estimateMM}
\bigl\llvert\mathbb{E}\operatorname{Tr}Q(z) -
\mathbb{E}
\operatorname{Tr}Q^\mathbb{C}(z) \bigr\rrvert\le K\frac
{|z|^3}{\operatorname{Im}(z)^7},
\end{equation}
where $K$ is independent from $N,n,z$.
\end{prop}
The proof is based on classical rank-one perturbation arguments and is
omitted (details can be found in Section~\ref
{secdiagonal-terms-variance} of the previous version
of this article---see footnote \ref{footpreprintv3}).

In order to transfer this bound to ${\mathcal B}_n(z)$, we invoke a
meta-model argument (cf. Section~\ref{secmeta-model}): Consider
matrix $\Sigma_n(M)$ and its counterpart $\Sigma_n^\mathbb{C}(M)$ as defined
in~\eqref{eqextended-model} and recall that in this case, we have a
genuine limit:
\begin{eqnarray*}
&& \mathbb{E} \operatorname{Tr} \bigl( \Sigma_n(M)
\Sigma_n^*(M) -zI_{NM} \bigr)^{-1} - \mathbb{E}
\operatorname{Tr} \bigl( \Sigma^\mathbb{C}_n(M) \bigl(
\Sigma_n^\mathbb{C} (M)\bigr)^*-zI_{NM}
\bigr)^{-1}
\\
&& \qquad\mathop{\mathop{\xrightarrow}_{M\to\infty}}_{N,n \ \mathrm{fixed}}
{\mathcal
B}_n(z).
\end{eqnarray*}
Since the estimate \eqref{estimateMM} remains true for all $M\ge1$,
we obtain
%
\begin{eqnarray}
\nonumber
\bigl\llvert{\mathcal B}_n(z)\bigr\rrvert&=&\lim
_{M\to\infty} \big| \mathbb{E} \operatorname{Tr} \bigl(
\Sigma_n(M)\Sigma_n(M)^* -zI_{MN}
\bigr)^{-1}
\\
\label{estimateBn}
&& {}- \mathbb{E} \operatorname{Tr} \bigl( \Sigma^\mathbb{C}(M)_n
\bigl(\Sigma_n^\mathbb{C}(M)\bigr)^*-zI_{NM}
\bigr)^{-1} \big|
\\
&\le& K\frac{|z|^3}{\operatorname{Im}(z)^7}.
\nonumber
\end{eqnarray}

Write
%
\begin{eqnarray}
 && \mathbb{E} \operatorname{Tr}f\bigl(\Sigma_n
\Sigma_n^*\bigr) - \mathbb{E} \operatorname{Tr} f\bigl(
\Sigma^\mathbb{C}_n\bigl(\Sigma_n^\mathbb{C}
\bigr)^*\bigr) -\frac{1}\pi\operatorname{Re} \int_{\mathbb{C}+}
\overline{\partial} \Phi(f) (z) {\mathcal B}_n(z)
\ell_2(dz)
\nonumber
\\[-8pt]
\label{eqDCT-argument}
\\[-8pt]
\nonumber
&&\qquad= \frac{1}\pi\operatorname{Re} \int_{\mathbb{C}^+} \overline{
\partial} \Phi(f) (z) \bigl\{ \mathbb{E}\operatorname{Tr}Q_n(z) -
\mathbb{E}\operatorname{Tr}Q_n^\mathbb{C}(z) - {\mathcal
B}_n(z) \bigr\} \ell_2(dz). 
\end{eqnarray}
In view of \eqref{eqDCT-argument}, we need a dominated convergence
argument in order to prove~\eqref{eqbias-interpolation}; such an
argument follows from
Proposition~\ref{propcompensation}, \eqref{estimateMM} and \eqref
{estimateBn} as long as $f$ is of class $C^8$ with large but bounded support.
This completes the proof of \eqref{eqbias-interpolation}.

\subsubsection{Proof of \texorpdfstring{\protect\eqref{eqbias-hard}}{(7.2)}} The gist of the proof
lies in the following proposition.
%
\begin{prop}\label{propbias-hard} Denote by $P_\ell(X)$ a polynomial
in $X$ with degree $\ell$ and positive coefficients, then
\[
\bigl\llvert\mathbb{E}\operatorname{Tr} \bigl( \Sigma_n^\mathbb{C}
\bigl( \Sigma_n^\mathbb{C} \bigr)^* - zI_N
\bigr)^{-1} - Nt_n(z)\bigr\rrvert\le\frac{1}n
P_{12}\bigl(|z|\bigr)P_{17}\bigl(\bigl|\operatorname{Im}(z)\bigr|^{-1}
\bigr).
\]
\end{prop}
The proof of Proposition~\ref{propbias-hard} builds upon techniques
borrowed from \cite{haagerup-thorbjornsen-2005,vallet-2012} and is
omitted. Details can be found in the previous version of this article;
see footnote \ref{footpreprintv3} and \cite{yao-PhD-2013}.

Using Helffer--Sj\"{o}strand's formula, Proposition~\ref
{propbias-hard} together with Proposition~\ref{propcompensation}
immediately yield \eqref{eqbias-hard} for any $f$ of class $C^{18}$
with large but bounded support.

\subsection{Proof of Proposition~\texorpdfstring{\protect\ref{propbias-as-distribution}}{4.4}}\label{proofprop-bias-distribution}
One can easily prove that $Z_n^2$ is a distribution on
$C_c^{18}(\mathbb{R})$
following the lines of proof of Proposition~\ref
{propidentification-distribution}.
Similarly, one can establish the boundary value representation \eqref
{eqboundary-bias}. It remains to prove that the singular points of
${\mathcal B}_n(z)$ are included in ${\mathcal S}_n$. Following the
definitions of ${\mathcal B}_{1,n}$ and ${\mathcal B}_{2,n}$ [cf.
\eqref{defB1} and \eqref{defB2}], we simply need to prove that the
quantities
\[
\biggl( 1 - z^2 \tilde t_n^2
\frac{1}n \operatorname{Tr}R_n^2
T_n^2 \biggr)\quad \mbox{and}\quad \biggl( 1 - |
\mathcal{V}|^2 z^2 \tilde t_n^2
\frac{1}n \operatorname{Tr}R_n^{1/2}
T_n(z) R_n^{1/2} \overline{R}_n^{1/2}
T^T_n(z) \overline{R}_n^{1/2} \biggr)
\]
are invertible for $z\notin{\mathcal S}_n$. We focus on the first one.
Assume first that $z\in\mathbb{C}\setminus\mathbb{R}$. Using the
inequality $|\operatorname{tr}
(AB)| \le(\operatorname{tr}(AA^*) \operatorname{tr}(BB^*))^{1/2}$ yields
\begin{eqnarray*}
&&\biggl\llvert z^2 \tilde t_n^2(z)
\frac{1}n \operatorname{Tr}R_n^2
T_n^2(z)\biggr\rrvert\le\frac{|z|^2 |\tilde t_n(z)|^2}n
\operatorname{tr}R_n T_n(z) R_n
T_n^*(z).
\end{eqnarray*}
Since $T_n^*(z)=T_n(\bar{z})$, we can assume without loss of
generality that $z_1,z_2\in\mathbb{C}^+$:
%
\begin{eqnarray}
\biggl\llvert1- z^2 \tilde t_n^2(z)
\frac{1}n \operatorname{Tr}R_n^2
T_n^2(z) \biggr\rrvert&\ge& 1- \frac{|z|^2 |\tilde t_n(z)|^2}n
\operatorname{tr}R_n T_n(z) R_n
T_n^*(z)
\nonumber
\\[-8pt]
\label{estimate-denominator}
\\[-8pt]
\nonumber
&=& \bigl|\tilde t_n(z)\bigr|^2 \frac{\operatorname{Im}(z)}{\operatorname
{Im}(\tilde
t_n(z))},
\end{eqnarray}
where the last identity can be found in the previous version of this
article [equation~(A.15)]; see footnote \ref{footpreprintv3}. In
order to extend the previous estimate to $z\in\mathbb{R}\setminus
{\mathcal
S}_n$, let
$z=x+\mathbf{i}y$ with $x\in\mathbb{R}\setminus{\mathcal S}_n$;
then a direct
computation yields
\[
\frac{\operatorname{Im}(\tilde t_n(z))}{\operatorname{Im}(z)} = \int
\frac
{{\tilde
F}_n(d\lambda)}{|\lambda- z|^2}\mathop{\xrightarrow}_{y\searrow0}
\int\frac{{\widetilde F}_n(d\lambda)}{|\lambda- x|^2}\neq0.
\]
Therefore, by continuity $(z)\mapsto1-z^2 \tilde t_n^2(z) \frac{1}n
\operatorname{Tr}R_n^2 T_n^2(z) $ does not vanish on $\mathbb
{C}\setminus{\mathcal
S}_n$ and ${\mathcal B}_{1,n}$
is analytic on this set. We can similarly prove that ${\mathcal
B}_{2,n}$ is also analytic on the same set. Consider now a function
$f\in C^{18}_c(\mathbb{R})$ whose support
is disjoint from ${\mathcal S}_n$, then it is straightforward to check
that $Z_n^2(f)=0$ and the proof of the proposition is completed.
%
%
%

\begin{appendix}
\section*{Appendix: Remaining proofs}

\subsection{Proof of Lemma~\texorpdfstring{\protect\ref{lemmaabstract}}{6.3}}\label
{prooflemmaabstract} By Proposition~\ref{propcompensation},
%
%
\begin{eqnarray*}
\mathbb{E}\biggl\llvert\int_D \overline\partial\Phi(g)
(z)\varphi_n(z) \ell_2(dz) \biggr\rrvert&\le& \int
_D \bigl|\overline{\partial} \Phi(g) (z)\bigr| \mathbb{E}\bigl|\varphi
_n(z)\bigr|\ell_2(dz)
\\
&\le& \| g\|_{k+1,\infty} \int_D \operatorname{Im}(z)^k
\bigl\{ \operatorname{var} \varphi_n(z) \bigr\}^{1/2}
\ell_2(dz)
\\
&<&\infty,
\end{eqnarray*}
by (iii) and (iv). Hence, $\frac{1}\pi\operatorname{Re}\int_D \overline
\partial\Phi(g)(z)\varphi_n(z)\ell_2(dz)$ is a well-defined a.s.
finite random variable. This estimate, uniform in $n$, readily implies
the tightness of
\[
\biggl( \frac{1}\pi\operatorname{Re}\int_D \overline
\partial\Phi(g) (z)\varphi_n(z)\ell_2(dz); n\in
\mathbb{N} \biggr).
\]
Notice that the integrals with $\psi_n$ instead of $\varphi_n$ are
similarly well defined and tight.

By conditions (i) and (ii), we obtain
\setcounter{equation}{0}
\begin{eqnarray}
&& d_{\mathcal{LP}} \biggl( \frac{1}\pi\operatorname
{Re}\int
_{D_\varepsilon} \overline\partial\Phi({\mathbf g}) (z)
\varphi_n(z)\ell_2(dz), \frac{1}\pi\operatorname{Re}
\int_{D_\varepsilon} \overline\partial\Phi({\mathbf g}) (z)
\psi_n(z)\ell_2(dz) \biggr)
\nonumber
\\[-8pt]
\label{eqconv-T}
\\[-8pt]
\nonumber
&&\qquad \mathop{\xrightarrow}_{N,n\to\infty} 0
\nonumber
\end{eqnarray}
(apply Lemma~\ref{lemmatightness}).

Let ${\mathbf g}=(g_\ell; 1\le\ell\le L)$ and $f:\mathbb{C}^L\to
\mathbb{C}$ be
bounded and continuous. Consider the following notation:
\begin{eqnarray*}
\xi_n &=& \frac{1}\pi\operatorname{Re}\int_{D}
\overline\partial\Phi({\mathbf g}) (z)\varphi_n(z)
\ell_2(dz), \qquad\xi_n^\varepsilon=\frac{1}\pi
\operatorname{Re}\int_{D_\varepsilon} \overline\partial\Phi({\mathbf
g}) (z)
\varphi_n(z)\ell_2(dz),
\\
\eta_n &=&\frac{1}\pi\operatorname{Re}\int_{D}
\overline\partial\Phi({\mathbf g}) (z)\psi_n(z)\ell_2(dz)
, \qquad\eta_n^\varepsilon=\frac{1}\pi\operatorname{Re}\int
_{D_\varepsilon
} \overline\partial\Phi({\mathbf g}) (z)
\psi_n(z)\ell_2(dz).
\end{eqnarray*}
We have
%
\begin{eqnarray}
 &&\qquad \bigl\llvert\mathbb{E}f ( \xi_n ) - \mathbb{E}f (
\eta_n ) \bigr\rrvert
\nonumber
\\[-8pt]
\label{eqcontrole-0}
\\[-8pt]
\nonumber
&&\qquad\qquad\le\bigl\llvert\mathbb{E}f (\xi_n ) - \mathbb{E}f \bigl( \xi
_n^\varepsilon\bigr) \bigr\rrvert+ \bigl\llvert\mathbb{E}f
\bigl(\xi_n^\varepsilon\bigr) - \mathbb{E}f \bigl( \eta
_n^\varepsilon\bigr)\bigr\rrvert+\bigl\llvert\mathbb{E}f \bigl(
\eta_n^\varepsilon\bigr) - \mathbb{E}f ( \eta_n )
\bigr\rrvert.
\end{eqnarray}
Given $\rho>0$, we first prove that for all $n\ge1$,
%
\begin{equation}
\label{eqcontrole-1}
\bigl\llvert\mathbb{E}f ( \xi_n ) - \mathbb{E}f
\bigl( \xi_n^\varepsilon\bigr) \bigr\rrvert\le\bigl(4\|f
\|_{\infty} +1\bigr)\rho
\end{equation}
for $\varepsilon$ small enough.

We have
%
\begin{equation}
\label{eqconv-proba}
\mathbb{P}\bigl\{ \bigl|\xi_n - \xi_n^\varepsilon\bigr|>
\delta\bigr\} \le\frac{1}\delta\biggl( \int_{[0,A]+\mathbf
{i}[0,\varepsilon]} \bigl|
\overline\partial\Phi({\mathbf g}) (z)\bigr| \mathbb{E}\bigl|\varphi_n(z)\bigr|
\ell_2(dz) \biggr)
\end{equation}
which can be made arbitrarily small if $\varepsilon$ is small enough,
independently from $n$. Now,
\begin{eqnarray*}
\bigl\llvert\mathbb{E}f ( \xi_n ) - \mathbb{E}f \bigl( \xi
_n^\varepsilon\bigr) \bigr\rrvert& \le& \bigl\llvert\mathbb{E}f
( \xi_n ) - \mathbb{E}f \bigl( \xi_n^\varepsilon
\bigr) \bigr\rrvert1_{\{ |\xi_n - \xi_n^\varepsilon|>\eta\} }
\\
&& {}+ \bigl\llvert\mathbb{E}f ( \xi_n ) - \mathbb{E}f \bigl( \xi
_n^\varepsilon\bigr) \bigr\rrvert1_{\{ |\xi_n - \xi_n^\varepsilon|\le
\eta, |\xi
_n|\vee|\xi_n^\varepsilon| >K\} }
\\
&& {}+ \bigl\llvert\mathbb{E}f ( \xi_n ) - \mathbb{E}f \bigl( \xi
_n^\varepsilon\bigr) \bigr\rrvert1_{\{ |\xi_n - \xi_n^\varepsilon|\le
\eta, |\xi
_n|\vee|\xi_n^\varepsilon| \le K\} }.
\end{eqnarray*}
First, invoke the tightness of $|\xi_n|\vee|\xi_n^\varepsilon|$ and
choose $K$ large enough so that the second term of the RHS is lower
than $2\|f\|_{\infty}\rho$; then choose $\eta>0$ small enough so
that $f$ being absolutely continuous over $\{ z\in\mathbb{C}^+,|z|\le
K\}$,
the third term of the RHS is lower that $\rho$; finally for such $K$
and $\eta$, take advantage of \eqref{eqconv-proba} and choose
$\varepsilon$ small enough so that the first term of the RHS is lower
than $2\|f\|_{\infty}\rho$. Equation~\eqref{eqcontrole-1} is proved.

One can similarly prove that $\llvert
\mathbb{E}f ( \eta_n ) - \mathbb{E}f ( \eta
_n^\varepsilon)
\rrvert\le(4\|f\|_{\infty} +1)\rho$ for $\varepsilon>0$ small
enough. Such $\varepsilon$ being fixed, it remains to control the
second term of the RHS of \eqref{eqcontrole-0}, but this immediately
follows from \eqref{eqconv-T}.

In order to prove that $\eta_n$ is multivariate Gaussian with
prescribed covariance~\eqref{eqcov-non-analytic}, we first consider
$\eta_n^\varepsilon$. Approximating the integral in $\eta
_n^\varepsilon$ by Riemann sums and using the fact that weak limits of
Gaussian vectors are Gaussian immediately yields that $\eta
_n^\varepsilon$ is a Gaussian vector with covariance matrix
\begin{eqnarray*}
&& \bigl[ \operatorname{cov}\bigl(\eta_n^\varepsilon\bigr)
\bigr]_{k\ell} \\
&&\qquad= \frac{1}{\pi^2} \mathbb{E} \biggl\{ \operatorname{Re}
\int
_{D^\varepsilon} \overline{\partial} \Phi(g_k) (z)
\psi_n(z) \ell_2(dz) \operatorname{Re} \int_{D^\varepsilon}
\overline{\partial} \Phi(g_\ell) (z) \psi_n(z)
\ell_2(dz) \biggr\}
\end{eqnarray*}
for $1\le k,\ell\le L$. Using the elementary identity
\[
\operatorname{Re}(z) \operatorname{Re}\bigl(z'\bigr)= \frac{
\operatorname{Re}(z\overline{z'})
+ \operatorname{Re}(zz')}2,
\]
we obtain
\begin{eqnarray*}
&&
\bigl[ \operatorname{cov}\bigl(\eta_n^\varepsilon\bigr)
\bigr]_{k\ell
}
\\
&&\qquad = \frac{1}{2\pi^2} \operatorname{Re} \int_{(D^\varepsilon)^2}
\overline{
\partial} \Phi(g_k) (z_1) \overline{\overline{\partial}
\Phi( g_\ell) (z_2)} \mathbb{E}\psi_n(z_1)
\overline{\psi_n(z_2)} \ell_2(
dz_1)\ell_2( dz_2)
\\
&&\qquad\quad{}+ \frac{1}{2\pi^2} \operatorname{Re} \int_{(D^\varepsilon)^2}
\overline{
\partial} \Phi(g_k) (z_1) {\overline{\partial}
\Phi(g_\ell) (z_2)} \mathbb{E}\psi_n(z_1){
\psi_n(z_2)} \ell_2( dz_1)
\ell_2( dz_2).
\end{eqnarray*}
Using the fact that $\overline{\psi_n(z_2)} = \psi_n(\overline
{z_2})$ yields
\begin{eqnarray*}
&&\bigl[ \operatorname{cov}\bigl(\eta_n^\varepsilon\bigr)
\bigr]_{k\ell
}
\\
&&\qquad= \frac{1}{2\pi^2} \operatorname{Re} \int_{(D^\varepsilon)^2}
\overline{
\partial} \Phi( g_k) (z_1) \overline{\overline{\partial}
\Phi( g_\ell) (z_2)} \kappa_n(z_1,
\overline{z_2}) \ell_2( dz_1)
\ell_2( dz_2)
\\
&&\qquad\quad{}+ \frac{1}{2\pi^2} \operatorname{Re} \int_{(D^\varepsilon)^2}
\overline{
\partial} \Phi( g_k) (z_1) {\overline{\partial}
\Phi(g_\ell) (z_2)} \kappa_n(z_1,z_2)
\ell_2( dz_1)\ell_2( dz_2).
\end{eqnarray*}
In order to lift the Gaussianity from $\eta_n^\varepsilon$ to $\eta
_n$ and to extend the covariance formula from the one above to formula
\eqref{eqcov-non-analytic}, we rely on the approximation theorem
\cite{book-kallenberg-second-edition}, Theorem~4.28, and on assumptions
(iv) and (v) on the variance estimates and on the regularity of
functions $g_k,g_\ell$ in Lemma~\ref{lemmaabstract}.

The proof of Lemma~\ref{lemmaabstract} is complete.

\subsection{Proof of Proposition~\texorpdfstring{\protect\ref{propvariance-estimate-sharp}}{6.4}(ii)}\label{appprop-variance-sharp}

We rely on a meta-model argument (cf. Section~\ref{secmeta-model}).
Denote by
\[
M_{n,M}^1(z)= \operatorname{Tr} \bigl(
\Sigma_n(M) \Sigma_n^*(M) - zI_N
\bigr)^{-1} - \mathbb{E}\operatorname{Tr} \bigl( \Sigma_n(M)
\Sigma_n^*(M) - zI_N \bigr)^{-1},
\]
then by Proposition~\ref{propvariance-estimate-sharp}(i), we get
\[
\operatorname{var} \bigl\{ \operatorname{tr} \bigl( \Sigma_n(M) \Sigma
_n^*(M) - zI_N \bigr)^{-1} \bigr\} \le
\frac{C}{\operatorname{Im}(z)^4},
\]
moreover $M_{n,M}^1(z)$ converges in distribution to $\psi_n(z)$ as
$M\to\infty$, $N$ and $n$ being fixed (see, e.g., the details
in Section~\ref{secgaussian-process}). Consider the continuous
bounded function $h_K(x)=|x|^2 \wedge K$, then
\[
\mathbb{E}h_K\bigl(\psi_n(z)\bigr) =\lim
_{M\to\infty} \mathbb{E}h_K\bigl(M_{n,M}^1(z)
\bigr) \le\limsup_{M\to\infty} \mathbb{E}\bigl|M_{n,M}^1(z)\bigr|^2
\le\frac{C}{\operatorname{Im}(z)^4}.
\]
Now letting $K\to\infty$ yields the desired bound by monotone
convergence theorem.
\end{appendix}

\section*{Acknowledgements}
We are particularly indebted to Charles Bordenave who drew our
attention to Helffer--Sj\"{o}strand's formula and related variance
estimates, which substantially shorten the initial proof of
fluctuations for nonanalytic functions; we thank Reinhold Meise for
his help to understand Tillmann's article; we also thank Djalil Chafa\"{\i},
Walid Hachem and Philippe Loubaton for fruitful discussions;
finally, we thank the referee whose remarks have substantially
simplified the proof of Proposition~\ref{propcumulant}.

%
%

%

\printaddresses
\end{document}